\documentclass[12pt]{article}
\usepackage{amsfonts,amssymb,bm,subeqnarray,eqsection,indent,epsf,cite}
\usepackage{graphicx}

\begin{document}

\title{\vspace*{-2.5cm} How to Clean a Dirty Floor: \\[2mm]
       Probabilistic Potential Theory \\
       and the Dobrushin Uniqueness Theorem}

\author{
  {\small Thierry de la Rue}                  \\[-1.7mm]
  {\small Roberto Fern\'andez}                  \\[-1.7mm]
  {\small\it Laboratoire de Math\'ematiques Rapha\"el Salem}       \\[-1.7mm]
  {\small\it UMR 6085 CNRS}         \\[-1.7mm]
  {\small\it Universit\'e de Rouen}          \\[-1.7mm]
  {\small\it Avenue de l'Universit\'e}          \\[-1.7mm]
  {\small\it F-76801 Saint Etienne du Rouvray, FRANCE}      \\[-1.7mm]
  {\small\tt thierry.de-la-rue@univ-rouen.fr},
  {\small\tt Roberto.Fernandez@univ-rouen.fr}               \\[-1.7mm]
  {\protect\makebox[5in]{\quad}}  
  \\[-2mm]
  {\small Alan D.~Sokal\thanks{Also at Department of Mathematics,
           University College London, London WC1E 6BT, England.}}  \\[-1.7mm]
  {\small\it Department of Physics}       \\[-1.7mm]
  {\small\it New York University}         \\[-1.7mm]
  {\small\it 4 Washington Place}          \\[-1.7mm]
  {\small\it New York, NY 10003 USA}      \\[-1.7mm]
  {\small\tt sokal@nyu.edu}               \\[-1.7mm]
  {\protect\makebox[5in]{\quad}}  
}

\date{\vspace*{-5mm} April 23, 2007}

\maketitle
\thispagestyle{empty}   

\vspace*{-6mm}

\begin{abstract}
Motivated by the Dobrushin uniqueness theorem in statistical mechanics,
we consider the following situation:
Let $\alpha$ be a nonnegative matrix over a finite or countably infinite
index set $X$, and define the ``cleaning operators''
$\beta_h = I_{1-h} + I_h \alpha$ for $h \colon\, X \to [0,1]$
(here $I_f$ denotes the diagonal matrix with entries $f$).
We ask: For which ``cleaning sequences'' $h_1,h_2,\ldots$
do we have $c \beta_{h_1} \cdots \beta_{h_n} \to 0$
for a suitable class of ``dirt vectors'' $c$?
We show, under a modest condition on $\alpha$,
that this occurs whenever $\sum_i h_i = \infty$ everywhere on $X$.
More generally, we analyze the cleaning of subsets $\Lambda \subseteq X$
and the final distribution of dirt on the complement of $\Lambda$.
We show that when ${\rm supp}(h_i) \subseteq \Lambda$ with
$\sum_i h_i = \infty$ everywhere on $\Lambda$,
the operators $\beta_{h_1} \cdots \beta_{h_n}$ converge as $n \to \infty$
to the ``balayage operator''
$\Pi_\Lambda = \sum_{k=0}^\infty (I_\Lambda \alpha)^k I_{\Lambda^c}$.
These results are obtained in two ways:
by a fairly simple matrix formalism,
and by a more powerful tree formalism that corresponds to working with
formal power series in which the matrix elements of $\alpha$
are treated as noncommuting indeterminates.
\end{abstract}


\bigskip
\noindent
{\bf Mathematics Subject Classification (MSC 2000) codes:}
60J99 (Primary);
15A48, 31C20, 31C99, 60J10, 60J45, 82B20 (Secondary).

\clearpage

\newtheorem{defin}{Definition}[section]
\newtheorem{definition}[defin]{Definition}
\newtheorem{prop}[defin]{Proposition}
\newtheorem{proposition}[defin]{Proposition}
\newtheorem{lem}[defin]{Lemma}
\newtheorem{lemma}[defin]{Lemma}
\newtheorem{guess}[defin]{Conjecture}
\newtheorem{ques}[defin]{Question}
\newtheorem{question}[defin]{Question}
\newtheorem{prob}[defin]{Problem}
\newtheorem{problem}[defin]{Problem}
\newtheorem{thm}[defin]{Theorem}
\newtheorem{theorem}[defin]{Theorem}
\newtheorem{cor}[defin]{Corollary}
\newtheorem{corollary}[defin]{Corollary}
\newtheorem{conj}[defin]{Conjecture}
\newtheorem{conjecture}[defin]{Conjecture}
\newtheorem{examp}[defin]{Example}
\newtheorem{example}[defin]{Example}
\newtheorem{claim}[defin]{Claim}



\renewcommand{\theenumi}{\alph{enumi}}
\renewcommand{\labelenumi}{(\theenumi)}
\def\prf{\par\noindent{\bf Proof.\enspace}\rm}
\def\rmk{\par\medskip\noindent{\bf Remark.\enspace}\rm}

\newcommand{\be}{\begin{equation}}
\newcommand{\ee}{\end{equation}}
\newcommand{\<}{\langle}
\renewcommand{\>}{\rangle}
\newcommand{\widebar}{\overline}
\def\reff#1{(\protect\ref{#1})}
\def\spose#1{\hbox to 0pt{#1\hss}}
\def\ltapprox{\mathrel{\spose{\lower 3pt\hbox{$\mathchar"218$}}
 \raise 2.0pt\hbox{$\mathchar"13C$}}}
\def\gtapprox{\mathrel{\spose{\lower 3pt\hbox{$\mathchar"218$}}
 \raise 2.0pt\hbox{$\mathchar"13E$}}}
\def\textprime{${}^\prime$}
\def\proof{\par\medskip\noindent{\sc Proof.\ }}
\newcommand{\qed}{\quad $\Box$ \medskip \medskip}
\def\proofof#1{\bigskip\noindent{\sc Proof of #1.\ }}
\def\half{ {1 \over 2} }
\def\third{ {1 \over 3} }
\def\twothird{ {2 \over 3} }
\def\smfrac#1#2{\textstyle{#1\over #2}}
\def\smhalf{ \smfrac{1}{2} }
\newcommand{\real}{\mathop{\rm Re}\nolimits}
\renewcommand{\Re}{\mathop{\rm Re}\nolimits}
\newcommand{\imag}{\mathop{\rm Im}\nolimits}
\renewcommand{\Im}{\mathop{\rm Im}\nolimits}
\newcommand{\sgn}{\mathop{\rm sgn}\nolimits}
\newcommand{\spr}{\mathop{\rm spr}\nolimits}
\newcommand\supp{\mathop{\rm supp}\nolimits}
\def\hboxscript#1{ {\hbox{\scriptsize\em #1}} }

\newcommand{\restrict}{\upharpoonright}
\renewcommand{\emptyset}{\varnothing}

\def\Z{{\mathbb Z}}
\def\ZZ{{\mathbb Z}}
\def\R{{\mathbb R}}
\def\C{{\mathbb C}}
\def\CC{{\mathbb C}}
\def\N{{\mathbb N}}
\def\NN{{\mathbb N}}
\def\Q{{\mathbb Q}}

\newcommand{\FH}{{\rm (FH)\ }}  

\newcommand{\scra}{{\mathcal{A}}}
\newcommand{\scrb}{{\mathcal{B}}}
\newcommand{\scrc}{{\mathcal{C}}}
\newcommand{\scrf}{{\mathcal{F}}}
\newcommand{\scrg}{{\mathcal{G}}}
\newcommand{\scrh}{{\mathcal{H}}}
\newcommand{\scrl}{{\mathcal{L}}}
\newcommand{\scro}{{\mathcal{O}}}
\newcommand{\scrp}{{\mathcal{P}}}
\newcommand{\scrr}{{\mathcal{R}}}
\newcommand{\scrs}{{\mathcal{S}}}
\newcommand{\scrt}{{\mathcal{T}}}
\newcommand{\scrv}{{\mathcal{V}}}
\newcommand{\scrw}{{\mathcal{W}}}
\newcommand{\scrz}{{\mathcal{Z}}}

\newcommand{\amat}{{\sf a}}

\newcommand{\bgamma}{\boldsymbol{\gamma}}
\newcommand{\bsigma}{\boldsymbol{\sigma}}
\renewcommand{\pmod}[1]{\;({\rm mod}\:#1)}


\newenvironment{sarray}{
	  \textfont0=\scriptfont0
	  \scriptfont0=\scriptscriptfont0
	  \textfont1=\scriptfont1
	  \scriptfont1=\scriptscriptfont1
	  \textfont2=\scriptfont2
	  \scriptfont2=\scriptscriptfont2
	  \textfont3=\scriptfont3
	  \scriptfont3=\scriptscriptfont3
	\renewcommand{\arraystretch}{0.7}
	\begin{array}{l}}{\end{array}}

\newenvironment{scarray}{
	  \textfont0=\scriptfont0
	  \scriptfont0=\scriptscriptfont0
	  \textfont1=\scriptfont1
	  \scriptfont1=\scriptscriptfont1
	  \textfont2=\scriptfont2
	  \scriptfont2=\scriptscriptfont2
	  \textfont3=\scriptfont3
	  \scriptfont3=\scriptscriptfont3
	\renewcommand{\arraystretch}{0.7}
	\begin{array}{c}}{\end{array}}

\newcommand{\bydef}{:=}
\newcommand{\defby}{=:}
\newcommand{\lexc}{\le^{w!}}
\newcommand{\norm}[1]{\left\|#1\right\|_w}
\newcommand{\nw}[1]{\left\|#1\right\|_{w}}
\newcommand{\triplenorm}{\| \hspace*{-0.3mm} |}
\newcommand{\cb}[1]{\succ_#1}
\newcommand{\cw}[1]{\prec_#1}
\newcommand{\ind}[1]{\chi_{#1}}
\newcommand{\lelambda}{\preccurlyeq_\Lambda}
\newcommand{\wsense}{\trianglelefteq}
\newcommand{\audessus}{\trianglelefteq}
\newcommand{\audessusreverse}{\trianglerighteq}
\newcommand{\ancestor}{\preccurlyeq}
\newcommand{\anc}{\preccurlyeq}
\newcommand{\ancneq}{\prec}
\newcommand{\ancreverse}{\succcurlyeq}
\newcommand{\prefix}{\preccurlyeq}
\newcommand{\leads}{\to}
\newcommand{\implies}{\Longrightarrow}
\newcommand{\doubleminus}{\setminus\!\!\!\setminus\,}
\newcommand{\binom}[2]{\left(#1\atop#2\right)}

\tableofcontents

\clearpage

\section{Introduction} \label{sec1}

Let $X$ be a finite or countably infinite set,
let $\scrt$ be a collection of nonnegative real matrices indexed by $X$,
and let $\scrc$ be a class of nonnegative real vectors indexed by $X$.
In this paper we want to ask variants of the following question:
Under what conditions does there exist a sequence of elements
$(T_i)_{i=1}^\infty$ in $\scrt$
such that $c T_1 \cdots T_n \to 0$ for all $c \in \scrc$?
(If $X$ is infinite, we must of course specify the topology
 in which this convergence is to be understood.)

Here is a homely but suggestive interpretation:
Think of the elements of $X$ as the ``sites'' of a dirty floor,
the nonnegative vectors $c = (c_x)_{x \in X}$ as ``distributions of dirt'',
and the matrices $T \in \scrt$ as the ``cleaning operators'' at our disposal.
Application of a cleaning operator $T$ transforms the
dirt distribution from $c$ to $cT$.
(Note that we always write our dirt vectors on the left,
 in analogy with probability distributions in Markov-chain theory.)
It is natural to ask:
Under what conditions can the floor be completely cleaned?
We will also ask:
Under what conditions can a subset $\Lambda \subseteq X$ be cleaned,
and in this case,
where in $\Lambda^c \bydef X \setminus \Lambda$ does the dirt go?

These questions arise in mathematical statistical mechanics
in connection with the Dobrushin
\cite{Dobrushin_68,Dobrushin_70,Vasershtein_69,Lanford_73,Follmer_82,%
Georgii_88,Simon_93}
and Dobrushin--Shlosman \cite{Dobrushin-Shlosman,Aizenman_unpub,Weitz_05}
uniqueness theorems.
Indeed, the simplest proofs of these theorems employ a ``cleaning'' process
of precisely the form just discussed.\footnote{
   We learned the ``cleaning'' interpretation
   of the Dobrushin uniqueness theorem
   from Michael Aizenman in the mid-1980s.
}
This led us to investigate the cleaning process in its own right.


In this paper we shall not treat the case of an arbitrary family $\scrt$
of cleaning operators\footnote{
   The general case leads, in fact, to interesting issues of computational
   complexity and decidability
  (see \cite{Blondel_00} for an excellent survey of closely related problems).
   Consider the following problem:
   \begin{quote}
      {\em Input:}\/ A finite set $T_1,\ldots,T_m$ of $n \times n$
           matrices with nonnegative rational entries.

      {\em Question:}\/  Does there exist a sequence of indices
           $i_1,\ldots,i_k$ such that the product $T_{i_1} \cdots T_{i_k}$
           has spectral radius $< 1$?
   \end{quote}
   This problem turns out to be {\em NP-hard}\/,
   even when restricted to matrices with elements 0 and 1
   \cite[Remarks 2 and 3 after Theorem 2]{Blondel_97},
   i.e.\ an oracle for solving it would permit the
   polynomial-time solution of any problem in the class NP
   (e.g.\ the traveling salesman problem).
   Even more strikingly, it has very recently been proven \cite{Blondel_05}
   that this problem is {\em algorithmically undecidable}\/
   even when restricted to $m=3$, $n=46$
   (indeed, even when $T_1,T_2$ are stochastic matrices
     and $T_3$ is a diagonal matrix with a single nonzero element).
   The proof uses a simple reduction from \cite[Theorem 2.1]{Blondel_03}.
},
but shall focus on the special case
of ``single-site'' cleaning operators:
for each $x \in X$ we are given exactly one cleaning operator $\beta_x$,
which leaves untouched the dirt on sites other than $x$
and which distributes the dirt on $x$ to sites $y$
with a weight factor $\alpha_{xy}$.
In other words,
\be
   (c\beta_x)_y  \;\bydef\;
   \cases{c_y + c_x \alpha_{xy}   & if $y \neq x$  \cr
          c_x \alpha_{xx}         & if $y = x$     \cr
         }
 \label{eq1.1}
\ee
This is the case that arises in the
proof of the Dobrushin uniqueness theorem.
The definition \reff{eq1.1} can trivially be rewritten as
\be
   \beta_x   \;=\;  I_{\{x\}^c} \,+\, I_{\{x\}} \alpha   \;,
 \label{eq1.2}
\ee
where $I_{\{x\}}$ and $I_{\{x\}^c}$ are the projection operators
on $\{x\}$ and its complement, respectively.
This way of writing the cleaning operators
brings out the close relations between our subject
and probabilistic potential theory
\cite{Revuz_75,Nummelin_84,Spitzer_76,Kemeny_66}.
Indeed, probabilistic potential theory in discrete time
and countable state space
can be interpreted as the theory of the algebra of operators
generated by a single nonnegative matrix
$\alpha = (\alpha_{xy})_{x,y\in X}$ together with
all the multiplication operators $I_f$
(where $f$ is a real-valued function on $X$) ---
or more specifically, as the theory of the multiplicative convex cone
of operators generated by $\alpha$ together with all the
nonnegative multiplication operators $I_f$.

Our main result (Theorem~\ref{thm.main.operators}) is that,
under mild conditions on the matrix $\alpha$ (see Section~\ref{sec2.2}),
{\em any}\/ sequence of cleaning operations inside $\Lambda$
that visits each site of $\Lambda$ infinitely many times
will lead, in the limit, to the {\em same}\/ result:
the dirt will be removed from $\Lambda$ and transferred to $\Lambda^c$
as specified by the ``balayage operator'' $\Pi_\Lambda$.

The plan of this paper is as follows:
In Section~\ref{sec2} we set forth the basic definitions
and state a few of our main results.
In Section~\ref{sec3} we analyze the cleaning operators
by deriving matrix identities and inequalities
in the spirit of probabilistic potential theory
\cite{Revuz_75,Nummelin_84,Spitzer_76,Kemeny_66};
here $\alpha = (\alpha_{xy})_{x,y\in X}$
is considered to be a fixed matrix of nonnegative real numbers.
In Section~\ref{sec4} we introduce an alternate approach
that we think clarifies the combinatorial structure
of these identities and inequalities:
it is based on the tree of finite sequences of elements of $X$.
In essence, we are now treating the matrix elements $\alpha_{xy}$
as noncommutative indeterminates;
or in physical terms, we are keeping track of the entire trajectory
of each particle of dirt, and not merely its endpoint.
This approach allows a much finer analysis of
the algebra of operators generated by $\alpha$
and the multiplication operators $I_f$.\footnote{
   Here our approach mirrors the spirit of
   modern enumerative combinatorics \cite{Stanley_99},
   where generating functions are considered in the first instance
   as {\em formal power series}\/, i.e.\ as an algebraic tool
   for efficiently manipulating collections of coefficients.
   Only at a second stage might one insert specific numerical values
   for the indeterminates and worry about convergence.
}
In Section~\ref{sec.complements} we present some alternative
sufficient conditions that guarantee the cleanability of $\Lambda$.
Finally, in Section~\ref{sec.converse}
we present some converses to our results.

\section{Basic set-up}  \label{sec2}

Let $X$ be a finite or countably infinite set (assumed nonempty),
and let $\alpha = (\alpha_{xy})_{x,y\in X}$ be a nonnegative matrix
indexed by $X$.

\subsection{Definition of operators}  \label{sec2.1}

We shall employ the following classes of matrices.
(We refer to them as ``operators'', but for now we treat them
 simply as matrices.  Later we shall make clear on what space
 of vectors they act.)

\bigskip

{\bf Multiplication operators.}
For each $\Lambda \subseteq X$, we denote by $I_\Lambda$
the projection on $\Lambda$, i.e.\ the matrix
\be
  (I_\Lambda)_{xy}   \;\bydef\;   \cases{1  & if $x=y \in \Lambda$  \cr
                                    0  & otherwise             \cr
                                   }
\ee
More generally, if $f$ is a real-valued function on $X$,
we denote by $I_f$ the operator of multiplication by $f$, i.e.\ the matrix
\be
  (I_f)_{xy}   \;\bydef\; f(x) \delta_{xy}
               \;=\;   \cases{f(x)   & if $x=y$    \cr
                              0      & otherwise   \cr
                             }
\ee
Clearly $I_\Lambda = I_{\chi_\Lambda}$,
where $\chi_\Lambda$ denotes the indicator function of $\Lambda$.

\bigskip

{\bf Cleaning operators.}
For each $x \in X$, we define the {\em cleaning operator}\/ $\beta_x$ by
\be
   \beta_x   \;\bydef\;  I_{\{x\}^c} \,+\, I_{\{x\}} \alpha   \;.
 \label{def.betax}
\ee
More generally, for each $\Lambda \subseteq X$,
we define the cleaning operator $\beta_\Lambda$ by
\be
   \beta_\Lambda   \;\bydef\;  I_{\Lambda^c} \,+\, I_{\Lambda} \alpha   \;,
 \label{def.betalambda}
\ee
where $\Lambda^c \bydef X \setminus \Lambda$.
More generally yet, for each function $f \colon\, X \to [0,1]$,
we define the cleaning operator $\beta_f$ by
\be
   \beta_f   \;\bydef\;  I_{1-f} \,+\, I_{f} \alpha   \;,
 \label{def.betaf}
\ee
so that $\beta_\Lambda = \beta_{\chi_\Lambda}$.
We also introduce, for later use, the ``dual'' cleaning operators
\be
   \beta^*_f   \;\bydef\;  I_{1-f} \,+\, \alpha I_{f}
 \label{def.betastarf}
\ee
(these have no obvious physical interpretation
 but will play an important role in our analysis).
Clearly all the operators $\beta_f$ and $\beta^*_f$ are nonnegative
(i.e.\ have nonnegative matrix elements).
We have $\beta_0 = \beta^*_0 = I$ and $\beta_1 = \beta^*_1 = \alpha$.

\medskip

{\bf Remark.}  In the ``single-site cleaning problem'' as formulated
in the Introduction, the only allowed cleaning operators are the $\beta_x$
($x \in X$).
When $\sum_{x \in X} f(x) = 1$ (resp.\ $\le 1$),
one can interpret $\beta_f$ as a convex combination of the $\beta_x$
(resp.\ of the $\beta_x$ and the identity operator),
hence as the expected output from a random choice
of single-site cleaning operator.
When $\sum_{x \in X} f(x) > 1$, the operators $\beta_f$
have no such interpretation in terms of single-site cleaning;
nevertheless, their introduction is natural from the point of view of
probabilistic potential theory, as we shall see.
Indeed, the quantity $\sum_{x \in X} f(x)$
plays no role whatsoever in our analysis;
the condition $0 \le f \le 1$, on the other hand, is crucial.
Whenever we write $\beta_f$ it will be assumed tacitly that $0 \le f \le 1$.

\bigskip

{\bf Balayage operators.}
For each $\Lambda \subseteq X$ and each $n \ge 0$, we define
$\Pi_\Lambda^{(n)}$ to be the result of
cleaning $n$ times the set $\Lambda$ and then keeping only the
dirt outside $\Lambda$:
\be
   \Pi_\Lambda^{(n)}  \;\bydef\;  \beta_\Lambda^n I_{\Lambda^c}
            \;=\;  \sum_{k=0}^n (I_\Lambda \alpha)^k I_{\Lambda^c}
   \;.
 \label{def.Pilambdan}
\ee
Let us note the identities
\begin{eqnarray}
   \beta_\Lambda^n   & = &
      \Pi_\Lambda^{(n)}  \,+\, (I_\Lambda \alpha I_\Lambda)^n I_\Lambda
     \\[2mm]
   \Pi_\Lambda^{(n)} \, \Pi_\Lambda^{(m)} & = &  \Pi_\Lambda^{(n)}
      \label{eq.Pilambdan.proj}
\end{eqnarray}
We then define the {\em balayage operator}\/ $\Pi_\Lambda$
by a limiting process:
\be
   \Pi_\Lambda  \;\bydef\;  \lim_{n\to\infty} \!\!\uparrow \Pi_\Lambda^{(n)}
                \;=\;  \sum_{k=0}^\infty (I_\Lambda \alpha)^k I_{\Lambda^c}
 \label{def.Pilambda}
\ee
Please note that, with no hypotheses other than the nonnegativity of
$\alpha$, the matrix elements of $\Pi_\Lambda$ are well-defined
(and nonnegative) but might be $+\infty$.
However, we shall soon introduce a condition
(the Fundamental Hypothesis)
under which the matrix elements of $\Pi_\Lambda$ are finite
and indeed $\Pi_\Lambda$ is a contraction operator on a suitable space
of dirt vectors (Lemma~\ref{lemma.PiLambda}).
It furthermore follows from \reff{eq.Pilambdan.proj}
that $\Pi_\Lambda$ is a projection operator,
i.e.\ $\Pi_\Lambda^2 = \Pi_\Lambda$
(see also Lemma~\ref{lemma.PiLambda.properties}).

\subsection{The fundamental hypothesis}   \label{sec2.2}

{}From now through the end of Section~\ref{sec4}
(with the exception of Sections~\ref{subsec.matrix.idineq}
 and \ref{subsec.matrix.converse}),
we shall make the following
\begin{quote}
  {\bf Fundamental Hypothesis.}
  There exists a vector $w = (w_x)_{x \in X}$ with $w_x > 0$ for all $x$,
  such that $\alpha w \le w$.
\end{quote}
In potential theory, a vector $w$ satisfying $\alpha w \le w$ is called 
{\em subinvariant}\/, {\em superharmonic}\/ or {\em excessive}\/
with respect to $\alpha$ \cite{Revuz_75,Nummelin_84,Spitzer_76,Kemeny_66}.
Thus, the Fundamental Hypothesis asserts the existence of a
strictly positive subinvariant vector.

{\bf For clarity, all results in this paper that assume the
Fundamental Hypothesis will be marked ``(FH)''.}

\bigskip

{\bf Remarks.}
1.  The Fundamental Hypothesis implies that the operator
$P \bydef I_w^{-1} \alpha I_w$ is submarkovian,
i.e.\ satisfies $P {\bf 1} \le {\bf 1}$.
Our results could therefore be given a more probabilistic flavor,
reminiscent of Markov-chain theory,
by rewriting them in terms of $P$
(this is easy because $I_w$ commutes with all the
 other multiplication operators $I_f$).
We shall not need this interpretation, however,
so we leave this translation to the reader.

2.  If $X$ is finite {\em and $\alpha$ is irreducible}\/,
the Fundamental Hypothesis means simply that the
Perron--Frobenius eigenvalue (= spectral radius) of $\alpha$ is $\le 1$.
In this case it is natural to take $w$ to be the
Perron--Frobenius eigenvector
(this is the unique choice if the spectral radius equals 1,
 but is nonunique otherwise).
But if $\alpha$ fails to be irreducible, the Fundamental Hypothesis
is {\em stronger}\/ than this assertion
about the spectral radius:  consider, for instance,
$\alpha = \left(\!\! \begin{array}{cc} 1 & a \\ 0 & b \end{array} \!\!\right)$
with $a > 0$ and $0 \le b \le 1$,
which has spectral radius 1 but does not satisfy the Fundamental Hypothesis.
It can be shown \cite{Rothblum_private} that, when $X$ is finite,
the Fundamental Hypothesis holds
if and only if the spectral radius of $\alpha$ is $\le 1$
{\em and}\/
each class $J$ for which the square submatrix $\alpha_{JJ}$
has spectral radius 1 is a final class
(i.e.\ $\alpha_{jk}=0$ whenever $j \in J$ and $k \notin J$).\footnote{
   We recall that the {\em classes}\/ of a nonnegative matrix $\alpha$
   are defined as the strongly connected components
   of the directed graph with vertex set $X$
   and edge set $\{ xy \colon\; \alpha_{xy} > 0 \}$:
   see e.g.\ \cite{Berman_79,Seneta_81}.
}

3. Example~\ref{example_sec6_FHfails} below shows that if $X$ is infinite,
then the Fundamental Hypothesis can fail even though each matrix
$I_\Lambda \alpha I_\Lambda$ for $\Lambda$ finite
has spectral radius 0 (i.e.\ is nilpotent).

\bigskip

{\em We fix once and for all a vector $w > 0$ satisfying $\alpha w \le w$.}\/
For any vector $c = (c_x)_{x \in X}$, we define
\be
   \| c \|_w  \;\bydef\;  \sum_{x \in X} |c_x| \, w_x
\ee
and we denote by $l^1(w)$ the space of vectors $c$ satisfying
$\| c \|_w < \infty$.

For any matrix $A = (A_{xy})_{x,y \in X}$, we define the operator norm
\be
   \| A \|_{w \to w}  \;\bydef\;
   \sup_{\begin{scarray}
           c\;\hbox{\scriptsize of finite support} \\
           c \neq 0
         \end{scarray}}
   {\| cA \|_w  \over \| c \|_w}   \;.
 \label{def.operatornorm}
\ee
If $\| A \|_{w \to w} < \infty$,
then $A$ defines a bounded operator on $l^1(w)$ of norm $\| A \|_{w \to w}$.
It is not hard to see that \reff{def.operatornorm} is equivalent to
\be
   \| A \|_{w \to w}  \;=\;
   \sup_{x \in X}  w_x^{-1} \sum\limits_{y \in X} |A_{xy}| \, w_y
   \;.
 \label{def.operatornorm_bis}
\ee
In particular, when $A$ is a nonnegative matrix,
\reff{def.operatornorm}/\reff{def.operatornorm_bis} reduces simply to
\be
   \| A \|_{w \to w}  \;=\;
   \sup_{x \in X}  {(Aw)_x \over w_x}   \;.
 \label{def.operatornorm2}
\ee
Thus, for nonnegative matrices,
$\| A \|_{w \to w} \le 1$ if and only if $Aw \le w$.

The norms of the multiplication operators are trivially given by
\be
   \| I_f \|_{w \to w}   \;=\;  \|f\|_\infty
        \;\bydef\;  \sup_{x \in X} |f(x)|  \;.
  \label{eq.norm.mult}
\ee
The Fundamental Hypothesis immediately gives
$\| \alpha \|_{w \to w} \le 1$.
More generally, for any $f \colon\, X \to [0,1]$,
the cleaning operator $\beta_f$ satisfies
\be
   \beta_f w  \;=\;  I_{1-f} w + I_f \alpha w  \;\le\;  I_{1-f} w + I_f w
        \;=\; w
  \label{eq.norm.beta}
\ee
and hence $\| \beta_f \|_{w \to w} \le 1$.
A similar argument shows, yet more generally,
that for any functions $(f_i)_{i=1}^n$ and operators $(A_i)_{i=1}^n$, we have
\be
   \left\| \sum\limits_{i=1}^n I_{f_i} A_i \right\|_{w \to w}
   \;\le\;
   \left\| \sum\limits_{i=1}^n \|A_i\|_{w \to w} \, |f_i| \right\| _\infty
   \;\le\;
   \sum\limits_{i=1}^n  \|A_i\|_{w \to w} \, \|f_i\|_\infty
   \;.
\ee

Finally, for the balayage operator $\Pi_\Lambda$ we have the easy result:
\begin{lemma}
 \label{lemma.PiLambda}
\FH
For each $\Lambda \subseteq X$,
the matrix $\Pi_\Lambda$ has finite matrix elements,
and indeed satisfies 
\be
   0 \;\le\; \Pi_\Lambda w \;\le\; w   \;,
\ee
so that $\| \Pi_\Lambda \|_{w \to w}  \,\le\,   1$.
More specifically, $\| \Pi_\Lambda \|_{w \to w}$
equals $1$ if $\Lambda \neq X$, and equals $0$ if $\Lambda = X$.
\end{lemma}

\proof
By \reff{eq.norm.beta} we have
$\Pi_\Lambda^{(n)} w = \beta_\Lambda^n I_{\Lambda^c} w
  \le \beta_\Lambda^n w \le w$.
Since $0 \le \Pi_\Lambda^{(n)} \uparrow \Pi_\Lambda$ elementwise,
we have $\Pi_\Lambda^{(n)} w \uparrow \Pi_\Lambda w$ elementwise
by the monotone convergence theorem, and hence $\Pi_\Lambda w \le w$.
Since $\Pi_\Lambda$ is nonnegative,
we conclude from \reff{def.operatornorm} or \reff{def.operatornorm_bis} that
$\| \Pi_\Lambda \|_{w \to w}  \,\le\,   1$.

On the other hand,
we have $c \Pi_\Lambda = c$ for any vector $c$ supported on $\Lambda^c$;
when $\Lambda \neq X$ (so that $c$ can be chosen nonzero)
this implies that $\| \Pi_\Lambda \|_{w \to w}  \,\ge\,   1$.
When $\Lambda = X$, by contrast, we have $\Pi_\Lambda = 0$.
\qed

\subsection{Main results}

Let us now state briefly a few of our main results,
just to give their flavor.
Our principal result on cleaning is the following:

\begin{theorem}[= Theorem~\ref{thm.main.operators}]
\FH
Consider a region $\Lambda \subseteq X$
and functions $0\le h_i\le \chi_\Lambda$ ($i\ge 1$) such that
\be
\sum_{i=1}^\infty h_i(x) \;=\; \infty\qquad\hbox{for all } x\in\Lambda\;.
\ee
If the nonnegative vector $c \in l^1(w)$ is such that
$\| c (I_\Lambda \alpha I_\Lambda)^\ell \|_w \to 0$
as $\ell\to\infty$, then
\be
\|c ( \beta_{h_1} \cdots \beta_{h_{N}}  \,-\, \Pi_\Lambda)\|_w
\;\mathop{\longrightarrow}\limits_{N\to\infty}\; 0\;.
  \label{eq.sec2.cleaning}
\ee
\end{theorem}

\noindent
Intuitively, this says that any sequence of cleaning operations
inside $\Lambda$ that covers $\Lambda$ infinitely many times
will lead, in the limit, to removing the dirt from $\Lambda$
and transferring it to $\Lambda^c$
as specified by the balayage operator $\Pi_\Lambda$.

We also have a variant of this result
(Theorem~\ref{thm.main.operators.uniform})
in which the vector-norm convergence \reff{eq.sec2.cleaning}
is strengthened to operator-norm convergence,
at the price of a stronger hypothesis on the functions $h_i$
(namely, requiring that they cover the set $\Lambda$ uniformly).
Finally, we will prove some far-reaching extensions of this result,
in terms of the tree formalism, in Section~\ref{sec4.8}.

We also have some sufficient conditions for cleanability
when the Fundamental Hypothesis is not assumed.  For instance:

\begin{theorem}[= Theorem~\ref{thm.complements.2}]
Let $\Lambda \subseteq X$,
and let $c \ge 0$ and $w \ge 0$ be vectors satisfying
\begin{itemize}
   \item[(a)] $c (I_\Lambda \alpha I_\Lambda)^k w \,<\, \infty$
      for all $k \ge 0$
\end{itemize}
and
\begin{itemize}
   \item[(b)] $\lim\limits_{k\to\infty} c (I_\Lambda \alpha I_\Lambda)^k w
               \,=\, 0$.
\end{itemize}
Then it is possible to find a sequence of sites $x_1,x_2,\ldots \in \Lambda$
and a sequence of numbers $\epsilon_1,\epsilon_2,\ldots \in (0,1]$
such that
\begin{equation}
  \lim\limits_{n \to\infty} c\beta_{\epsilon_1 \delta_{x_1}}
     \cdots \beta_{\epsilon_n \delta_{x_n}} I_\Lambda w
  \;=\;  0
  \;.
\end{equation}
\end{theorem}

\noindent
Although this result refers only to matrices,
its proof uses the tree formalism (see Section~\ref{sec.complements}).

Finally, we have a converse result when the Fundamental Hypothesis
is not assumed:

\begin{theorem}[= Theorem~\ref{thm.matrix.converse}]
Let $X$ be a finite or countably infinite set,
let $\Lambda \subseteq X$, and let $c \ge 0$ and $w \ge 0$
be vectors that are strictly positive on $\Lambda$.
Consider the following conditions on a matrix $\alpha$:
\begin{itemize}
  \item[(a)]  $\sum_{k=0}^\infty c (I_\Lambda \alpha I_\Lambda)^k w < \infty$.
  \item[(b)]  For all $h \colon\, X \to [0,1]$ with $\supp h = \Lambda$
      such that $h \ge \epsilon \chi_\Lambda$ for some $\epsilon > 0$,
      we have
      $\sum_{k=0}^\infty c (I_\Lambda \beta_h I_\Lambda)^k w < \infty$.
  \item[(b${}'$)]  There exists $h \colon\, X \to [0,1]$
      with $\supp h = \Lambda$ such that
      $\sum_{k=0}^\infty c (I_\Lambda \beta_h I_\Lambda)^k w < \infty$.
  \item[(c)]  For every finite sequence $f_1,\ldots,f_m$
     of functions $X \to [0,1]$ with $\supp(f_i) \subseteq \Lambda$
     such that $\sum_i f_i \ge \epsilon \chi_\Lambda$
     for some $\epsilon > 0$, we have
     $\sum_{k=0}^\infty c (I_\Lambda \beta_{f_1} \cdots
                                     \beta_{f_m} I_\Lambda)^k w < \infty$.
  \item[(c${}'$)]  There exists a finite sequence $f_1,\ldots,f_m$
     of functions $X \to [0,1]$ with $\supp(f_i) \subseteq \Lambda$
     such that
     $\sum_{k=0}^\infty c (I_\Lambda \beta_{f_1} \cdots
                                     \beta_{f_m} I_\Lambda)^k w < \infty$.
\end{itemize}
Then (a)--(c) are all equivalent and imply (c${}'$);
and for matrices $\alpha$ satisfying the additional hypothesis
\be
  \hbox{There exists a constant $C < \infty$ such that
        $I_\Lambda \alpha I_\Lambda w \le C w$}
\ee
all five conditions are equivalent.
\end{theorem}

\noindent
Intuitively, this says (at least if $\Lambda$ is a finite set)
that if the spectral radius of
$I_\Lambda \alpha I_\Lambda$ is $\ge 1$, then there is no way to clean the
set $\Lambda$ completely.
In fact, we shall prove a stronger version of this result
in Section~\ref{sec.converse}, using the tree formalism
(see Theorem~\ref{thm.converse.2.NEW} and
 Corollary~\ref{cor.converse.finite.lambda}).

\section{Matrix approach to balayage}   \label{sec3}

In this section we study the algebra of operators generated
by $\alpha$ and the multiplication operators $I_f$,
by deriving matrix identities and inequalities
in the spirit of probabilistic potential theory
\cite{Revuz_75,Nummelin_84,Spitzer_76,Kemeny_66}.
Our main goal is to study the convergence of a product
of cleaning operators $\beta_{h_1} \cdots \beta_{h_n}$
where all the $h_i$ have support contained in $\Lambda$.
We shall show that, under very general conditions,
any such product must converge to $\Pi_\Lambda$
(see Theorems~\ref{thm.main.operators.uniform} and \ref{thm.main.operators}).
We shall also prove a converse result
when the Fundamental Hypothesis is not assumed
(Theorem~\ref{thm.matrix.converse}).


\subsection{Identities for cleaning operators}

\begin{proposition}[Fundamental identities]
  \label{prop3.1}
For all $h,h_1,h_2 \colon\, X \to [0,1]$, the following identities hold:
\begin{itemize}
\item[(i)] \emph{Intertwining I:}
  \begin{equation}
    \label{eq:9}
    (I-\alpha)\,\beta_h \;=\; \beta^*_h\, (I-\alpha)
  \end{equation}

\item[(ii)] \emph{Intertwining II:}
  \begin{eqnarray}
    (\beta_{h_1} - I)\, I_{h_2}  & = &   I_{h_1}\, (\beta^*_{h_2} - I)
        \label{eq.intertwiningIIa}  \\[2mm]
    \beta_h\,I_h   & = &   I_h\,\beta^*_h
        \label{eq:41}
  \end{eqnarray}


\item[(iii)] \emph{Comparison:}
  \begin{eqnarray}
    \label{eq:11}
    \beta_{h_1} - \beta_{h_2} & = &  I_{h_2-h_1}\, (I-\alpha)  \\[2mm]
    I - \beta_h  & = & I_h (I-\alpha)
    \label{eq:11a}
  \end{eqnarray}

\item[(iv)] \emph{Collapse:}
  \begin{equation}
    \label{eq:13}
    \beta_{h_1}\,\beta_{h_2} \;=\; \beta_{1-(1-h_1)(1-h_2)} - 
I_{h_1}\,\alpha\,I_{h_2}\,(I-\alpha)
  \end{equation}
\end{itemize}
\end{proposition}

\proof
\noindent{(i)}
\begin{equation}
  \label{eq:10}
  \alpha\,\beta_h-\beta^*_h\,\alpha \;=\; 
\alpha\,I_{1-h} - I_{1-h}\,\alpha
\;=\; I_h\,\alpha - \alpha\, I_h \;=\; \beta_h-\beta^*_h \;.
\end{equation}

\noindent{(ii)}
Both sides of \reff{eq.intertwiningIIa} are equal to
$I_{h_1} (\alpha-I) I_{h_2}$.
Equation \reff{eq:41} follows by setting $h_1 = h_2 = h$.

\bigskip 

\noindent{(iii)}
\begin{equation}
  \label{eq:12}
  \beta_{h_1}-\beta_{h_2} \;=\; I_{1-h_1}-I_{1-h_2} +
(I_{h_1}-I_{h_2})\,\alpha \;=\;
I_{h_2-h_1}\, (I-\alpha)  \;.
\end{equation}
Equation \reff{eq:11a} follows by setting $h_1 = 0$, $h_2 = h$.

\noindent{(iv)}
\begin{eqnarray}
  \label{eq:14}
  \beta_{h_1}\,\beta_{h_2} +
  I_{h_1}\,\alpha\,I_{h_2}\,(I-\alpha) &=&
I_{(1-h_1)(1-h_2)} + I_{h_1}\,\alpha\,I_{1-h_2} +
I_{1-h_1}\,I_{h_2}\,\alpha + I_{h_1}\,\alpha\,I_{h_2}
\nonumber\\
&=& I_{(1-h_1)(1-h_2)} + I_{h_1}\,\alpha +
I_{1-h_1}\,I_{h_2}\,\alpha \nonumber\\
&=& \beta_{1-(1-h_1)(1-h_2)}  \;.
\end{eqnarray}
\qed

\begin{corollary}[Telescoping comparison]
   \label{cor_telescoping}
For all $g_1,\ldots,g_n,h_1,\ldots,h_n \colon\, X \to [0,1]$, we have
\begin{subeqnarray}
   \beta_{g_1}\cdots\beta_{g_n} -  \beta_{h_1}\cdots\beta_{h_n} 
   & = &  
   \sum_{i=1}^n   \beta_{g_1}\cdots\beta_{g_{i-1}} 
          I_{h_i - g_i} (I-\alpha) \beta_{h_{i+1}}\cdots\beta_{h_n} 
       \slabel{eq.telescoping_identity}  \\[2mm]
   & = &
   \sum_{i=1}^n   \beta_{g_1}\cdots\beta_{g_{i-1}} 
          I_{h_i - g_i} \beta^*_{h_{i+1}}\cdots\beta^*_{h_n} (I-\alpha) 
\end{subeqnarray}
\end{corollary}

\proof
This is an immediate consequence of the telescopic decomposition
\be
   \beta_{g_1}\cdots\beta_{g_n} -  \beta_{h_1}\cdots\beta_{h_n}
\;=\; \sum_{i=1}^n   \beta_{g_1}\cdots\beta_{g_{i-1}}
(\beta_{g_i}-\beta_{h_i})   \beta_{h_{i+1}}\cdots\beta_{h_n}
   \;,
\ee
the comparison identity \reff{eq:11},
and the intertwining relation \reff{eq:9}.
\qed

\begin{lemma}[Cleaners with restricted support]
  \label{lemma_betagi}
Suppose that $\supp(g_i) \subseteq \Lambda$ for $i=1,\ldots,n$,
and let $h_i \colon\, X \to \R$ be any functions satisfying
$h_i \restrict \Lambda \equiv 1$ for $i=1,\ldots,n$.
Then
\be
   \beta_{g_1} \,\cdots\, \beta_{g_n} \, I_\Lambda
   \;=\;
   I_{h_1} \, \beta_{g_1} \,\cdots\, I_{h_n} \, \beta_{g_n} \, I_\Lambda
 \label{eq.betagi.1}
\ee
and
\be
   I_\Lambda \, \beta^*_{g_1} \,\cdots\, \beta^*_{g_n} \;=\;
   I_\Lambda \, \beta^*_{g_1} \, I_{h_1} \,\cdots\, \beta^*_{g_n} \, I_{h_n}
   \;.
 \label{eq.betagi.2}
\ee
Furthermore,
\be\label{roma.r10}
I_{\Lambda^c} \beta_{g_1} \cdots \beta_{g_n}
\;=\; I_{\Lambda^c} \;=\;
 \beta^*_{g_1} \cdots \beta^*_{g_n} I_{\Lambda^c}
\;.
\ee
\end{lemma}

\proof
Note first that if $\supp(g) \subseteq \Lambda$,
we have $I_{\Lambda^c} \beta_g I_\Lambda = 0$;
therefore, for any $h$ satisfying $h \restrict \Lambda \equiv 1$, we have
\be
   \beta_g I_\Lambda  \;=\;
   I_\Lambda \beta_g I_\Lambda  \;=\;
   I_h \beta_g I_\Lambda  \;.
 \label{eq.betagi.3}
\ee
So, starting with $\beta_{g_1} \cdots \beta_{g_n} I_\Lambda$,
we successively use the first equality in \reff{eq.betagi.3},
working from right to left in the product, to transform it into
$I_\Lambda \beta_{g_1} \cdots I_\Lambda \beta_{g_n} I_\Lambda$.
We then successively use the second equality in \reff{eq.betagi.3},
working from left to right in the product, to transform it into
$I_{h_1} \beta_{g_1} \cdots I_{h_n} \beta_{g_n} I_\Lambda$.
This proves \reff{eq.betagi.1}.
An analogous argument gives \reff{eq.betagi.2}.
Identities \reff{roma.r10} are immediately verified by induction.
%
\qed

\subsection{Inequalities for cleaning operators}

We now turn our attention to proving inequalities that say,
roughly speaking, that one operator ``cleans better'' than another.

\begin{lemma}
   \label{lemma.2and3}
\FH
Let $\Lambda$ be any subset of $X$,
and let $f$ be any function satisfying $\chi_\Lambda \le f \le 1$.
Then
\be
   I_\Lambda \, (I-\alpha) \, I_f \, w   \;\ge\;   0   \;.
 \label{eq.lemma2}
\ee
If, in addition, $0 \le h_i \le \chi_\Lambda$ for $i=1,\ldots,n$, then
\be
   I_\Lambda \, (I-\alpha) \beta_{h_1} \cdots \beta_{h_n} \, I_f \, w
   \;\ge\;   0   \;.
 \label{eq.lemma3}
\ee
\end{lemma}


\proof
We have
\begin{eqnarray}
   I_\Lambda \, (I-\alpha) \, I_f \, w
          & = &  I_\Lambda f w \,-\, I_\Lambda \alpha f w    \nonumber \\
          &\ge&  I_\Lambda f w \,-\, I_\Lambda \alpha w
                            \qquad \hbox{[since $f \le 1$]}  \nonumber \\
          &\ge&  I_\Lambda f w \,-\, I_\Lambda w
                 \qquad\:\:\, \hbox{[since $\alpha w \le w$]} \nonumber \\
          & = &  I_\Lambda (f-1) w                           \nonumber \\
          & = &  0     
           \qquad\qquad\qquad\qquad\! \hbox{[since $f \equiv 1$ on $\Lambda$]}
\end{eqnarray}
This proves \reff{eq.lemma2}.  Then
\begin{eqnarray}
   I_\Lambda (I-\alpha) \beta_{h_1} \cdots \beta_{h_n} I_f w
      & = &   I_\Lambda \beta^*_{h_1} \cdots \beta^*_{h_n} (I-\alpha) I_f w
                      \nonumber \\[1mm]
      & = &   I_\Lambda \beta^*_{h_1} \cdots \beta^*_{h_n} I_\Lambda
                  (I-\alpha) I_f w
                      \nonumber \\[1mm]
      & \ge &  0
\end{eqnarray}
where the first equality uses the intertwining relation \reff{eq:9},
the second uses the identity \reff{eq.betagi.2},
and the final inequality uses \reff{eq.lemma2}
and the nonnegativity of $\beta^*_{h_1},\ldots,\beta^*_{h_n}$.
\qed

%


\begin{proposition}[Multi-monotonicity]
  \label{prop.multi-monotonicity}
\FH
Suppose that $0 \le g_i \le h_i \le \chi_\Lambda \le f \le 1$
for $i=1,\ldots,n$.  Then
\be
   \beta_{h_1} \cdots \beta_{h_n} I_f w  \;\le\;
   \beta_{g_1} \cdots \beta_{g_n} I_f w  \;.
 \label{eq.multi1}
\ee
\end{proposition}

\proof
We will prove \reff{eq.multi1} by proving separately that
$I_{\Lambda^c}({\rm LHS}) \le I_{\Lambda^c}({\rm RHS})$
and that
$I_{\Lambda}({\rm LHS}) \le I_{\Lambda}({\rm RHS})$.
The former is in fact equality, since
$I_{\Lambda^c} \beta_{h_1} \cdots \beta_{h_n} =
 I_{\Lambda^c} \beta_{g_1} \cdots \beta_{g_n} = I_{\Lambda^c}$.
The latter follows from the
telescoping comparison identity \reff{eq.telescoping_identity}
and the inequality \reff{eq.lemma3},
along with the nonnegativity of the operators
$\beta_{g_1},\ldots,\beta_{g_{i-1}}$
that lie on the left in \reff{eq.telescoping_identity}.
\qed

\begin{proposition}[Collapse inequality]
  \label{prop.collapse}
\FH
Suppose that $0 \le g_i, h_j \le \chi_\Lambda \le f \le 1$
for $i=1,\ldots,n$ and $j=1,\ldots,m$.  Then
\be
  \beta_{g_1}\cdots\beta_{g_n} \beta_{h_1}\cdots\beta_{h_m} I_f w
  \;\le\;
  \beta_{1-\prod_{i=1}^n (1-g_i)}  \beta_{h_1}\cdots\beta_{h_m} I_f w
  \;.
 \label{eq.coll1}
\ee
\end{proposition}

\proof
For $n=1$ this is trivial.
For $n=2$ it follows from the collapse identity \reff{eq:13}
together with the inequality \reff{eq.lemma3}.
The cases $n \ge 3$ are obtained by an easy induction from the case $n=2$.
\qed

Combining Propositions~\ref{prop.multi-monotonicity} and
\ref{prop.collapse}, we obtain the following comparison result:

\begin{corollary}[Multi-monotonicity + collapse]
   \label{cor.multi-monotonicity+collapse}
\FH
Suppose that $0 \le h_i, g_j \le \chi_\Lambda \le f \le 1$
for $i=1,\ldots,N$ and $j=1,\ldots,k$.
Suppose further that there exist integers
$0 \le n_0 < n_1 < \ldots < n_k \le N$ such that
\be
   1 - \!\prod\limits_{i=n_{j-1}+1}^{n_j} (1-h_i)   \;\ge\;   g_j
\ee
for all $j$.
Then
\be
   \beta_{h_1} \cdots \beta_{h_N} I_f w  \;\le\;
   \beta_{g_1} \cdots \beta_{g_k} I_f w  \;.
 \label{eq.multi+collapse1}
\ee
\end{corollary}

\bigskip

Now we develop some analogous inequalities going in the {\em reverse}\/
direction provided that we look only at the dirt {\em outside}\/ $\Lambda$;
moreover, these inequalities hold pointwise.

\begin{lemma}
   \label{lemma.4and5}
Let $\Lambda$ be any subset of $X$.  Then
\be
   I_\Lambda \, (I-\alpha) \, I_{\Lambda^c}  \;\le\;  0   \;.
 \label{eq.lemma4}
\ee
If, in addition, $0 \le h_i \le \chi_\Lambda$  for $i=1,\ldots,n$, then
\be
   I_\Lambda \, (I-\alpha) \beta_{h_1} \cdots \beta_{h_n} \, I_{\Lambda^c}
   \;\le\;   0   \;.
 \label{eq.lemma5}
\ee
\end{lemma}

\proof
The inequality \reff{eq.lemma4} is trivial since $\alpha \ge 0$.
To prove \reff{eq.lemma5}, we make a computation analogous to that
in the proof of Lemma~\ref{lemma.2and3}:
\begin{eqnarray}
   I_\Lambda (I-\alpha) \beta_{h_1} \cdots \beta_{h_n} I_{\Lambda^c}
         & = &   I_\Lambda \beta^*_{h_1} \cdots \beta^*_{h_n} (I-\alpha)
                     I_{\Lambda^c}    \nonumber \\[1mm]
         & = &   I_\Lambda \beta^*_{h_1} \cdots \beta^*_{h_n} I_\Lambda
                     (I-\alpha) I_{\Lambda^c}    \nonumber \\[1mm]
         & \le &  0
\end{eqnarray}
where the final step uses \reff{eq.lemma4}.
\qed

\begin{proposition}[Reverse multi-monotonicity]
  \label{prop.reverse.multi-monotonicity}
Suppose that $0 \le h_i \le g_i \le \chi_\Lambda$ for $i=1,\ldots,n$.  Then
\be
   \beta_{g_1} \cdots \beta_{g_n} I_{\Lambda^c}  \;\ge\;
   \beta_{h_1} \cdots \beta_{h_n} I_{\Lambda^c}  \;.
 \label{eq.multi3}
\ee
\end{proposition}

\proof
We will prove \reff{eq.multi3} by proving separately that
$I_{\Lambda^c}({\rm LHS}) \le I_{\Lambda^c}({\rm RHS})$
and that
$I_{\Lambda}({\rm LHS}) \le I_{\Lambda}({\rm RHS})$.
The former is in fact equality because of \reff{roma.r10}.
The latter follows from the
telescoping comparison identity \reff{eq.telescoping_identity}
and the inequality \reff{eq.lemma5},
along with the nonnegativity of the operators
$\beta_{h_1},\ldots,\beta_{h_{i-1}}$
that lie on the left in \reff{eq.telescoping_identity}.
\qed

\begin{proposition}[Reverse collapse inequality]
  \label{prop.reverse.collapse}
Suppose that $0 \le g_i, h_j \le \chi_\Lambda$
for $i=1,\ldots,n$ and $j=1,\ldots,m$.  Then
\be
 \label{eq.collrev}
  \beta_{g_1}\cdots\beta_{g_n} \beta_{h_1}\cdots\beta_{h_m} I_{\Lambda^c}
  \;\ge\;
  \beta_{1-\prod_{i=1}^n (1-g_i)}  \beta_{h_1}\cdots\beta_{h_m} I_{\Lambda^c}
  \;.
\ee
\end{proposition}

\proof
For $n=1$ this is trivial.
For $n=2$ it follows from the collapse identity \reff{eq:13}
together with the inequality \reff{eq.lemma4}.
The cases $n \ge 3$ are obtained by an easy induction from the case $n=2$.
\qed

Combining Propositions~\ref{prop.reverse.multi-monotonicity} and
\ref{prop.reverse.collapse}, we obtain:

\begin{corollary}[Reverse multi-monotonicity + collapse]
   \label{cor.reverse.multi-monotonicity+collapse}
Suppose that $0 \le h_i, g_j \le \chi_\Lambda$
for $i=1,\ldots,N$ and $j=1,\ldots,k$.
Suppose further that there exist integers
$0 \le n_0 < n_1 < \ldots < n_k \le N$ such that
\be
   1 - \!\prod\limits_{i=n_{j-1}+1}^{n_j} (1-h_i)   \;\ge\;   g_j
\ee
for all $j$.
Then
\be
   \beta_{h_1} \cdots \beta_{h_N} I_{\Lambda^c}  \;\ge\;
   \beta_{g_1} \cdots \beta_{g_k} I_{\Lambda^c}  \;.
 \label{eq.reverse.multi+collapse1}
\ee
\end{corollary}

\subsection{Identities related to $\Pi_\Lambda$}

\begin{lemma}[Properties of \mbox{\protect\boldmath $\Pi_\Lambda$}]
  \label{lemma.PiLambda.properties}
Let $\Lambda \subseteq X$.  Then the following identities hold:
\begin{itemize}
\item[(i)] \emph{Basic properties of $\Pi_\Lambda$:}
 \begin{eqnarray}
    \Pi_\Lambda^2     & = &   \Pi_\Lambda               \label{eq.Pi1} \\[2mm]
    I_\Lambda (I-\Pi_\Lambda)   & = &    I-\Pi_\Lambda  \label{eq.Pi2} 
 \end{eqnarray}
\item[(ii)] \emph{Absorption of cleaning operators:}
   If $\supp(h) \subseteq \Lambda$, then
   \begin{eqnarray}
      I_h \Pi_\Lambda  & = &  I_h \alpha \Pi_\Lambda   \label{eq.Pi3}  \\[2mm]
      \beta_h \, \Pi_\Lambda  & = &  \Pi_\Lambda   \label{eq.Pi4}  \\[2mm]
      \Pi_\Lambda \, \beta_h  & = &  \Pi_\Lambda   \label{eq.Pi5}
   \end{eqnarray}
\end{itemize}
\end{lemma}

\proof
It follows immediately from the definition \reff{def.Pilambda}
of $\Pi_\Lambda$ that
$\Pi_\Lambda I_{\Lambda^c} = \Pi_\Lambda$ and
$I_{\Lambda^c} \Pi_\Lambda = I_{\Lambda^c}$,
hence
$\Pi_\Lambda \Pi_\Lambda = \Pi_\Lambda I_{\Lambda^c} \Pi_\Lambda =
 \Pi_\Lambda I_{\Lambda^c} = \Pi_\Lambda$.
This proves \reff{eq.Pi1}.
Equation~\reff{eq.Pi2} is trivially equivalent to
$I_{\Lambda^c} \Pi_\Lambda = I_{\Lambda^c}$.
It also follows immediately from the definition of $\Pi_\Lambda$ that
\be
    I_\Lambda \Pi_\Lambda  \;=\;   \Pi_\Lambda \,-\, I_{\Lambda^c}
    \;=\;
    I_\Lambda \alpha \Pi_\Lambda
   \;,
        \label{eq.Pi.star1}
\ee
and premultiplying this by $I_h$ yields \reff{eq.Pi3}.
The identity \reff{eq.Pi4} follows immediately from \reff{eq.Pi3}
and the definition of $\beta_h$.
The identity \reff{eq.Pi5} follows from
$\Pi_\Lambda \beta_h = \Pi_\Lambda I_{\Lambda^c} (I_{1-h} + I_h \alpha)
 = \Pi_\Lambda I_{\Lambda^c} = \Pi_\Lambda$.
\qed

\begin{proposition}[Convergence-to-balayage identity]
  \label{prop.convergence-to-balayage}
Suppose that $0 \le h_i \le \chi_\Lambda$ for $i=1,\ldots,n$.  Then
\begin{subeqnarray}
   \beta_{h_1} \cdots \beta_{h_n}  \,-\, \Pi_\Lambda
   & = &
   \beta_{h_1} \cdots \beta_{h_n} I_\Lambda (I - \Pi_\Lambda)
       \slabel{eq.betaPi.1}  \\[2mm]
   & = &
   (I_\Lambda \beta_{h_1} I_\Lambda) \,\cdots\,
   (I_\Lambda \beta_{h_n} I_\Lambda) \, I_\Lambda (I - \Pi_\Lambda)  \;.
       \slabel{eq.betaPi.2}
       \label{eq.betaPi.all}
\end{subeqnarray}
\end{proposition}

\proof
We shall prove \reff{eq.betaPi.1} by induction on $n$.
It is true for $n=0$, by \reff{eq.Pi2}.
So assume it is true for $n-1$, i.e.\ that
\be
   \beta_{h_2} \cdots \beta_{h_n}  \,-\, \Pi_\Lambda
   \;=\;
   \beta_{h_2} \cdots \beta_{h_n} I_\Lambda (I - \Pi_\Lambda)
\ee
Left-multiplying by $\beta_{h_1}$ and using \reff{eq.Pi4},
the desired identity follows.
The alternative form \reff{eq.betaPi.2} is then an immediate consequence
of Lemma~\ref{lemma_betagi}.
\qed

{\bf Remark.}  In the absence of the Fundamental Hypothesis,
some of the matrix elements of $\Pi_\Lambda$ could be $+\infty$,
but the identities \reff{eq.Pi1}--\reff{eq.Pi5} and \reff{eq.betaPi.all}
continue to hold
(with some matrix elements possibly $+\infty$ or $-\infty$).

\subsection{Comparison of cleaning operators with $\Pi_\Lambda$}


\begin{definition}
Let $\Lambda \subseteq X$.
We say that an operator $A \ge 0$ is
{\em absorbed by $\Pi_\Lambda$}\/
in case $A \Pi_\Lambda = \Pi_\Lambda$.
\end{definition}

\noindent
By Lemma~\ref{lemma.PiLambda.properties},
operators $\beta_{h_1} \cdots \beta_{h_n}$
with $\supp(h_i) \subseteq \Lambda$
are absorbed by $\Pi_\Lambda$,
as are all convex combinations thereof.

Operators absorbed by $\Pi_\Lambda$ obey some elementary
but remarkable identities and inequalities:

\begin{lemma}[Comparison with \mbox{\protect\boldmath $\Pi_\Lambda$}]
   \label{lemma.A_minus_Pilambda}
\FH
Let $\Lambda \subseteq X$,
and let $A \ge 0$ be an operator absorbed by $\Pi_\Lambda$.
Then:
\begin{itemize}
   \item[(i)]  $0 \,\le\, \Pi_\Lambda w \,\le\, Aw$.
      [In particular, if $\Lambda \neq X$, we have
       $\| A \| _{w \to w} \ge 1$.]
   \item[(ii)]  The operator $A - \Pi_\Lambda$ can be decomposed in the form
\begin{subeqnarray}
   A - \Pi_\Lambda
   & = &
   (A - \Pi_\Lambda) I_\Lambda \;+\; (A - \Pi_\Lambda) I_{\Lambda^c}  \\[1mm]
   & = &
   A I_\Lambda  \;-\; (\Pi_\Lambda - A I_{\Lambda^c})
       \slabel{lemma.A_minus_Pilambda.star1b}  \\[1mm]
   & = &
   A I_\Lambda (I - \Pi_\Lambda)
 \label{lemma.A_minus_Pilambda.star1}
\end{subeqnarray}
where
\begin{subeqnarray}
   A I_\Lambda  & \ge &   0
       \slabel{lemma.A_minus_Pilambda.star2a}  \\[1mm]
\hspace*{-4cm}
   \Pi_\Lambda - A I_{\Lambda^c}  \;=\;  A I_\Lambda \Pi_\Lambda
      & \ge &    0
       \slabel{lemma.A_minus_Pilambda.star2b}  \\[1mm]
   (\Pi_\Lambda - A I_{\Lambda^c}) w   & \le &  A I_\Lambda w
       \slabel{lemma.A_minus_Pilambda.star2c}
       \label{lemma.A_minus_Pilambda.star2}
\end{subeqnarray}
   \item[(iii)]  The following norm relations hold:
\be
   \| A - \Pi_\Lambda \|_{w \to w}
   \;=\;
   \|  A I_\Lambda (I-\Pi_\Lambda) \|_{w \to w}
   \;\le\;
   2 \| A I_\Lambda \|_{w \to w}
 \label{rouen.r1}
\ee
and, for every vector $c \ge 0$,
\be
   \| c (A - \Pi_\Lambda) \|_w
   \;=\;
   \| c A I_\Lambda \|_w  \,+\,
   \| c (\Pi_\Lambda - A I_{\Lambda^c}) \|_w
   \;\le\;
   2 \| c A I_\Lambda \|_w
   \;.
 \label{lemma.A_minus_Pilambda.star3}
\ee
\end{itemize}
\end{lemma}

\proof
By Lemma~\ref{lemma.PiLambda} we have $0 \le \Pi_\Lambda w \le w$.
Applying $A$ on the left and using $A \Pi_\Lambda = \Pi_\Lambda$,
we obtain (i).
The remark in brackets is obtained by left-multiplying with
any vector $c \ge 0$ ($c \neq 0$) supported on $\Lambda^c$.

(\ref{lemma.A_minus_Pilambda.star1}a,b) are trivial,
and (\ref{lemma.A_minus_Pilambda.star1}c) follows from
$A \Pi_\Lambda = \Pi_\Lambda$ using \reff{eq.Pi2}:
$A - \Pi_\Lambda = A (I - \Pi_\Lambda) = A I_\Lambda (I - \Pi_\Lambda)$.
\reff{lemma.A_minus_Pilambda.star2a} is trivial.
The hypothesis $A \Pi_\Lambda = \Pi_\Lambda$ yields
$\Pi_\Lambda - A I_{\Lambda^c} = A (\Pi_\Lambda - I_{\Lambda^c})
 = A I_\Lambda \Pi_\Lambda \ge 0$,
which is \reff{lemma.A_minus_Pilambda.star2b}.
Finally, \reff{lemma.A_minus_Pilambda.star2c}
is a rewriting of (i).

The identity in \reff{rouen.r1} is an application of
(\ref{lemma.A_minus_Pilambda.star1}c), and the inequality follows from
the fact that $\| I - \Pi_\Lambda \| _{w\to w} \le 2$
(cf.\ Lemma~\ref{lemma.PiLambda}).

Since $c A I_\Lambda$ is supported on $\Lambda$
while $c (\Pi_\Lambda - A I_{\Lambda^c})$ is supported on $\Lambda^c$,
the equality in \reff{lemma.A_minus_Pilambda.star3}
follows from \reff{lemma.A_minus_Pilambda.star1}.
The inequality follows from \reff{lemma.A_minus_Pilambda.star2c}
together with the nonnegativity of the operators
$\Pi_\Lambda - A I_{\Lambda^c}$ and $A I_\Lambda$.
\qed

We can now compare the ``efficiency of cleaning'' of
two operators $A$ and $B$.

\begin{corollary}[Comparison of cleaners]
   \label{cor.AB_minus_Pilambda}
\FH
Let $\Lambda \subseteq X$,
and let $A,B \ge 0$ be operators absorbed by $\Pi_\Lambda$.
Suppose further that
\begin{eqnarray}
   A I_\Lambda w       & \le &   B I_\Lambda w       \label{rouen.r10}  \\
   A I_{\Lambda^c} w   & \ge &   B I_{\Lambda^c} w   \label{rouen.r11}
\end{eqnarray}
Then
\be\label{rouen.r6}
   \| A - \Pi_\Lambda \|_{w\to w}
   \;\le\;
   \| B - \Pi_\Lambda \|_{w\to w}
\ee
and, for every vector $c \ge 0$,
\be \label{rouen.r7}
   \| c (A - \Pi_\Lambda) \|_w
   \;\le\;
   \| c (B - \Pi_\Lambda) \|_w   \;.
\ee
\end{corollary}

\proof
The vector-norm inequality \reff{rouen.r7}
is an immediate consequence of
the equality in \reff{lemma.A_minus_Pilambda.star3}
together with the hypotheses \reff{rouen.r10}/\reff{rouen.r11}.

For the operator-norm inequality \reff{rouen.r6},
note that the decomposition \reff{lemma.A_minus_Pilambda.star1b} yields
\be
  (A-\Pi_\Lambda)_{xy}  \;=\;
  \cases{ A_{xy}                &  for $y\in\Lambda$   \cr
          \noalign{\vskip 2mm}
          -(\Pi_\Lambda-A)_{xy} &  for $y\in\Lambda^c$ \cr
        }
\ee
Inequalities \reff{lemma.A_minus_Pilambda.star2a} and
\reff{lemma.A_minus_Pilambda.star2b} therefore imply the componentwise
identity
\begin{equation}
  \label{eq:rz.1}
  |A-\Pi_\Lambda| \;=\; A I_\Lambda + (\Pi_\Lambda - A I_{\Lambda^c}) \;.
\end{equation}
Inequality \reff{rouen.r6} is an immediate consequence
of this identity and the hypotheses \reff{rouen.r10}/\reff{rouen.r11}
together with \reff{def.operatornorm_bis}.
\qed

The applications of interest to us follow from the multi-monotonicity bounds
\reff{eq.multi1} and \reff{eq.multi3} and the collapse inequalities
\reff{eq.coll1} and \reff{eq.collrev}:

\begin{corollary}\label{roma.cr1}
\FH
Let $\Lambda \subseteq X$.
\begin{itemize}
\item[(i)]  Suppose that $0 \le g_i \le h_i \le \chi_\Lambda$
for $i=1,\ldots,n$.  Then
\be\label{roma.r30}
\| \beta_{h_1} \cdots \beta_{h_n}  \,-\, \Pi_\Lambda \|_{w\to w}
\;\le\;
\| \beta_{g_1} \cdots \beta_{g_n}  \,-\, \Pi_\Lambda \|_{w\to w}
\ee
and, for each vector $c\ge 0$,
\be\label{roma.r31}
\|c( \beta_{h_1} \cdots \beta_{h_n}  \,-\, \Pi_\Lambda) \|_w
\;\le\;
\|c( \beta_{g_1} \cdots \beta_{g_n}  \,-\, \Pi_\Lambda) \|_w   \;.
\ee

\item[(ii)] Suppose that $0\le h_i\le\chi_\Lambda$
for $i=1,\ldots,N$.  Then, for every choice
of integers $0\le n_1\le \cdots\le n_k\le N$, we have
\be\label{roma.r32}
\| \beta_{h_1} \cdots \beta_{h_N}  \,-\, \Pi_\Lambda \|_{w\to w}
\;\le\;
\| \beta_{1-\prod_{i=1}^{n_1} (1-h_i)}\cdots
\beta_{1-\prod_{i=n_{k-1}+1}^{n_k} (1-h_i)}
\,-\, \Pi_\Lambda \|_{w\to w}
\ee
and, for each vector $c\ge 0$,
\be\label{roma.r33}
\| c(\beta_{h_1} \cdots \beta_{h_N}  \,-\, \Pi_\Lambda) \|_w
\;\le\;
\|c( \beta_{1-\prod_{i=1}^{n_1} (1-h_i)}\cdots
\beta_{1-\prod_{i=n_{k-1}+1}^{n_k} (1-h_i)}
\,-\, \Pi_\Lambda) \|_w  \;.
\ee
\end{itemize}

\end{corollary}


\subsection{Convergence of cleaning operators to $\Pi_\Lambda$}
  \label{sec.operator.convergence}

We are now ready to study the convergence of cleaning operators
$\beta_{h_1} \cdots \beta_{h_n}$ to $\Pi_\Lambda$.
We shall prove the main result in two versions:
a uniform (operator-norm) version,
and a dust-dependent (vector-norm) version.
A central role in these analyses will be played, respectively,
by the quantities
\be
   \rho_\Lambda(\ell)   \;\bydef\;  
\| (I_\Lambda \alpha I_\Lambda)^\ell \| _{w \to w}
\ee
and, for each vector $c \ge 0$,
\be
   \rho_\Lambda(\ell;c)   \;\bydef\;  
\| c (I_\Lambda \alpha I_\Lambda)^\ell \| _w   \;.
\ee
Clearly we have $0 \le \rho_\Lambda(\ell) \le 1$
and $0\le \rho_\Lambda(\ell;c)\le \|c\|_w$.
Note also that both $\rho_\Lambda(\ell)$ and
[since $c \ge 0$] $\rho_\Lambda(\ell;c)$
are increasing functions of $\Lambda$.

For brevity let us denote
\begin{equation}
\beta_{\epsilon\Lambda} \;\bydef\; \beta_{\epsilon\chi_\Lambda}
\label{roma.r1}
\end{equation}
for $0 < \epsilon < 1$ and a set $\Lambda \subseteq X$.

\begin{lemma}\mbox{}\label{roma.le2}
Fix a region $\Lambda \subseteq X$ and a number $\epsilon > 0$.
\begin{itemize}
\item[(i)] If $\rho_\Lambda(\ell)\to 0$ as $\ell\to\infty$, then
\be\label{roma.r22}
\|(I_{\Lambda'} \beta_{\epsilon\Lambda'} I_{\Lambda'})^N\|_{w\to w}
\;\mathop{\longrightarrow}\limits_{N\to\infty}\; 0
\ee
uniformly for all regions $\Lambda'\subseteq\Lambda$.
\item[(ii)] If the nonnegative vector $c \in l^1(w)$ is such that
$\rho_\Lambda(\ell;c)\to 0$ as $\ell\to\infty$, then
\be
\|c (I_{\Lambda'} \beta_{\epsilon\Lambda'} I_{\Lambda'})^N\|_w
\;\mathop{\longrightarrow}\limits_{N\to\infty}\; 0
\ee
uniformly for all regions $\Lambda'\subseteq\Lambda$.
\end{itemize}
\end{lemma}

\proof
The obvious identity
\begin{eqnarray}
(I_{\Lambda'} \beta_{\epsilon\Lambda'} I_{\Lambda'})^N &=&
\Bigl[ (1-\epsilon) I_{\Lambda'} + 
\epsilon I_{\Lambda'}\alpha I_{\Lambda'} \Bigr]^N \nonumber\\
&=& \sum_{\ell=0}^N {N \choose \ell} 
(1-\epsilon)^{N-\ell}\epsilon^\ell 
(I_{\Lambda'}\alpha I_{\Lambda'})^\ell  I_{\Lambda'}
\end{eqnarray}
(valid for $N \ge 1$)
allows us to write
\be\label{roma.r20}
\left.\begin{array}{l}
\|(I_{\Lambda'} \beta_{\epsilon\Lambda'} I_{\Lambda'})^N\|_{w\to w}\\[5pt]
\|c (I_{\Lambda'} \beta_{\epsilon\Lambda'} I_{\Lambda'})^N\|_w
      \end{array}
\right\} \;\le\; E(F_{\Lambda'}(X_N))
\ee
where $X_N$ is a Binomial($N,\epsilon$) random variable and
\begin{equation}
F_{\Lambda'}(\ell) \;=\; \left\{\begin{array}{ll}
\rho_{\Lambda'}(\ell) & \mbox{in case (i)}\\
\rho_{\Lambda'}(\ell;c) & \mbox{in case (ii)}
				\end{array}\right.
\end{equation}
In both cases,
\be
F_{\Lambda'}(\ell) \;\le\; F_\Lambda(\ell) 
\;\mathop{\longrightarrow}\limits_{\ell\to\infty}\; 0\;.
\ee
Denoting
\be
M_K \;\bydef\; \sup_{\ell\ge K} F_\Lambda(\ell)  \;,
\ee
we can decompose
\be \label{roma.r21}
E(F_{\Lambda'}(X_N)) \;\le\; M_0\,P(X_N\le K) + M_K\;,
\ee
where $M_0 < \infty$ [in case (ii) because $c \in l^1(w)$].
By hypothesis, given $\delta>0$ we can choose $K$ so that $M_K\le
\delta/2$.  For such $K$,
\be
P(X_N\le K) \;\le\; (1-\epsilon)^{N-K} N^K
\sum_{\ell=0}^K (1-\epsilon)^{K-\ell}\epsilon^\ell 
 \;\mathop{\longrightarrow}\limits_{N\to\infty}\; 0\;.
\ee
Thus for $N$ large enough the first term on the right-hand-side of
\reff{roma.r21} is also smaller than $\delta/2$.  \qed

We remark that in the situation (i), the inequality
\be
\|(I_{\Lambda'} \beta_{\epsilon\Lambda'} I_{\Lambda'})^{N_1+N_2}\|_{w\to w}
\;\le\;
\|(I_{\Lambda'} \beta_{\epsilon\Lambda'} I_{\Lambda'})^{N_1}\|_{w\to w} \,
\|(I_{\Lambda'} \beta_{\epsilon\Lambda'} I_{\Lambda'})^{N_2}\|_{w\to w}
\ee
implies that the convergence in \reff{roma.r22} is actually
exponentially fast in $N$.

\begin{theorem}[Uniform cleaning]
 \label{thm.main.operators.uniform}
\FH
Consider a region $\Lambda \subseteq X$
and functions $0\le h_i\le \chi_\Lambda$ ($i\ge 1$)
such that
\be\label{roma.r23}
   \sum_{i=1}^\infty h_i(x) \;=\; \infty
     \qquad\hbox{uniformly for } x\in\Lambda\;.
\ee
\begin{itemize}
\item[(i)] If $\rho_\Lambda(\ell) \equiv
               \| (I_\Lambda \alpha I_\Lambda)^\ell \| _{w \to w} \to 0$
as $\ell\to\infty$, then
\be\label{roma.r26}
\| \beta_{h_1} \cdots \beta_{h_N}  \,-\, \Pi_\Lambda \|_{w\to w}
\;\mathop{\longrightarrow}\limits_{N\to\infty}\; 0\;.
\ee
\item[(ii)] If the nonnegative vector $c \in l^1(w)$ is such that
$\rho_\Lambda(\ell;c) \equiv \| c (I_\Lambda \alpha I_\Lambda)^\ell \|_w \to 0$
as $\ell\to\infty$, then
\be
\|c ( \beta_{h_1} \cdots \beta_{h_N}  \,-\, \Pi_\Lambda)\|_w
\;\mathop{\longrightarrow}\limits_{N\to\infty}\; 0\;.
\ee
\end{itemize}
\end{theorem}

We remark that if $\Lambda$ is a {\em finite}\/ set,
then the hypothesis \reff{roma.r23} is equivalent to the
apparently weaker hypothesis that
$\sum_{i=1}^\infty h_i(x) = \infty$ for all $x\in\Lambda$
[cf.\ \reff{roma.r43} below].
But if  $\Lambda$ is an {\em infinite}\/ set,
then \reff{roma.r23} is stronger.

\proofof{Theorem~\ref{thm.main.operators.uniform}}
By Proposition~\ref{prop.convergence-to-balayage}
and the fact that $\| I - \Pi_\Lambda \| _{w\to w} \le 2$
[or alternatively by  Lemma~\ref{lemma.A_minus_Pilambda}(iii)],
it is enough to show that
\begin{subeqnarray}
\|  \beta_{h_1} \cdots \beta_{h_{N}} I_\Lambda \|_{w \to w}
 & \mathop{\longrightarrow}\limits_{N\to\infty} & 0
 \qquad\hbox{in case (i)}  \\[2mm]
\|c  \beta_{h_1} \cdots \beta_{h_{N}} I_\Lambda \|_w
 & \mathop{\longrightarrow}\limits_{N\to\infty} & 0
 \qquad\hbox{in case (ii)}
\end{subeqnarray}
%
Now, since the $h_i$ are bounded, \reff{roma.r23}
is equivalent to the existence of $\delta>0$ and a sequence
$0 = n_0 < n_1 < n_2 < \ldots$ satisfying
\be\label{roma.r23a}
\sum_{i=n_{j-1}+1} ^{n_j} h_i \;\ge\; \delta\,\chi_\Lambda
\ee
for all $j$.
This, in turn, is equivalent to the existence of
$\epsilon>0$ such that 
\be\label{roma.r24}
1-\prod_{i=n_{j-1}+1}^{n_j}(1-h_i) \;\ge\;\epsilon\chi_\Lambda
\ee
for every $j$
(indeed, we can set $\epsilon = 1 - e^{-\delta}$).
Therefore, by Corollary~\ref{cor.multi-monotonicity+collapse},
if $N \ge n_k$ we have
\begin{subeqnarray}
\| \beta_{h_1} \cdots \beta_{h_N}  I_\Lambda \|_{w\to w}
  & \le & \|\beta_{\epsilon\Lambda}^k I_\Lambda \|_{w\to w}  \\[2mm]
\| c \beta_{h_1} \cdots \beta_{h_N}  I_\Lambda \|_{w}
  & \le & \| c \beta_{\epsilon\Lambda}^k I_\Lambda \|_{w}
\end{subeqnarray}
The theorem then follows from
Lemmas~\ref{lemma_betagi} and \ref{roma.le2}.
\qed

{\bf Remark.}
The standard ``cleaning'' proof
of the Dobrushin uniqueness theorem
\cite{Vasershtein_69,Lanford_73,Follmer_82,Georgii_88,Simon_93}
proves a very special case of Theorem~\ref{thm.main.operators.uniform}(ii):
namely, one assumes the very strong hypothesis $\alpha w \le (1-\epsilon)w$
for some $\epsilon > 0$,
and one takes $h_i = \delta_{x_i}$,
where $x_1,x_2,\ldots$ is a sequence that visits each site of $\Lambda$
infinitely many times.
The correlations between $\Lambda$ and $\Lambda^c$
can then be bounded in terms of $\Pi_\Lambda$
(see \cite{Follmer_82,Georgii_88,Simon_93} for variants of this idea).

Let us also remark that, in the application to the Dobrushin uniqueness
theorem, it appears to be necessary to take $w = {\bf 1}$.
This choice plays no role in the ``cleaning'' proof itself,
but plays a role in the final step of the argument,
where the total oscillation of a function of many variables
is bounded by the sum of its single-variable oscillations
--- {\em not}\/ a weighted sum.
\qed

\bigskip

The non-uniform (dust-dependent) version of the previous result relies
on the following decomposition:

\begin{lemma}
\label{roma.pr1}
Fix $\epsilon>0$ and a set $\Lambda \subseteq X$.
Then, for all integers $N \ge K \ge 0$,
\begin{equation}\label{roma.r2}
\beta^N_{\epsilon\Lambda} \;=\; 
\beta^K_{\epsilon\Lambda} I_{\Lambda^{\rm c}} +
(I_\Lambda\beta_{\epsilon\Lambda} I_\Lambda)^K I_\Lambda
\beta^{N-K}_{\epsilon\Lambda} \;.
\end{equation}
\end{lemma}

\proof
By \reff{roma.r10} we have
\begin{equation}\label{roma.r3}
I_{\Lambda^{\rm c}} \beta^N_{\epsilon\Lambda} \;=\; 
I_{\Lambda^{\rm c}}
\end{equation}
and hence
\begin{equation}\label{roma.r4} 
\beta^N_{\epsilon\Lambda}\;=\; I_{\Lambda^{\rm c}}
+ I_\Lambda \beta^N_{\epsilon\Lambda}\;.
\end{equation}
It follows that
\begin{eqnarray}
\beta^N_{\epsilon\Lambda} &=& 
\beta^K_{\epsilon\Lambda}\beta^{N-K}_{\epsilon\Lambda}\nonumber\\[1mm]
&=&
\beta^K_{\epsilon\Lambda} (I_{\Lambda^{\rm c}}
              + I_\Lambda \beta^{N-K}_{\epsilon\Lambda})
       \qquad\qquad\qquad \hbox{[by \reff{roma.r4}]}  \nonumber \\[1mm]
&=&
\beta^K_{\epsilon\Lambda} I_{\Lambda^{\rm c}} +
(I_\Lambda\beta_{\epsilon\Lambda} I_\Lambda)^K I_\Lambda
\beta^{N-K}_{\epsilon\Lambda}
       \qquad \hbox{[by \reff{eq.betagi.1}]}
\end{eqnarray}
\qed

\begin{theorem}[Dust-dependent cleaning]
 \label{thm.main.operators}
\FH
Consider a region $\Lambda \subseteq X$
and functions $0\le h_i\le \chi_\Lambda$ ($i\ge 1$) such that 
\be\label{roma.r43}
\sum_{i=1}^\infty h_i(x) \;=\; \infty\qquad\hbox{for all } x\in\Lambda\;.
\ee
If the nonnegative vector $c \in l^1(w)$ is such that
$\rho_\Lambda(\ell;c) \equiv \| c (I_\Lambda \alpha I_\Lambda)^\ell \|_w \to 0$
as $\ell\to\infty$, then
\be
\|c ( \beta_{h_1} \cdots \beta_{h_{N}}  \,-\, \Pi_\Lambda)\|_w
\;\mathop{\longrightarrow}\limits_{N\to\infty}\; 0\;.
  \label{eq.thm.main.operators}
\ee
\end{theorem}

\proof
As in the proof of Theorem \ref{thm.main.operators.uniform},
it is enough to show that
\be \label{roma.r54}
\|c  \beta_{h_1} \cdots \beta_{h_{N}} I_\Lambda \|_w
\;\mathop{\longrightarrow}\limits_{N\to\infty}\; 0\;.
\ee
Fix $\delta > 0$ and $0 < \epsilon < 1$.
By Lemma~\ref{roma.le2} we can choose $K$ so that
\be
   \|c(I_{\Lambda'}\beta_{\epsilon\Lambda'} I_{\Lambda'})^K\|_w
   \;\le\;  \delta/2
 \label{eq.roberto.3} 
\ee
uniformly for all regions $\Lambda' \subseteq \Lambda$.
On the other hand, since $c \in l^1(w)$, we can choose
a finite set $\Lambda' \subseteq \Lambda$ so that
\be\label{roma.r44}
\|c(I_\Lambda\alpha I_\Lambda)^\ell I_{\Lambda\setminus\Lambda'} \|_w
\;\le\;  \delta/[2(K+1)]
 \label{eq.roberto.4} 
\ee
for $\ell = 0,1,\ldots,K$.
Then hypothesis \reff{roma.r43} guarantees that
there exists a sequence of integers $0 = n_0 < n_1 < n_2 < \ldots$ such that 
\be
1-\prod_{i=n_{j-1}+1}^{n_j}(1-h_i) \;\ge\;\epsilon\chi_{\Lambda'}
\ee
for all $j$.
Therefore, by Corollary~\ref{cor.multi-monotonicity+collapse},
if $N \ge n_k$ we have
\be
 \|c \beta_{h_1} \cdots \beta_{h_{N}} I_\Lambda \|_w
\;\le\; \| c \beta^k_{\epsilon\Lambda'} I_\Lambda \|_w\;.
 \label{eq.roberto.5}
\ee
Now, by the decomposition \reff{roma.r2}, whenever $k \ge K$ we have
\begin{equation}
\beta^k_{\epsilon\Lambda'} I_\Lambda \;=\; 
\beta^K_{\epsilon\Lambda'} I_{\Lambda\setminus\Lambda'} +
(I_{\Lambda'}\beta_{\epsilon\Lambda'} I_{\Lambda'})^K  I_{\Lambda'}
\beta^{k-K}_{\epsilon\Lambda'} I_\Lambda \;.
 \label{eq.roberto.1}
\end{equation}  
Applying this to the vector $w$, we deduce the vector inequalities
\begin{eqnarray}
\beta^k_{\epsilon\Lambda'} I_\Lambda w
 &\le &
\beta^K_{\epsilon\Lambda'} I_{\Lambda\setminus\Lambda'} w +
(I_{\Lambda'}\beta_{\epsilon\Lambda'} I_{\Lambda'})^K I_{\Lambda'} w
   \qquad\qquad  \hbox{[by \reff{eq.norm.beta}]}  \nonumber\\[1mm]
 &\le &
\beta^K_{\Lambda'} I_{\Lambda\setminus\Lambda'} w +
(I_{\Lambda'}\beta_{\epsilon\Lambda'} I_{\Lambda'})^K I_{\Lambda'} w
   \qquad\qquad\; \hbox{[by \reff{eq.multi3}]}  \nonumber\\[1mm]
 & = & \sum_{\ell =0}^{K} (I_{\Lambda'} \alpha)^\ell
I_{\Lambda\setminus\Lambda'} w +
(I_{\Lambda'}\beta_{\epsilon\Lambda'} I_{\Lambda'})^K I_{\Lambda'} w
     \nonumber\\[1mm]
&\le & \sum_{\ell =0}^{K} (I_\Lambda\alpha I_\Lambda)^\ell
I_{\Lambda\setminus\Lambda'} w +
(I_{\Lambda'}\beta_{\epsilon\Lambda'} I_{\Lambda'})^K I_{\Lambda'} w\;.
 \label{eq.roberto.2}
\end{eqnarray}
Thus, 
\begin{eqnarray}
\|c \beta^k_{\epsilon\Lambda'} I_\Lambda \|_w
  & \le &
\sum_{\ell =0}^{K} \|c (I_\Lambda\alpha I_\Lambda)^\ell
I_{\Lambda\setminus\Lambda'}\|_w  \,+\,
\|c(I_{\Lambda'}\beta_{\epsilon\Lambda'} I_{\Lambda'})^K\|_w  \nonumber\\[1mm]
  & \le &  \delta  \;.
\end{eqnarray}
\qed

Let us observe that, in the absence of some unformity hypothesis
on the $h_i$ [like \reff{roma.r23}],
the convergence in \reff{eq.thm.main.operators}
can be {\em arbitrarily slow}\/,
even if $\rho_\Lambda(\ell;c) \to 0$ arbitrarily rapidly:

\begin{example}
 \label{example.thm.main.operators}
\rm
Let $\Lambda$ be countably infinite (say, $\Lambda= \{1,2,3,\ldots\}$);
let $w = {\bf 1}$ and choose any $c>0$ with $c \in l^1$.
Let us consider the {\em best possible case}\/ for cleaning,
namely $\alpha=0$, so that $\rho_\Lambda(\ell;c) = 0$ for all $\ell \ge 1$.
Then $\| c (\beta_h^N - \Pi_\Lambda) \|_w =
      \sum_{i=1}^\infty c_i [1-h(i)]^N \defby \epsilon_N$.
It is a fairly simple analysis exercise to show that the $(h_i)_{i \ge 1}$
can be chosen so that $\epsilon_N$ decays more slowly with $N$
than any specified convergent-to-zero sequence
$(\delta_N)_{N \ge 0}$.\footnote{
   Let $(c_i)_{i \ge 1}$ be any strictly positive sequence,
   and let $(\delta_N)_{N \ge 0}$ be any sequence of nonnegative
   numbers converging to zero.
   We claim that one can choose a sequence $(h_i)_{i \ge 1}$ 
   of numbers in $(0,1]$ such that
   $\epsilon_N \bydef \sum_{i=1}^\infty c_i (1-h_i)^N \ge \delta_N$
   for all but finitely many $N$.
   {\sc Proof.}
   Choose $N_1$ such that $\delta_N \le c_1/2$ for all $N \ge N_1$;
   and for $i \ge 2$, inductively choose $N_i > N_{i-1}$ such that
   $\delta_N \le c_i/2$ for all $N \ge N_i$.
   Then, for each $i \ge 1$, choose $h_i$ small enough so that
   $(1-h_i)^{N_{i+1}} \ge 1/2$.
   It follows that, for $N_i \le N \le N_{i+1}$ we have
   $$ \delta_N  \;\le\;  c_i/2  \;\le\;  c_i (1-h_i)^{N_{i+1}}
      \;\le\;  c_i (1-h_i)^N  \;\le\;  \epsilon_N  \;.$$
   Since every $N \ge N_1$ belongs to some such interval, we are done.
   \qed
}
\end{example}

\subsection{Some further identities and inequalities}
   \label{subsec.matrix.idineq}

Let us now prove a beautiful identity for the sum of a 
geometric series $\sum_{N=0}^\infty (I_\Lambda \beta_h I_\Lambda)^N$.
This identity will play a central role in the next subsection
in the proof of the converse theorem on cleaning
(Theorem~\ref{thm.matrix.converse}).

\begin{lemma}
   \label{lemma.12Apr07.identity}
Let $\Lambda \subseteq X$, and let $h \colon\, X \to [0,1]$
be strictly positive on $\Lambda$ and zero outside $\Lambda$.
(In other words, $\supp h = \Lambda$.)  Then
\be
   \sum_{N=0}^\infty (I_\Lambda \beta_h I_\Lambda)^N I_h
   \;=\;
   \sum_{k=0}^\infty (I_\Lambda \alpha I_\Lambda)^k I_\Lambda
   \;.
 \label{eq.lemma.12Apr07.identity}
\ee
\end{lemma}

\proof
Everything occurs within $\Lambda$,
so for notational simplicity let us suppose that $\Lambda = X$.
Write $\beta_h = I_{1-h} + I_h \alpha$, expand out the $N$th power,
and sum over $N$.  We get a sum over all finite sequences
(including the empty sequence) of factors $I_{1-h}$ and $I_h \alpha$.
Now let us treat the matrix elements of $\alpha$ as
noncommuting indeterminates and extract the coefficient of a monomial
$\alpha_{x_0 x_1} \alpha_{x_1 x_2} \cdots \alpha_{x_{k-1} x_k}$
with $k \ge 0$
(it is easy to see that these are the only monomials that arise).
To the left of $\alpha_{x_0 x_1}$ we have an arbitrary number (including zero)
of factors $1-h(x_0)$ followed by one factor $h(x_0)$:
this gives $\sum_{n=0}^\infty [1-h(x_0)]^n h(x_0) = 1$
since $0 < h(x_0) \le 1$.
Likewise to the immediate left of each $\alpha_{x_{i-1} x_i}$.
Finally, to the right of $\alpha_{x_{k-1} x_k}$
we have an arbitrary number (including zero) of factors $1-h(x_k)$:
this gives $\sum_{n=0}^\infty [1-h(x_k)]^n = 1/h(x_k)$,
and this factor is cancelled by the $I_h$
on the left-hand side of \reff{eq.lemma.12Apr07.identity}.
So each monomial
$\alpha_{x_0 x_1} \alpha_{x_1 x_2} \cdots \alpha_{x_{k-1} x_k}$
gets a coefficient 1, which corresponds exactly
to the right-hand side of \reff{eq.lemma.12Apr07.identity}.
\qed

{\bf Important Remark.}
By treating the matrix elements of $\alpha$
as noncommuting indeterminates, we are in essence using the tree formalism
that will be described in detail in Section~\ref{sec4}.

\bigskip

If we use a product $\beta_{f_1} \cdots \beta_{f_m}$ in place of
the single cleaning operator $\beta_h$,
then we can obtain an {\em inequality}\/ in place of the identity
\reff{eq.lemma.12Apr07.identity}:

\begin{lemma}
   \label{lemma.12Apr07.inequality}
Let $\Lambda \subseteq X$, and let $f_1,\ldots,f_m \colon\, X \to [0,1]$
be supported on $\Lambda$.  Define $h \bydef 1 - \prod\limits_{i=1}^m (1-f_i)$.
Then
\be
   \sum_{k=0}^\infty (I_\Lambda \beta_{f_1} \cdots \beta_{f_m} I_\Lambda)^k I_h
   \;\le\;
   \sum_{k=0}^\infty (I_\Lambda \alpha I_\Lambda)^k I_\Lambda
 \label{eq.lemma.12Apr07.inequality}
\ee
provided that the right-hand side is elementwise finite.
\end{lemma}

\proof
As in the previous lemma, everything here occurs within $\Lambda$,
so we can assume for notational simplicity that $\Lambda=X$.
We apply Corollary~\ref{cor_telescoping} with
$g_1 = \ldots = g_m = 0$ and $h_i = f_i$, to obtain
\begin{subeqnarray}
   I - \beta_{f_1} \cdots \beta_{f_m}
   & = &
   \Biggl( \sum_{i=1}^m  I_{f_i} \beta^*_{f_{i+1}} \cdots \beta^*_{f_m} \Biggr)
      \,(I-\alpha)
        \\[2mm]
   & \bydef &  P (I-\alpha)
 \label{eq.18Apr07.niceversion}
\end{subeqnarray}
where $P$ is a sum of products of the operators
$\alpha$, $I_{f_j}$ and $I_{1-f_j}$ $(1 \le j \le m)$.
Furthermore, the term in $P$ containing no factors of $\alpha$ is
the operator of multiplication by
\be
   \sum_{i=1}^m f_i \prod_{j=i+1}^m (1-f_j)
   \;=\;
   1 - \prod\limits_{i=1}^m (1-f_i)
   \;=\;
   h \;.
 \label{eq.18Apr07.Ih}
\ee
Since the other terms are nonnegative, we have $P \ge I_h$.

Let us now abbreviate $B = \beta_{f_1} \cdots \beta_{f_m}$.
Under the assumption that $\sum_{k=0}^\infty \alpha^k$
is elementwise finite, we have
\be
   (I-\alpha)
   \, \Biggl( \sum\limits_{k=0}^\infty \alpha^k \Biggr)
   \;=\;
   I
   \;.
\ee
Therefore, right-multiplying \reff{eq.18Apr07.niceversion}
by $\sum_{k=0}^\infty \alpha^k$ yields
\be
   I_h \;\le\; P
   \;=\;
   (I-B) \sum\limits_{k=0}^\infty \alpha^k
   \;.
\ee
We now left-multiply this inequality by $B^k$
and sum from $k=0$ to $N$:
since $B^k (I-B) = B^k - B^{k+1}$, the sum telescopes and we have
\be
   \sum_{k=0}^N B^k I_h
   \;\le\;
   (I - B^{N+1}) \sum\limits_{k=0}^\infty \alpha^k
   \;\le\;
   \sum\limits_{k=0}^\infty \alpha^k
   \;.
\ee
Taking $N \to\infty$ gives the result.
\qed

Lemma~\ref{lemma.12Apr07.inequality}
is a special case of a result to be proven in Section~\ref{sec4.9}
using the tree formalism 
[see Lemmas~\ref{lemma.12Apr07.cloud}(b) and \ref{lemma.16Apr07}(a)].
Indeed, the mysterious operator $P$ in \reff{eq.18Apr07.niceversion}
will correspond to the cloud $\mu$ in Lemma~\ref{lemma.16Apr07}(a).
Furthermore, in the tree context the summability condition on
$\sum_{k=0}^\infty \alpha^k$ can be removed.

\bigskip

Finally, we have a reverse inequality:

\begin{lemma}
   \label{lemma.18Apr07.hardest}
Let $\Lambda \subseteq X$, and let $f_1,\ldots,f_m \colon\, X \to [0,1]$
be supported on $\Lambda$.
Then
\be
   \sum_{k=0}^\infty (I_\Lambda \alpha I_\Lambda)^k
   \;\le\;
   \left( \sum_{n=0}^{\infty} (I_\Lambda \beta_{f_1} \cdots \beta_{f_m}
                               I_\Lambda)^n \right)
   \left( \sum_{k=0}^{m-1} (I_\Lambda \alpha I_\Lambda)^k \right)
 \label{eq.lemma.18Apr07.hardest}
\ee
provided that
$\sum_{n=0}^{\infty} (I_\Lambda \beta_{f_1} \cdots \beta_{f_m} I_\Lambda)^n$
is elementwise finite.
\end{lemma}

The main tool in the proof is the following bound,
which uses only the nonnegativity of the matrix $\alpha$:

\begin{lemma}
   \label{lemma.18Apr07.betastarsum}
Let $f_1,\ldots,f_m \colon\, X \to [0,1]$.
Then
\be
   0
   \;\le\;
   \sum_{i=1}^m  I_{f_i} \beta^*_{f_{i+1}} \cdots \beta^*_{f_m}
   \;\le\;
   \sum_{k=0}^{m-1} \alpha^k
   \;.
 \label{eq.lemma.18Apr07.betastarsum}
\ee
\end{lemma}

\proof
Write $\beta^*_f = I_{1-f} + \alpha I_f$
and expand out the left-hand side of \reff{eq.lemma.18Apr07.betastarsum}.
Let us once again treat the matrix elements of $\alpha$ as
noncommuting indeterminates and extract the coefficient of a monomial
$\alpha_{x_0 x_1} \alpha_{x_1 x_2} \cdots \alpha_{x_{k-1} x_k}$
with $0 \le k \le m-1$
(it is easy to see that these are the only monomials that arise).
We need to show that each such coefficient is $\le 1$
(the nonnegativity is obvious).
We have already computed in \reff{eq.18Apr07.Ih}
the term with no powers of $\alpha$:
it is $1 - \prod_{i=1}^m [1-f_i(x_0)] \le 1$.
Now suppose that there are $k$ $\alpha$'s,
occurring at positions $j_1,\ldots,j_k$
with $i+1 \le j_1 < j_2 < \ldots < j_k \le m$.
The coefficient of such a term will be
\begin{eqnarray}
   [f_i (1-f_{i+1}) \cdots (1-f_{j_1-1})](x_0)
   & \;\times\; &
   [f_{j_1} (1-f_{{j_1}+1}) \cdots (1-f_{j_2-1})](x_1)  \qquad
         \nonumber   \\
   \;\times\;\;\;\cdots
   & \;\times\; &   [f_{j_k} (1-f_{{j_k}+1}) \cdots (1-f_m)](x_k)
\end{eqnarray}
and we then need to sum over all choices of indices
$1 \le i < j_1 < j_2 < \ldots < j_k \le m$.
First fix $j_1, \ldots, j_k$ and sum over $i$: one gets
\be
   \sum_{i=1}^{j_1-1} [f_i (1-f_{i+1}) \cdots (1-f_{j_1-1})](x_0)
   \;=\;
   1 \,-\, \prod_{i=1}^{j_1-1} [1-f_i(x_0)]
   \;\le\;
   1
\ee
[just as in \reff{eq.18Apr07.Ih}].
We then sum over $j_1$ and so forth,
each time bounding the sum by 1.
\qed

\proofof{Lemma~\ref{lemma.18Apr07.hardest}}
Once again, we can assume that $\Lambda = X$.
Let us abbreviate $B = \beta_{f_1} \cdots \beta_{f_m}$.
{}From \reff{eq.18Apr07.niceversion} we have
\be
   I - B  \;=\;  P (I-\alpha)
 \label{eq.proof.lemma.18Apr07.hardest}
\ee
where Lemma~\ref{lemma.18Apr07.betastarsum} gives
\be
   0  \;\le\; P  \;\le\;  \sum_{k=0}^{m-1} \alpha^k
   \;.
\ee
Under the assumption that $\sum_{n=0}^\infty B^n$
is elementwise finite, we have
\be
   \Biggl( \sum\limits_{n=0}^\infty B^n \Biggr)
   \, (I-B)
   \;=\;
   I
   \;,
\ee
so we can left-multiply \reff{eq.proof.lemma.18Apr07.hardest} by
$\sum_{n=0}^\infty B^n$ to obtain
\be
   I   \;=\;  \Biggl( \sum\limits_{n=0}^\infty B^n \Biggr) \, P (I-\alpha)
   \;.
\ee
Now right-multiply this inequality by $\alpha^k$
and sum from $k=0$ to $N$:  the sum on the right telescopes and we obtain
\begin{subeqnarray}
   \sum\limits_{k=0}^N \alpha^k
   & = &
   \Biggl( \sum\limits_{n=0}^\infty B^n \Biggr) \, P (I-\alpha^{N+1})
        \\[2mm]
   & \le &
   \Biggl( \sum\limits_{n=0}^\infty B^n \Biggr) \, P
        \\[2mm]
   & \le &
   \Biggl( \sum\limits_{n=0}^\infty B^n \Biggr) \, 
   \Biggl( \sum\limits_{k=0}^{m-1} \alpha^k \Biggr)
      \;.
\end{subeqnarray}
Taking $N \to\infty$ gives the result.
\qed


Lemma~\ref{lemma.18Apr07.hardest} is a special case
of a result to be proven in Section~\ref{sec4.9}
using the tree formalism [see Lemma~\ref{lemma.16Apr07}(c,d)],
where moreover the summability condition on
$\sum_{n=0}^{\infty} (I_\Lambda \beta_{f_1} \cdots \beta_{f_m} I_\Lambda)^n$
can be removed.
The pair of inequalities
\reff{eq.lemma.12Apr07.inequality}/\reff{eq.lemma.18Apr07.hardest}
will play a crucial role in the proof of the converse theorem on cleaning
(Theorem~\ref{thm.matrix.converse}),
just as their tree generalizations will do in the proof of the strong form
of this result (Theorem~\ref{thm.converse.2.NEW}).

\subsection{Converse results}   \label{subsec.matrix.converse}

{\em In this subsection we do {\bf not} assume the Fundamental Hypothesis.}\/
Rather, our goal is to study what happens in case
the Fundamental Hypothesis fails.
Here we are entitled to use the algebraic {\em identities}\/
that were proven in the preceding subsections,
since such identities are valid irrespective of the Fundamental Hypothesis.
Furthermore, we are entitled to use those
{\em inequalities that do not refer to the vector $w$}\/,
since they are based only on the nonnegativity
of the matrix elements of $\alpha$.
But we must be very careful to avoid using
any result that depends on the Fundamental Hypothesis.



Our main result is the following:

\begin{theorem}
   \label{thm.matrix.converse}
Let $X$ be a finite or countably infinite set,
let $\Lambda \subseteq X$, and let $c \ge 0$ and $w \ge 0$
be vectors that are strictly positive on $\Lambda$.
Consider the following conditions on a matrix $\alpha$:
\begin{itemize}
  \item[(a)]  $\sum_{k=0}^\infty c (I_\Lambda \alpha I_\Lambda)^k w < \infty$.
  \item[(b)]  For all $h \colon\, X \to [0,1]$ with $\supp h = \Lambda$
      such that $h \ge \epsilon \chi_\Lambda$ for some $\epsilon > 0$,
      we have
      $\sum_{k=0}^\infty c (I_\Lambda \beta_h I_\Lambda)^k w < \infty$.
  \item[(b${}'$)]  There exists $h \colon\, X \to [0,1]$
      with $\supp h = \Lambda$ such that
      $\sum_{k=0}^\infty c (I_\Lambda \beta_h I_\Lambda)^k w < \infty$.
  \item[(c)]  For every finite sequence $f_1,\ldots,f_m$
     of functions $X \to [0,1]$ with $\supp(f_i) \subseteq \Lambda$
     such that $\sum_i f_i \ge \epsilon \chi_\Lambda$
     for some $\epsilon > 0$, we have
     $\sum_{k=0}^\infty c (I_\Lambda \beta_{f_1} \cdots
                                     \beta_{f_m} I_\Lambda)^k w < \infty$.
  \item[(c${}'$)]  There exists a finite sequence $f_1,\ldots,f_m$
     of functions $X \to [0,1]$ with $\supp(f_i) \subseteq \Lambda$
     such that
     $\sum_{k=0}^\infty c (I_\Lambda \beta_{f_1} \cdots
                                     \beta_{f_m} I_\Lambda)^k w < \infty$.
\end{itemize}
Then (a)--(c) are all equivalent and imply (c${}'$);
and for matrices $\alpha$ satisfying the additional hypothesis
\be
  \hbox{There exists a constant $C < \infty$ such that
        $I_\Lambda \alpha I_\Lambda w \le C w$}
 \label{hyp.converse.2.matrix}
\ee
all five conditions are equivalent.
\end{theorem}

\proof
(b) $\implies$ (b${}'$) is trivial,
while (b${}'$) $\implies$ (a) $\implies$ (b)
are immediate consequences of Lemma~\ref{lemma.12Apr07.identity}.

(a) $\implies$ (c) is an immediate consequence of
Lemma~\ref{lemma.12Apr07.inequality},
and (c) $\implies$ (b) is trivial, as is
(c) $\implies$ (c${}'$).

Finally, Lemma~\ref{lemma.18Apr07.hardest} entails
(c${}'$) $\implies$ (a) under the hypothesis \reff{hyp.converse.2.matrix}.
\qed

\medskip

Please note that the hypothesis \reff{hyp.converse.2.matrix} is automatic
whenever $\Lambda$ is a {\em finite}\/ set.
On the other hand, the following two examples show that,
when $\Lambda$ is infinite, hypothesis \reff{hyp.converse.2.matrix}
cannot be dispensed with in proving that (c${}'$) $\implies$ (a):

\begin{example}
    \label{example_thm.converse.2.NEW_1}
\rm
Let $X = \{0,1,2,3,\ldots\}$;
set $\alpha_{0j} = 1$ for all $j \ge 1$,
and set all other matrix elements of $\alpha$ to 0.
Set $w = {\bf 1}$ and let $c$ be any strictly positive vector in $l^1$.
Now let $\Lambda$ be any infinite subset of $X$ containing 0.
We have
\be
   c (I_\Lambda \alpha I_\Lambda)^k w  \;=\;
   \cases{  \|c\|_1 \in (0,\infty)   & for $k=0$  \cr
            +\infty                  & for $k=1$  \cr              
            0                        & for $k \ge 2$ \cr
         }
 \label{eq.example_thm.converse.2.NEW_1}
\ee
so that $\sum_{k=0}^\infty c (I_\Lambda \alpha I_\Lambda)^k w = +\infty$.
On the other hand, if we take $f_1 = \ldots = f_m = \chi_\Lambda$
for any $m \ge 2$, we have
$\sum_{k=0}^\infty c (I_\Lambda \beta_{f_1} \cdots \beta_{f_m} I_\Lambda)^k w
 = \sum_{k=0}^\infty c (I_\Lambda \alpha I_\Lambda)^{mk} w = \|c\|_1 < \infty$.
Indeed, we have $\beta_\Lambda^\ell I_\Lambda = 0$ for all $\ell \ge 2$.
\qed
\end{example}

\begin{example}
    \label{example_thm.converse.2.NEW_1.bis}
\rm
In the preceding example, one of the components of the vector
$I_\Lambda \alpha I_\Lambda w$ was $+\infty$.
Here is a variant in which
the vector $I_\Lambda \alpha I_\Lambda w$ is pointwise finite
but is not bounded by any multiple of $w$.
Take $\Lambda = X = \{x_1,x_2,x_3,\ldots\} \cup \{y_1,y_2,y_3,\ldots\}$;
set $\alpha_{x_i y_j} = 1$ if $j \le i$,
and set all other matrix elements of $\alpha$ to 0.
Set $w = {\bf 1}$ and $c_{x_i} = c_{y_i} = 1/i^2$.
Then $(\alpha w)_{x_i} = i$, $(\alpha w)_{y_i} =0$
and $\alpha^k w = 0$ for $k \ge 2$;
so \reff{eq.example_thm.converse.2.NEW_1} again holds
and the same choice of $f_1,\ldots,f_m$ provides a counterexample.
\qed
\end{example}

We will return to these questions in Section~\ref{sec.converse},
where we will prove a significant extension of
Theorem~\ref{thm.matrix.converse}
(see Theorem~\ref{thm.converse.2.NEW}).

\section{Tree approach to balayage}   \label{sec4}

In this section we introduce an alternate approach to studying
the algebra of operators generated by $\alpha$
and the multiplication operators $I_f$,
which brings out more clearly its underlying combinatorial structure
and which permits a far-reaching generalization of the results
obtained in the preceding section.
This approach is based on considering the tree of finite sequences
of elements of $X$.\footnote{
   We use the word ``tree'' even though it turns out to be convenient,
   in our approach, to suppress the root of the tree
   (that is, the empty sequence).
}
Another way of phrasing matters is to say that we are working in
the algebra of formal power series in noncommuting indeterminates
$\{\alpha_{xy}\}_{x,y \in X}$ subject to the relations
$\alpha_{xy} \alpha_{uv} = 0$ whenever $y \neq u$.

The plan of this section is as follows:
In Section~\ref{sec4.1} we introduce the formalism of
``markers'' and ``clouds'',
and we analyze its relation with the operator formalism of Section~\ref{sec3}.
In particular, we introduce the key operation of convolution of clouds,
which corresponds to multiplication of operators.
In Section~\ref{sec4.2} we define clouds to represent each of
the special operators $\alpha$, $I_f$, $\beta_f$ and $\Pi_\Lambda$.
In Section~\ref{sec4.3} we introduce a very important partial ordering
$\audessus$ on clouds, which formalizes (roughly speaking)
the comparison of operators by their ``efficiency of cleaning''.
We study the circumstances under which
the partial ordering $\audessus$ is preserved by
convolution from the left or the right,
and we introduce several important subclasses of clouds
($\scrb \subsetneq \scrp \subsetneq \scrr \subsetneq \scrs$).
In Section~\ref{sec4.4} we prove a fundamental comparison inequality,
which substantiates our assertion that the
partial ordering $\audessus$ is related to the efficiency of cleaning.
In Section~\ref{sec4.5} we introduce the notion of a cloud being
``carried'' by a subset $\Lambda$,
and in Section~\ref{sec4.6} we introduce the stricter notion of being
``$\Lambda$-regular''.
In Section~\ref{sec4.7} we show that the cloud $\pi_\Lambda$
(which is associated to the balayage operator $\Pi_\Lambda$)
plays a special role among $\Lambda$-regular clouds,
by virtue of its minimality with respect to $\audessus$.
In Section~\ref{sec4.8} we put all these tools together,
and study the convergence of cleaning operators
$\beta_{h_1} \cdots \beta_{h_n}$ to $\Pi_\Lambda$.
In the cloud context we can shed additional light on this convergence,
by distinguishing convergence of clouds from convergence of the
corresponding operators.
In Section~\ref{sec4.9} we prove some further identities and inequalities
that will play a crucial role in the converse results
of Section~\ref{sec.converse}.

\subsection{Markers, clouds and operators}  \label{sec4.1}

A nonempty finite sequence $\eta = (x_0,x_1,\ldots,x_k)$
of elements of $X$ will be called a {\em marker}\/.
We denote by $X^{[\infty]} = \bigcup_{k=0}^\infty X^{k+1}$
the set of all markers.
Given a marker $\eta = (x_0,x_1,\ldots,x_k)$, we define
\begin{eqnarray}
   \hbox{level}(\eta)  & \bydef &  k   \\[1mm]
   \hbox{first}(\eta)  & \bydef &  x_0 \\[1mm]
   \hbox{last}(\eta)   & \bydef &  x_k \\[1mm]
   \eta_i^j            & \bydef &  (x_i,x_{i+1},\ldots,x_j)
\end{eqnarray}
An {\em ancestor}\/ (or {\em prefix}\/) of $\eta$ is
any one of the markers $\eta_0^j$ ($0 \le j \le k$);
we write $\eta' \ancestor \eta$ to denote that
$\eta'$ is an ancestor of $\eta$.
We write $\eta' \ancneq \eta$ to denote that
$\eta' \anc \eta$ and $\eta' \neq \eta$.
A {\em child}\/ of $\eta$ is any marker
of the form $(x_0,x_1,\ldots,x_k,x_{k+1})$ for some $x_{k+1} \in X$.
A {\em suffix}\/ of $\eta$ is
any one of the markers $\eta_j^k$ ($0 \le j \le k$).
A {\em subsequence}\/ of $\eta = (x_0,x_1,\ldots,x_k)$
is any marker of the form $\eta' = (x_{j_0},x_{j_1},\ldots,x_{j_l})$
for some choice of indices $0 \le j_0 < j_1 < \ldots < j_l \le k$.

A {\em cloud}\/ $\nu = (\nu_\eta)$ is a real-valued
function on the set $X^{[\infty]}$ of markers.
[We shall sometimes write $\nu(\eta)$ as a synonym for $\nu_\eta$.]
We say that $\nu$ has {\em finite support}\/ if $\nu_\eta = 0$
for all but finitely many markers $\eta$.
We say that $\nu$ is {\em supported on levels $\le N$}\/
if $\nu_\eta = 0$ whenever $\hbox{level}(\eta) > N$.
We define
\be
   \triplenorm \nu \triplenorm  \;\bydef\;
   \sup_{\eta} \sum_{\sigma \ancestor \eta} |\nu_\sigma|  \;,
 \label{def.cloudnorm}
\ee
and we say that $\nu$ has {\em finite norm}\/ if
$\triplenorm \nu \triplenorm < \infty$.
The set of clouds of finite norm forms a Banach space
with the norm $\triplenorm \,\cdot\, \triplenorm$.

Clouds allow us to give an abstract combinatorial representation
of the algebra of operators generated by $\alpha$ and the
multiplication operators, independently of any particular choice
of the matrix $\alpha$.
To see this, let us associate to the marker $\eta = (x_0,x_1,\ldots,x_k)$
the operator
\be
   T_\eta  \;\bydef\; I_{\{x_0\}} \alpha I_{\{x_1\}} \alpha \cdots
                         \alpha I_{\{x_{k-1}\}} \alpha I_{\{x_k\}}
   \;.
\ee
[The level of a marker thus corresponds to the number of
 factors $\alpha$ in the corresponding operator.
 In physical terms, a marker $\eta = (x_0,x_1,\ldots,x_k)$
 represents a piece of dirt that has traveled from $x_0$ to $x_k$
 via the path $x_0 \to x_1 \to \ldots \to x_k$.]
More generally, to a cloud $\nu$ we associate the operator $T_\nu$
defined by
\be
   T_\nu  \;\bydef\;  \sum_{\eta}  \nu_\eta \, T_\eta   \;.
\ee
(Initially this formula makes sense only for clouds of finite support,
but we will soon extend the definition to clouds of finite norm.)
Conversely, given any finite sum of operators
of the form $I_{f_0} \alpha I_{f_1} \cdots I_{f_{k-1}} \alpha I_{f_k}$,
we can expand it out as a (possibly infinite) sum of terms
$I_{\{x_0\}} \alpha I_{\{x_1\}} \cdots I_{\{x_{k-1}\}} \alpha I_{\{x_k\}}$
and thus represent it in the form $T_\nu$ for some cloud $\nu$
supported on a finite number of levels;
furthermore, $\nu$ has finite norm if the functions $f_i$ are bounded.
Of course, to make these considerations precise we will have
to deal with possibly infinite sums and specify the exact classes
of clouds and operators under consideration.

\begin{proposition}[Extension of \mbox{\protect\boldmath $T_\nu$}]
   \label{prop.Tnu.norm}
\FH
Let $\nu$ be a cloud of finite support.
Then
\be
   \| T_\nu \|_{w \to w}   \;\le\;   \triplenorm \nu \triplenorm
   \;.
 \label{eq.prop.Tnu.norm}
\ee
Therefore, the definition of the operator $T_\nu$
can be extended by linearity and continuity to all clouds
of finite norm,
and the map $\nu \mapsto T_\nu$ is a contraction.
\end{proposition}

\proof
Because $\alpha$ has nonnegative entries,
it suffices to consider the case in which $\nu_\eta\ge 0$ for all $\eta$. 
We assume that the cloud $\nu$ is supported on levels $\le N$,
and we shall prove the proposition by induction on $N$.
If $N=0$, \reff{eq.prop.Tnu.norm} is straightforward.
Suppose therefore that $N>0$, and consider a vector $c=(c_x)_{x\in X}$.
(Again, we can suppose $c_x\ge 0$ for all $x\in X$.) 
We have
\begin{equation}
  \label{eq:normTnu}
  \nw{cT_\nu}  \;=\;  \sum_{k=0}^N \sum_{\eta=(x_0,\ldots,x_k)}
      \nu_\eta \, c_{x_0}\alpha_{x_0 x_1}\alpha_{x_1 x_2}\cdots 
                  \alpha_{x_{k-1}x_k} w_{x_k}  \;.
\end{equation}
In this sum, the contribution of the markers $\eta$
such that $\hbox{level}(\eta)\in\{N-1,N\}$ is
\be
  \sum_{\eta=(x_0,\ldots,x_{N-1})}
  c_{x_0}\alpha_{x_0 x_1} \cdots \alpha_{x_{N-2}x_{N-1}} 
  \left( \nu_\eta w_{x_{N-1}} +
         \sum_{x_N\in X} \nu_{\eta \circ x_N} \alpha_{x_{N-1}x_N} w_{x_N}
  \right)
 \label{eq:normTnu.contribution}
\ee
where $\eta \circ x_N$ denotes the marker obtained by appending
the element $x_N$ to $\eta$.
Using the Fundamental Hypothesis in the form
\be
    \sum_{x_N\in X} \alpha_{x_{N-1}x_N} w_{x_N} \;\le\; w_{x_{N-1}}  \;,
\ee
the contribution \reff{eq:normTnu.contribution} can be bounded by
\be
  \sum_{\eta=(x_0,\ldots,x_{N-1})}
  c_{x_0}\alpha_{x_0 x_1} \cdots \alpha_{x_{N-2}x_{N-1}} w_{x_{N-1}} 
  \left(\nu_\eta+\sup_{x_N\in X}\nu_{\eta \circ x_N}\right)
  \;.
\ee
Thus, we obtain
\be
    \nw{cT_\nu}\ \le\ \nw{cT_{\nu'}}
\ee
where $\nu'$ is the cloud defined by
\be
    \nu'_\eta   \;=\;
    \cases{ 0  & if $\hbox{level}(\eta) = N$       \cr
            \nu_\eta + \sup_{x_N\in X}\nu_{\eta \circ x_N}
               & if $\hbox{level}(\eta)=N-1$       \cr
            \nu_\eta    & otherwise                \cr
          }
 \label{def_of_nuprime}
\ee            
The cloud $\nu'$ is supported on levels $\le N-1$
and satisfies $\triplenorm \nu' \triplenorm = \triplenorm \nu \triplenorm$;
this completes the proof.
\qed

Since our goal is to represent combinatorially
the algebra of operators generated by $\alpha$ and the $I_f$,
we need to introduce an operation on clouds that corresponds
to the multiplication of operators.
We do this as follows:
Given two markers $\eta$ and $\eta'$, we say that
$\eta$ {\em leads to}\/ $\eta'$ (and write $\eta \to \eta'$)
in case $\hbox{last}(\eta) = \hbox{first}(\eta')$.
If $\eta \to \eta'$, we define $\eta * \eta'$ to be the
concatentation of $\eta$ with $\eta'$ with the proviso that the
element $\hbox{last}(\eta) = \hbox{first}(\eta')$ is not repeated.
That is, if $\eta=(x_0,x_1,\ldots,x_k)$ and
$\eta'=(y_0,y_1,\ldots,y_{k'})$ with $x_k = y_0$,
then $\eta * \eta' \bydef (x_0,x_1,\ldots,x_k,y_1,\ldots,y_{k'})$.
Note that
$\hbox{level}(\eta * \eta') = \hbox{level}(\eta) + \hbox{level}(\eta')$.
Finally, given clouds $\nu_1$ and $\nu_2$, we define
the {\em convolution}\/ $\nu_1 * \nu_2$ by
\be
   (\nu_1 * \nu_2)(\eta)  \;=\;
   \!\!
   \sum_{\begin{scarray}
            \eta_1 \to \eta_2 \\
            \eta_1 * \eta_2 = \eta
         \end{scarray}}
    \nu_1(\eta_1) \, \nu_2(\eta_2)   \;.
 \label{def.convolution}
\ee
(Note that this convolution is associative but non-commutative.)
Since the sum \reff{def.convolution} is finite for each $\eta$,
convolution is well-defined for arbitrary pairs of clouds $\nu_1,\nu_2$.
Moreover, if $\nu_1$ and $\nu_2$ are of finite support
(resp.\ supported on finitely many levels),
then so is $\nu_1 * \nu_2$.
Furthermore, we have:

\begin{proposition}[Norm boundedness of convolution]
   \label{prop.norm.convolution}
Let $\nu_1$ and $\nu_2$ be clouds of finite norm.
Then
\be
   \triplenorm \nu_1 * \nu_2 \triplenorm
   \;\,\le\;\,
   \triplenorm \nu_1 \triplenorm  \;\, \triplenorm \nu_2 \triplenorm
   \;.
\ee
\end{proposition}

\proof
For any marker $\eta$ of level $k$, we have
\begin{eqnarray}
   \sum_{\sigma \ancestor \eta} |(\nu_1 * \nu_2)(\sigma)|
   & = &
   \sum_{j=0}^k |(\nu_1 * \nu_2)(\eta_0^j)|
            \nonumber \\[1mm]
   & = &
   \sum_{j=0}^k \;
       \left| \, \sum_{i=0}^j \nu_1(\eta_0^i) \, \nu_2(\eta_i^j) \, \right|
            \nonumber \\[1mm]
   & \le &
   \sum_{i=0}^k  \sum_{j=i}^k  |\nu_1(\eta_0^i)| \; |\nu_2(\eta_i^j)|
            \nonumber \\[1mm]
   & = &
   \sum_{i=0}^k  |\nu_1(\eta_0^i)|  \,
       \sum_{\sigma \ancestor \eta_i^k} |\nu_2(\sigma)|
            \nonumber \\[1mm]
   & \le &
   \sum_{i=0}^k  |\nu_1(\eta_0^i)| \;\, \triplenorm \nu_2 \triplenorm
            \nonumber \\[1mm]
   & \le &
   \triplenorm \nu_1 \triplenorm  \;\, \triplenorm \nu_2 \triplenorm
   \;.
\end{eqnarray}
\qed

It is now easily verified that convolution of clouds
corresponds to multiplication of operators, i.e.
\be
   T_{\nu_1 * \nu_2}  \;=\; T_{\nu_1} \, T_{\nu_2}
 \label{eq.nu1starnu2}
\ee
whenever $\nu_1,\nu_2$ are clouds of finite norm.
Indeed, this formula is easily seen to be true for clouds of finite support;
it then extends by continuity to clouds of finite norm
as a consequence of   
Propositions~\ref{prop.Tnu.norm} and \ref{prop.norm.convolution}.
Of course, the identity \reff{eq.nu1starnu2} is no accident:
we defined convolution so that \reff{eq.nu1starnu2} would hold!
Note in particular the importance of the ``leads to'' relation
and the non-repetition of $\hbox{last}(\eta) = \hbox{first}(\eta')$:
these implement the identity $I_{\{x\}} I_{\{y\}} = \delta_{xy} I_{\{x\}}$.

\subsection{Some special clouds}  \label{sec4.2}

Let us now define clouds to represent
each of the special operators introduced in Section~\ref{sec2.1}.

\bigskip

{\bf Level indicators.}  For each $k \ge 0$, we define a cloud $\rho^k$ by
\be
   (\rho^k)_\eta   \;=\;   \cases{ 1  & if $\hbox{level}(\eta) = k$  \cr
                                  0  & otherwise                    \cr
                                }
\ee
It is immediate that $\triplenorm \rho^k \triplenorm = 1$
and that $T_{\rho^k} = \alpha^k$;  in particular, $T_{\rho^0} = I$.
Note also that $\rho^0$ is the two-sided identity for convolution:
$\rho^0 * \nu = \nu * \rho^0 = \nu$ for every cloud $\nu$.
Furthermore, $\rho^k * \rho^\ell = \rho^{k+\ell}$.

For any cloud $\nu$, we shall denote its convolution powers
by $\nu^{*n}$ (or simply $\nu^n$),
with the convention that $\nu^{*0} = \rho^0$.
In particular, we have
$(\rho^k)^{*\ell} = \rho^{k\ell}$ for all $k,\ell \ge 0$.

\bigskip

{\bf Clouds associated to multiplication operators.}
To each function $f \colon\, X \to \R$,
we associate a cloud that we shall call
(by slight abuse of notation) $I_f$:
\be
    (I_f)_{(x_0,\ldots,x_k)}
    \;=\;  \cases{ f(x_0)     &  if $k=0$  \cr
                   0          &  if $k \ge 1$  \cr
                 }
\ee
We also write $I_\Lambda$ as a shorthand for $I_{\chi_\Lambda}$.
It is easy to verify that $\triplenorm I_f \triplenorm = \| f \|_\infty$
and $T_{I_f} = I_f$.

\bigskip

{\bf Clouds associated to cleaning operators.}
To each function $f \colon\, X \to [0,1]$,
we associate a cloud that we shall call
(by slight abuse of notation) $\beta_f$:
\be
    (\beta_f)_{(x_0,\ldots,x_k)}
    \;=\;  \cases{ 1-f(x_0)     &  if $k=0$  \cr
                   f(x_0)       &  if $k=1$  \cr
                   0            &  if $k \ge 2$  \cr
                 }
 \label{def.cloud.betaf}
\ee
We also write $\beta_\Lambda$ as a shorthand for $\beta_{\chi_\Lambda}$.
It is easy to verify that $\triplenorm \beta_f \triplenorm = 1$,
$\beta_f = I_{1-f} + I_f * \rho^1$
and $T_{\beta_f} = \beta_f$.

\bigskip

{\bf Clouds associated to balayage operators.}
For each $\Lambda \subseteq X$, we denote by $\partial\Lambda$
the set of markers of the form $\eta = (x_0,\ldots,x_k)$ [$k \ge 0$]
for which $x_0,\ldots,x_{k-1} \in \Lambda$ and $x_k \in \Lambda^c$.
We then define a cloud $\pi_\Lambda$ by
\be
   (\pi_\Lambda)_\eta
   \;=\;
   \cases{ 1   & if $\eta \in \partial\Lambda$  \cr
           0   & otherwise                      \cr
         }
\ee
It is easy to verify that
$\triplenorm \pi_\Lambda \triplenorm = 1$ whenever $\Lambda \neq X$,
$\pi_X = 0$,
and $T_{\pi_\Lambda} = \Pi_\Lambda$.
Note also that $\pi_\Lambda * \pi_\Lambda = \pi_\Lambda$.

\bigskip

{\bf Indicator clouds.}
We denote by ${\bf 1}$ the cloud that takes the value 1 on all markers.
More generally, for $\Lambda \subseteq X$, we denote by ${\bf 1}_\Lambda$
the cloud that takes the values
\be
   ({\bf 1}_\Lambda)_\eta
   \;=\;
   \cases{ 1   & if $\eta$ has all its entries in $\Lambda$ \cr
           0   & otherwise                                  \cr
         }
\ee
{\em Formally}\/ we have
$T_{\bf 1} = \sum_{k=0}^\infty \alpha^k = (I-\alpha)^{-1}$ and
$T_{{\bf 1}_\Lambda} =
 I_\Lambda \sum_{k=0}^\infty (I_\Lambda \alpha I_\Lambda)^k I_\Lambda =
 I_\Lambda (I - I_\Lambda \alpha I_\Lambda)^{-1} I_\Lambda$;
but since $\triplenorm {\bf 1} \triplenorm =
           \triplenorm {\bf 1}_\Lambda \triplenorm = \infty$
(for $\Lambda \neq \emptyset$),
there is no guarantee that these series converge.
The advantage of the tree formalism is that it makes sense
to speak of the clouds ${\bf 1}$ and ${\bf 1}_\Lambda$
without worrying about convergence questions.

\subsection{A partial ordering on clouds}   \label{sec4.3}

The space of clouds is obviously endowed with the ``pointwise''
partial ordering defined by
\be
   \mu \,\le\, \nu  \qquad\hbox{if and only if}\qquad
   \mu_\eta \,\le\, \nu_\eta \;\hbox{for all}\; \eta  \;.
\ee
(In particular, $\nu \ge 0$ means that a cloud is nonnegative.)
We would now like to introduce a weaker partial ordering
that will play an important role in the sequel;
roughly speaking, it formalizes the comparison of operators
by their ``efficiency of cleaning''.

To each cloud $\nu$, we associate another cloud $\widetilde{\nu}$,
called its {\em cumulative distribution}\/, defined by
\be
   \widetilde{\nu}_\eta  \;=\;
   \sum_{\sigma \ancestor \eta}  \nu_\sigma
\ee
or in other words by
\be
   \widetilde{\nu}  \;=\;  \nu * {\bf 1}   \;.
\ee
(Note also that
 $\triplenorm \nu \triplenorm = \| \widetilde{\;|\nu|\;} \|_\infty$.)
We then introduce a partial ordering $\audessus$ on the space of clouds,
defined by:
\begin{definition}[Definition of \mbox{\protect\boldmath $\audessus$}]
   \label{def.audessus}
$\mu \audessus \nu$ if and only if $\widetilde{\mu} \le \widetilde{\nu}$,
i.e.
\be
   \mu \,\audessus\, \nu  \qquad\hbox{if and only if}\qquad
   \widetilde{\mu}_\eta \,\le\, \widetilde{\nu}_\eta \;\hbox{for all}\; \eta
   \;.
\ee
\end{definition}
Obviously $\mu \le \nu$ implies $\mu \audessus \nu$ but not conversely.

Let us begin by proving two useful formulae
for the cumulative distribution of a convolution.
If $\eta=(x_0,\ldots,x_k)$ is any marker
and $\sigma=(x_0,\ldots,x_j)$ is any ancestor of $\eta$,
we denote by $\eta \doubleminus \sigma$ the marker $(x_j,\ldots,x_k)$;
it is the unique marker $\xi$ satisfying $\sigma * \xi = \eta$.
Also, if $\eta$ is any marker of level $k \ge 1$,
we define $\eta^- = \eta_0^{k-1}$ and $\eta_- = \eta_1^k$.
If $\eta$ is a marker of level 0,
we define $\eta^- = \eta_- = \emptyset$ (the empty sequence).

\begin{lemma}[Cumulative distribution of a convolution]
  \label{lemma.cumul.convolution}
Let $\mu,\nu$ be clouds.  Then
\be
   \widetilde{(\mu * \nu)}  \;=\;  \mu * \widetilde{\nu}
 \label{eq.cumul.convolution.0}
\ee
or equivalently
\begin{eqnarray}
   \widetilde{(\mu * \nu)}(\eta)
   & = &
   \sum_{\sigma \anc \eta}  \mu(\sigma) \;
      \widetilde{\nu}(\eta \doubleminus \sigma)
          \label{eq.cumul.convolution}  \\[3mm]
   & = &
   \sum_{\sigma \anc \eta}  \widetilde{\mu}(\sigma) \,
      \Bigl[ \widetilde{\nu}(\eta \doubleminus \sigma) \,-\,
             \widetilde{\nu}((\eta \doubleminus \sigma)_-)
      \Bigr]
          \label{eq.cumul.convolution.2}
\end{eqnarray}
where in \reff{eq.cumul.convolution.2} we make the convention
that $\widetilde{\nu}(\emptyset) = 0$.
\end{lemma}

\proof
\reff{eq.cumul.convolution.0} is an immediate consequence of
the associativity of convolution:
\be
   \widetilde{(\mu * \nu)}  \;=\;
   (\mu * \nu) * {\bf 1}   \;=\;
   \mu * (\nu * {\bf 1})   \;=\;
   \mu * \widetilde{\nu}   \;.
\ee
\reff{eq.cumul.convolution} is just a rewriting of this.

We now insert in \reff{eq.cumul.convolution} the identity
$\mu(\sigma) = \widetilde{\mu}(\sigma) - \widetilde{\mu}(\sigma^-)$
with the convention $\widetilde{\mu}(\emptyset) = 0$.
This yields
\be
   \widetilde{(\mu * \nu)}(\eta)
   \;=\;
   \sum_{\sigma \anc \eta}  \widetilde{\mu}(\sigma) \;
      \widetilde{\nu}(\eta \doubleminus \sigma)
   \;-\;
   \sum_{\begin{scarray}
            \sigma \anc \eta \\
            \hbox{\scriptsize level}(\sigma) \ge 1
         \end{scarray}}
      \widetilde{\mu}(\sigma^-) \;
      \widetilde{\nu}(\eta \doubleminus \sigma)
   \;.
\ee
In the second sum on the right, we make the change of variables
$\xi = \sigma^-$
and observe that $\eta \doubleminus \sigma = (\eta \doubleminus \xi)_-$;
this sum therefore becomes
\be
   \sum_{\xi \anc \eta^-} \widetilde{\mu}(\xi) \;
      \widetilde{\nu}((\eta \doubleminus \xi)_-)
   \;.
\ee
The term $\xi = \eta$ can now be adjoined to this sum,
thanks to the convention $\widetilde{\nu}(\emptyset) = 0$.
This proves \reff{eq.cumul.convolution.2}.
\qed

In the special case $\nu = \beta_h$, we have the following formulae:

\begin{lemma}[Cumulative distribution of \mbox{\protect\boldmath $\beta_h$}]
   \label{lemma.cumul.betah}
Let $h \colon\, X \to [0,1]$.  Then
\be
   \widetilde{\beta_h}  \;=\;  {\bf 1} \,-\, I_h   \;.
 \label{eq.lemma.cumul.betah}
\ee
\end{lemma}

\proof
{}From the definition \reff{def.cloud.betaf} we easily obtain
\be
    (\widetilde{\beta_h})_{(x_0,\ldots,x_k)}
    \;=\;  \cases{ 1-h(x_0)     &  if $k=0$  \cr
                   1            &  if $k \ge 1$  \cr
                 }
 \label{eq.cumul.betah}
\ee
which is equivalent to \reff{eq.lemma.cumul.betah}.
\qed

\begin{corollary}[Cumulative distribution of a convolution with
                          \mbox{\protect\boldmath $\beta_h$}]
   \label{cor.somme_ancetres}
Let $\mu$ be a cloud, and let $h \colon\, X \to [0,1]$.
Then
\be
   \widetilde{\mu * \beta_h}   \;=\;  \widetilde{\mu} \,-\, \mu * I_h
 \label{somme_ancetres3a}
\ee
or equivalently
\begin{eqnarray}
   (\widetilde{\mu * \beta_h})(\eta)
   & = &
      \widetilde{\mu}(\eta) \;-\;
      h(\hbox{\rm last}(\eta)) \: \mu(\eta)
               \label{somme_ancetres0}  \\[3mm]
   & = &
   \widetilde{\mu}(\eta^-) 
       \;+\; [1-h(\hbox{\rm last}(\eta))] \; \mu(\eta)
               \label{somme_ancetres}  \\[3mm]
   & = &
      h(\hbox{\rm last}(\eta)) \; \widetilde{\mu}(\eta^-) \;+\;
      [1-h(\hbox{\rm last}(\eta))] \; \widetilde{\mu}(\eta)
      \;.
               \label{somme_ancetres2}
\end{eqnarray}
\end{corollary}

\proof
The formula \reff{somme_ancetres3a} is an immediate consequence of
\reff{eq.cumul.convolution.0} and \reff{eq.lemma.cumul.betah}.
The alternate forms \reff{somme_ancetres0}--\reff{somme_ancetres2}
are trivial rewritings.
\qed

\begin{corollary}
   \label{cor_to_lemma.somme_ancetres}
If $\mu \ge 0$ and $0 \le g \le h \le 1$,
then $\mu * \beta_h \audessus \mu * \beta_g$.
In particular, $\beta_h \audessus \beta_g$.
\end{corollary}

\proof
This is an immediate consequence of \reff{somme_ancetres}.
\qed

\begin{corollary}[Convolution from the right by
                          \mbox{\protect\boldmath $\beta_h$}]
   \label{cor2_to_lemma.somme_ancetres}
Let $\mu_1,\mu_2$ be clouds and let $0 \le h \le 1$.
If $\mu_1 \audessus \mu_2$, then $\mu_1 * \beta_h \audessus \mu_2 * \beta_h$.
\end{corollary}

\proof
This is an immediate consequence of \reff{somme_ancetres2}.
\qed

For general clouds $\nu$, one can obtain an
{\em inequality}\/ analogous to \reff{somme_ancetres}:

\begin{lemma}[Inequality for cumulative distribution of a convolution]
   \label{lemma.10Apr07.a}
\quad\par\noindent
Let $\mu,\nu \ge 0$ be clouds with $\triplenorm \nu \triplenorm < \infty$.
Then, for every marker $\eta$,
\be
   \widetilde{(\mu * \nu)}(\eta)
    \;\,\le\;\,
   \widetilde{\mu}(\eta^-) \: \triplenorm \nu \triplenorm 
    \;+\;
   \mu(\eta) \: \nu\Bigl((\hbox{\rm last}(\eta) )\Bigr)
   \;.
 \label{eq.lemma.10Apr07.a}
\ee
\end{lemma}

\proof
Rewrite \reff{eq.cumul.convolution} as
\be
  \widetilde{(\mu*\nu)}(\eta)
  \;=\;
  \sum_{\sigma \ancneq \eta}  \mu(\sigma) \;
      \widetilde{\nu}(\eta \doubleminus \sigma)
  \;+\;
  \mu(\eta) \, \nu\Bigl((\hbox{last}(\eta) )\Bigr)
  \;.
\ee
Using
$\widetilde{\nu}(\eta \doubleminus \sigma) \le \triplenorm \nu \triplenorm$
in the first term, we obtain \reff{eq.lemma.10Apr07.a}.
\qed

See Lemma~\ref{lemma.12Apr07.cloud} below
for an interesting application of Lemma~\ref{lemma.10Apr07.a}.

\bigskip

It is also useful to know under what conditions a cloud $\mu$
is majorized by a cloud $\beta_h$.  The easy answer is the following:

\begin{lemma}[Majorization by \mbox{\protect\boldmath $\beta_h$}]
   \label{lemma.mubetah}
Let $\mu \ge 0$ be a cloud and let $h \colon\, X \to [0,1]$.
Then the following are equivalent:
\begin{itemize}
   \item[(a)] $\mu \audessus \beta_h$.
   \item[(b)] $\triplenorm \mu \triplenorm \le 1$
      and $\mu\bigl((x)\bigr) \le 1-h(x)$ for every marker $(x)$ of level 0.
\end{itemize}
\end{lemma}

\proof
$\mu \audessus \beta_h$ means that
$\widetilde{\mu}(\eta) \le \widetilde{\beta}_h(\eta)$
for every marker $\eta$.
By \reff{eq.cumul.betah}, this means
precisely that $\mu\bigl((x)\bigr) \le 1-h(x)$ for every $x \in X$
and that $\widetilde{\mu}(\eta) \le 1$ for all $\eta$.
\qed

\begin{corollary}[Comparison of \mbox{\protect\boldmath $\beta_g$}
                             with \mbox{\protect\boldmath $\beta_h$}]
   \label{cor.mubetah}
Let $g,h \colon\, X \to [0,1]$.
Then $\beta_g \audessus \beta_h$ if and only if $g \ge h$.
\end{corollary}


For the cumulative distribution of $\beta_g * \beta_h$
we obtain the following important identity:

\begin{lemma}[Collapse identity]
  \label{lemma.collapse.clouds}
Let $g,h \colon\, X \to [0,1]$.  Then
\be
   \widetilde{\beta_g * \beta_h}  \;=\;
   \widetilde{\beta}_{1-(1-g)(1-h)}  \,-\, I_g * \rho^1 * I_h   \;.
 \label{identity.collapse.clouds}
\ee
\end{lemma}

\proof
We use the Kronecker delta notation $\delta_{ij} = 1$ if $i = j$
and 0 otherwise.
Let $\eta = (x_0,\ldots,x_k)$.
Then
\begin{eqnarray}
   (\widetilde{\beta_g * \beta_h})(\eta)
   & = &
   \widetilde{\beta_g}(\eta^-) \;+\; [1-h(x_k)] \, \beta_g(\eta)
         \nonumber \\[2mm]
   & = &
   1 \,-\, g(x_0) \delta_{k1} \;+\;
       [1-h(x_k)] \Bigl[ [1-g(x_0)] \delta_{k0} + g(x_0) \delta_{k1} \Bigr]
         \nonumber \\[2mm]
   & = &
   1 \,-\, [1-g(x_0)] [1-h(x_0)] \delta_{k0}
     \,-\, g(x_0) h(x_1) \delta_{k1}
         \nonumber \\[2mm]
   & = &
   \widetilde{\beta}_{1-(1-g)(1-h)}(\eta)  \;-\; (I_g * \rho^1 * I_h)(\eta)
\end{eqnarray}
where the first equality uses \reff{somme_ancetres},
the second uses \reff{eq.cumul.betah} and the definition of $\beta_g$,
the third is an easy rewriting,
and the fourth again uses \reff{eq.cumul.betah}.
\qed

\begin{corollary}[Collapse inequality]
\label{cor.collapse.clouds}
Let $g_i \colon\, X \to [0,1]$ for $i=1,\ldots,n$.
Then
\be
   \beta_{g_1} * \cdots * \beta_{g_n}
   \;\,\audessus\;\,
   \beta_{1-\prod_{i=1}^n(1-g_i)} \;.
\ee
\end{corollary}

\proof
The claim is trivial for $n=1$.
For $n=2$ it follows immediately from the identity
\reff{identity.collapse.clouds}
and the fact that $I_{g_1} * \rho^1 * I_{g_2} \ge 0$.
The cases $n \ge 3$ can be proven by an elementary induction
using the case $n=2$ along with Corollary~\ref{cor2_to_lemma.somme_ancetres}.
\qed

\bigskip

We next need to know under what circumstances the
partial ordering $\audessus$ is preserved by convolution.
Convolution from the left is easy:

\begin{proposition}[Convolution from the left]
  \label{prop.convolution.left}
Let $\mu,\nu_1,\nu_2$ be clouds satisfying
$\mu \ge 0$ and $\nu_1 \audessus \nu_2$.
Then $\mu * \nu_1 \audessus \mu * \nu_2$.
\end{proposition}

\proof
This is an immediate consequence of the formula \reff{eq.cumul.convolution}.
%
\qed

\noindent
Physically, Proposition~\ref{prop.convolution.left} says that
if $\nu_1$ is a better cleaner than $\nu_2$,
then this same relation holds if both cleaners are preceded
by an arbitrary nonnegative operation $\mu$.

The behavior of the partial ordering $\audessus$
under convolution from the right is considerably more subtle;
this reflects the fact that the partial ordering between cleaners
is not preserved by arbitrary {\em subsequent}\/ nonnegative operations.
Rather, we need to limit the class of subsequent operations
that are allowed:

\begin{proposition}[Convolution from the right]
  \label{prop.convolution.right}
Let $\nu \ge 0$ be a cloud.  Then the following are equivalent:
\begin{itemize}
   \item[(a)]  For all pairs of clouds $\mu_1,\mu_2$ satisfying
       $\mu_1 \audessus \mu_2$, we have
       $\mu_1 * \nu \audessus \mu_2 * \nu$.
   \item[(b)]  $\rho^k * \nu \audessus \nu$ for all $k \ge 0$.
   \item[(b$\,{}'$)]  $\rho^1 * \nu \audessus \nu$.
   \item[(c)]  $\widetilde{\nu}_{\eta'} \le \widetilde{\nu}_{\eta}$
       whenever $\eta'$ is a suffix of $\eta$.
   \item[(c$\,{}'$)]  $\widetilde{\nu}_{\eta'} \le \widetilde{\nu}_{\eta}$
       whenever $\eta$ is a marker of level $k \ge 1$ and $\eta' = \eta_1^k$.
\end{itemize}
\end{proposition}

\proof
The implications  (b) $\implies$ (b${}'$) and  (c) $\implies$ (c${}'$)
are trivial, and (c${}'$) $\implies$ (c) is easy.

{}From $\rho^1*\nu \audessus\nu$ it follows,
using Proposition~\ref{prop.convolution.left}
and the associativity of convolution, that
\be
    \rho^{k+1}*\nu \;=\; \rho^k*(\rho^1*\nu) \; \audessus\; \rho^k*\nu  \;.
\ee
Induction on $k$ then gives (b${}'$) $\implies$ (b).



Now let $\eta$ be a marker of level $\ell$.
By \reff{eq.cumul.convolution} we have
\be
    \widetilde{(\rho^k*\nu)}(\eta)  \;=\;
 \cases{ \widetilde{\nu}(\eta_{k}^\ell)  & if  $\ell\ge k$  \cr
         0                               & if  $\ell<k$     \cr
       }
\ee
{}From this, we easily get  (b) $\Longleftrightarrow$ (c)
and  (b${}'$) $\Longleftrightarrow$ (c${}'$).
Hence (b), (b${}'$), (c) and (c${}'$) are all equivalent.

The implication (a) $\implies$ (b) is also trivial,
because $\rho^0*\nu = \nu$ and
$\rho^k \audessus \rho^0$ for all $k\ge 0$. 

Finally, (c${}'$) $\implies$ (a) is an immediate consequence of
\reff{eq.cumul.convolution.2}, since (c${}'$) ensures that
the square brackets in \reff{eq.cumul.convolution.2}
are always nonnegative.
%
%
\qed

\begin{definition}
   \label{def_after_prop.convolution.right}
We denote by $\scrr$ the class of clouds $\nu \ge 0$
satisfying any one (hence all) of the equivalent conditions
of Proposition~\ref{prop.convolution.right}.
\end{definition}

\medskip

{\bf Remark.}
Given a cloud $\nu$, one can define a ``dual'' cloud $\nu^*$
by the intertwining relation
$(\rho^0 - \rho^1) * \nu = \nu^* * (\rho^0 - \rho^1)$
[compare \reff{eq:9}].
Convoluting on the right with ${\bf 1}$, we obtain
$\nu^*(\eta) = \widetilde{\nu}(\eta) - \widetilde{\nu}(\eta_-)$.
Thus,
the ``dual'' of a nonnegative cloud $\nu$
is nonnegative if and only if $\nu \in \scrr$.
(Note, however, that $(\nu^*)^* \neq \nu$,
 so this is not a true duality.)
\qed

\begin{corollary}[Multi-monotonicity]
   \label{cor.powers}
For $1 \le i \le n$, let $\mu_i,\nu_i \ge 0$ be clouds satisfying
$\mu_i \audessus \nu_i$ and $\nu_i \in \scrr$.
Then $\mu_1 * \ldots * \mu_n \audessus \nu_1 * \ldots * \nu_n$.
\end{corollary}

\proof
We use the telescoping identity
\be
   \mu_1 * \ldots * \mu_n \,-\, \nu_1 * \ldots * \nu_n
   \;=\;
   \sum_{i=1}^n \mu_1 * \ldots * \mu_{i-1} * (\mu_i - \nu_i) *
                \nu_{i+1} * \ldots * \nu_n
   \;.
\ee
By hypothesis, we have $\mu_i - \nu_i \audessus 0$ for all $i$.
By Proposition~\ref{prop.convolution.left},
we can convolve on the left with $\mu_1 * \ldots * \mu_{i-1}$;
and by Proposition~\ref{prop.convolution.right}(a)
and Definition~\ref{def_after_prop.convolution.right},
we can convolve on the right successively by $\nu_{i+1}, \ldots, \nu_n$.
\qed

\smallskip

{\bf Remark.}
In Corollary~\ref{cor.powers}, only $\nu_2,\ldots,\nu_n$ really need
to belong to $\scrr$, as we never convolve on the right with $\nu_1$.

\begin{corollary}[Multi-monotonicity for cleaners]
\label{cor.multi-monotonicity.clouds}
Let $0 \le g_i \le h_i \le 1$ for $i=1,\ldots,n$.
Then
\be
  \beta_{h_1} * \cdots * \beta_{h_n} 
  \;\,\audessus\;\,
  \beta_{g_1} * \cdots * \beta_{g_n} \;.
\ee
\end{corollary}

\proof
An immediate consequence of Corollaries~\ref{cor.mubetah}
and \ref{cor.powers}.
\qed

%

We now resume our study of the class $\scrr$:

\begin{proposition}
The class $\scrr$ forms a multiplicative convex cone.  That is,
\begin{itemize}
   \item[(a)]  If $\mu, \nu \in \scrr$ and $a,b \ge 0$,
       then $a\mu + b\nu \in \scrr$.
   \item[(b)]  If $\mu, \nu \in \scrr$, then $\mu * \nu \in \scrr$.
\end{itemize}
\end{proposition}

\proof
This is immediate from property (a) of
Proposition~\ref{prop.convolution.right}.
\qed

We know (Corollary~\ref{cor2_to_lemma.somme_ancetres})
that all the clouds $\beta_f$ ($0 \le f \le 1$) belong to $\scrr$,
as do all sums of convolutions thereof.
But the latter turn out to constitute a strictly smaller class,
as they satisfy a condition like that of
Proposition~\ref{prop.convolution.right}(c)
not only for suffixes but also for arbitrary subsequences:

\begin{proposition}
 \label{prop.clouds.Rbeta}
Let $\mu$ be a finite sum of clouds
of the form $a \beta_{f_1}*\ldots*\beta_{f_n}$
with $n \ge 0$, $0 \le f_1,\ldots,f_n \le 1$ and $a \ge 0$.
Then $\widetilde{\mu}_{\eta'} \le \widetilde{\mu}_\eta$
whenever $\eta'$ is a subsequence of $\eta$.
\end{proposition}

\proof
It suffices to prove the result for $\mu = \beta_{f_1}*\ldots*\beta_{f_n}$.
The case $n=0$ (i.e.\ $\mu = \rho^0$) is trivial, so assume $n \ge 1$.
First we consider the special case in which each $f_i=\ind{\Lambda_i}$.
In this case, we have $\widetilde{\mu}_\eta=0$ or 1
according to the following rule: 
\begin{eqnarray}
  \widetilde{\mu}_{(x_0,\ldots,x_k)} \,=\, 0
  & \Longleftrightarrow &
     \mbox{there exist indices } 1\le i_0< \ldots <i_k\le n
       \nonumber \\
  & & \quad
     \mbox{ such that } x_j \in \Lambda_{i_j}\mbox{ for }j=0,\ldots,k \,.
  \label{eq:convolution_of_beta_Lambda}
\end{eqnarray}
This is easily seen by induction on $n$:
recall that $\beta_\Lambda$ is supported on levels 0 and 1,
with $\beta_\Lambda(x_0)=\ind{\Lambda^c}(x_0)$
and $\beta_\Lambda(x_0,x_1)=\ind{\Lambda}(x_0)$,
so that \reff{eq:convolution_of_beta_Lambda} is clear for $n=1$. 
For $n\ge2$, let $\mu=\beta_{\Lambda_1}*\ldots*\beta_{\Lambda_n}$
and $\nu=\beta_{\Lambda_1}*\ldots*\beta_{\Lambda_{n-1}}$.
By \reff{somme_ancetres2}, we have for $\eta=(x_1,\ldots,x_k)$
\be
    \widetilde{\mu}_\eta   \;=\;
      \ind{\Lambda_n}(x_k) \, \widetilde{\nu}_{\eta_0^{k-1}} \,+\,
      \ind{\Lambda_n^c}(x_k) \, \widetilde{\nu}_{\eta}
    \;,
\ee
from which we can see that if $\nu$ satisfies
\reff{eq:convolution_of_beta_Lambda}, then so does $\mu$. 

Now, \reff{eq:convolution_of_beta_Lambda} obviously implies that if
$\widetilde{\mu}_\eta=0$ and $\eta'$ is a subsequence of $\eta$,
then $\widetilde{\mu}_{\eta'}=0$. 
This proves the Proposition in the special case $f_i=\ind{\Lambda_i}$.

To handle the general case,
we note that each function $f\colon\, X \to [0,1]$ can be written
in the form
\be
    \label{convex_decomposition_of_f}
    f  \;=\; \sum_{k\ge 0} a_k \ind{\Lambda_k}
\ee
where $a_k\ge 0$, $\sum_k a_k=1$,
and $(\Lambda_k)$ is a sequence of (possibly empty) subsets of $X$. 
[One way to get such a decomposition is to use the binary expansion of $f(x)$,
\begin{equation}
  f(x)  \;=\;  \sum_{k=1}^{\infty} d_k(x)\, 2^{-k}
\end{equation}
with $d_k(x)\in\{0,1\}$,
and then to set $a_k\bydef 2^{-k}$
and $\Lambda_k\bydef\{x\in X \colon\, d_k(x)=1\}$.]
{}From~\reff{convex_decomposition_of_f} together with the fact that
$\beta_f$ is affine in $f$,
it follows that each cloud of the form
$\beta_{f_1}*\ldots*\beta_{f_n}$ is a convex combination
(with a finite or countably infinite number of terms)
of clouds of the form
$\beta_{\ind{\Lambda_1}}*\ldots*\beta_{\ind{\Lambda_n}}$. 
Since the set of clouds satisfying the conclusion of the Proposition is
obviously stable under convex combination, the proof is complete. 
\qed

\begin{definition}
   \label{def_after_prop.clouds.Rbeta}
We denote by $\scrb$ the class consisting of finite sums of clouds
of the form $a \beta_{f_1}*\ldots*\beta_{f_n}$
with $n \ge 0$, $0 \le f_1,\ldots,f_n \le 1$ and $a \ge 0$.

We denote by $\scrp$ the class of clouds $\mu \ge 0$
for which $\widetilde{\mu}_{\eta'} \le \widetilde{\mu}_\eta$
whenever $\eta'$ is a subsequence of $\eta$.
\end{definition}

We have just shown that $\scrb \subseteq \scrp \subseteq \scrr$.
The following examples show that both these inclusions are strict:

\begin{example}
\rm
Let $|X| \ge 2$ and $\Lambda \subseteq X$ with $\Lambda \neq \emptyset, X$.
Then $\pi_\Lambda \in \scrp$ but $\pi_\Lambda \notin \scrb$.
One might worry that this example is somehow ``pathological''
because $\pi_\Lambda$ is not supported on finitely many levels.
But the next example avoids this objection \ldots
\end{example}

\begin{figure}
\begin{center}
\begin{picture}(0,0)%
\includegraphics{P_not_B.pstex}%
\end{picture}%
\setlength{\unitlength}{4144sp}%
\begingroup\makeatletter\ifx\SetFigFont\undefined%
\gdef\SetFigFont#1#2#3#4#5{%
  \reset@font\fontsize{#1}{#2pt}%
  \fontfamily{#3}\fontseries{#4}\fontshape{#5}%
  \selectfont}%
\fi\endgroup%
\begin{picture}(4621,1818)(586,-2531)
\put(2026,-1771){\makebox(0,0)[lb]{\smash{{\SetFigFont{10}{12.0}{\rmdefault}{\mddefault}{\updefault}$xy$}}}}
\put(1351,-1096){\makebox(0,0)[lb]{\smash{{\SetFigFont{10}{12.0}{\rmdefault}{\mddefault}{\updefault}$xxy$}}}}
\put(4553,-1771){\makebox(0,0)[lb]{\smash{{\SetFigFont{10}{12.0}{\rmdefault}{\mddefault}{\updefault}$yy$}}}}
\put(4823,-1096){\makebox(0,0)[lb]{\smash{{\SetFigFont{10}{12.0}{\rmdefault}{\mddefault}{\updefault}$yyy$}}}}
\put(1441,-2491){\makebox(0,0)[lb]{\smash{{\SetFigFont{10}{12.0}{\rmdefault}{\mddefault}{\updefault}$x$}}}}
\put(3968,-2476){\makebox(0,0)[lb]{\smash{{\SetFigFont{10}{12.0}{\rmdefault}{\mddefault}{\updefault}$y$}}}}
\put(871,-1786){\makebox(0,0)[lb]{\smash{{\SetFigFont{10}{12.0}{\rmdefault}{\mddefault}{\updefault}$xx$}}}}
\put(3398,-1786){\makebox(0,0)[lb]{\smash{{\SetFigFont{10}{12.0}{\rmdefault}{\mddefault}{\updefault}$yx$}}}}
\put(586,-1096){\makebox(0,0)[lb]{\smash{{\SetFigFont{10}{12.0}{\rmdefault}{\mddefault}{\updefault}$xxx$}}}}
\put(3985,-1132){\makebox(0,0)[lb]{\smash{{\SetFigFont{10}{12.0}{\rmdefault}{\mddefault}{\updefault}$yyx$}}}}
\end{picture}%
\end{center}
\caption{
   The cloud $\mu$ in Example~\ref{example_cloud_P_not_B}
   takes the value 1 (resp.\ 0) on the markers
   indicated by a full (resp.\ empty) circle.
   This cloud belongs to $\scrp$ but not to $\scrb$.
}
  \label{fig_cloud_P_not_B}
\end{figure}

\begin{figure}
\begin{center}
\begin{picture}(0,0)%
\includegraphics{contre_exemple.pstex}%
\end{picture}%
\setlength{\unitlength}{4144sp}%
\begingroup\makeatletter\ifx\SetFigFont\undefined%
\gdef\SetFigFont#1#2#3#4#5{%
  \reset@font\fontsize{#1}{#2pt}%
  \fontfamily{#3}\fontseries{#4}\fontshape{#5}%
  \selectfont}%
\fi\endgroup%
\begin{picture}(4407,2490)(586,-2531)
\put(2026,-1771){\makebox(0,0)[lb]{\smash{\SetFigFont{10}{12.0}{\rmdefault}{\mddefault}{\updefault}$xy$}}}
\put(1351,-1096){\makebox(0,0)[lb]{\smash{\SetFigFont{10}{12.0}{\rmdefault}{\mddefault}{\updefault}$xxy$}}}
\put(586,-1096){\makebox(0,0)[lb]{\smash{\SetFigFont{10}{12.0}{\rmdefault}{\mddefault}{\updefault}$xxx$}}}
\put(2296,-1096){\makebox(0,0)[lb]{\smash{\SetFigFont{10}{12.0}{\rmdefault}{\mddefault}{\updefault}$xyy$}}}
\put(1853,-1096){\makebox(0,0)[lb]{\smash{\SetFigFont{10}{12.0}{\rmdefault}{\mddefault}{\updefault}$xyx$}}}
\put(4553,-1771){\makebox(0,0)[lb]{\smash{\SetFigFont{10}{12.0}{\rmdefault}{\mddefault}{\updefault}$yy$}}}
\put(3878,-1096){\makebox(0,0)[lb]{\smash{\SetFigFont{10}{12.0}{\rmdefault}{\mddefault}{\updefault}$yxy$}}}
\put(3113,-1096){\makebox(0,0)[lb]{\smash{\SetFigFont{10}{12.0}{\rmdefault}{\mddefault}{\updefault}$yxx$}}}
\put(4823,-1096){\makebox(0,0)[lb]{\smash{\SetFigFont{10}{12.0}{\rmdefault}{\mddefault}{\updefault}$yyy$}}}
\put(4380,-1096){\makebox(0,0)[lb]{\smash{\SetFigFont{10}{12.0}{\rmdefault}{\mddefault}{\updefault}$yyx$}}}
\put(1441,-2491){\makebox(0,0)[lb]{\smash{\SetFigFont{10}{12.0}{\rmdefault}{\mddefault}{\updefault}$x$}}}
\put(3968,-2476){\makebox(0,0)[lb]{\smash{\SetFigFont{10}{12.0}{\rmdefault}{\mddefault}{\updefault}$y$}}}
\put(871,-1786){\makebox(0,0)[lb]{\smash{\SetFigFont{10}{12.0}{\rmdefault}{\mddefault}{\updefault}$xx$}}}
\put(3398,-1786){\makebox(0,0)[lb]{\smash{\SetFigFont{10}{12.0}{\rmdefault}{\mddefault}{\updefault}$yx$}}}
\end{picture}
\end{center}
\caption{
   The cloud $\mu$ in Example~\ref{example_cloud1}
   takes the value 1 (resp.\ 0) on the markers
   indicated by a full (resp.\ empty) circle.
   This cloud belongs to $\scrr$ but not to $\scrp$.
}
  \label{fig_cloud1}
\end{figure}

\begin{example}
 \label{example_cloud_P_not_B}
\rm
Let $X=\{x,y\}$, and let $\mu$ be the cloud that takes the value 1
on the markers $xy$, $xxx$, $xxy$, $yx$, $yyx$ and $yyy$
and takes the value 0 elsewhere (Figure~\ref{fig_cloud_P_not_B}).
This is a cloud of norm 1, supported on levels $\le 2$.
To verify that $\mu \in \scrp$,
it suffices to check that for each $\eta$
satisfying $\widetilde{\mu}_{\eta} = 0$,
one has $\widetilde{\mu}_{\eta'} = 0$ for each subsequence $\eta'$ of $\eta$.
We leave this verification to the reader.

On the other hand, we claim that $\mu \notin \scrb$.
Indeed, suppose that we could write $\mu$ in the form
\be
\label{sum_of_nu_i}
\mu \;=\;  \sum_{i=1}^n  a_i \, \beta_{f_1^{(i)}}*\ldots*\beta_{f_{n_i}^{(i)}}
\ee
with all $a_i > 0$.
Then all markers on which $\mu$ takes the value 0
must also be given mass 0 by the cloud
$\nu_i \bydef \beta_{f_1^{(i)}}*\ldots*\beta_{f_{n_i}^{(i)}}$.
Furthermore, it is easily seen that $\nu_i$ satsifies
$\widetilde{\nu_i}(\eta) = 1$ for each marker $\eta$ of level $\ge n_i$
(see also Corollary~\ref{cor.constant.sums} below).
Since from Figure~\ref{fig_cloud_P_not_B} we see that
each marker $\eta$ has at most one ancestor lying
in the support of $\mu$,
we conclude that we must have $\nu_i=\mu$ for all $i$,
so that the right-hand side of~(\ref{sum_of_nu_i})
can be reduced to a single term.
We would then have $\mu = \beta_{f_1}*\ldots*\beta_{f_{n}}$,
where we can assume that none of the functions $f_i$ are identically 0.
Since $\mu$ charges markers of level 2 but no higher, we must have $n=2$;
and since $\mu$ takes only the values 0 and 1,
we must have $f_1=\ind{\Lambda_1}$ and $f_2=\ind{\Lambda_2}$
for some subsets $\Lambda_1, \Lambda_2 \subseteq \{x,y\}$.
Since $\mu(xxx)=1$,
we must have $x \in \Lambda_1$ and $x \in \Lambda_2$;
likewise, since $\mu(yyy)=1$,
we must have $y \in \Lambda_1$ and $y \in \Lambda_2$;
but then $\mu$ should take the value 0 on the marker $xy$, which it does not.
This proves that $\mu\notin\scrb$.
\qed
\end{example}

\begin{example}
 \label{example_cloud1}
\rm
Let $X=\{x,y\}$, and let $\mu$ be the cloud that takes the value 1
on the markers $xx$, $xyxx$, $xyxy$, $xyy$, $yxx$, $yxy$, $yxx$ and $yyy$
and takes the value 0 elsewhere (Figure~\ref{fig_cloud1}).
This is a cloud of norm 1, supported on levels $\le 3$.
To verify condition (c) of Proposition~\ref{prop.convolution.right},
it suffices to check that for each $\eta$
satisfying $\widetilde{\mu}_{\eta} = 0$,
one has $\widetilde{\mu}_{\eta'} = 0$ for each suffix $\eta'$ of $\eta$.
We leave this verification to the reader.
Hence $\mu \in \scrr$.

On the other hand, $\mu$ does not satisfy the
analogous condition for arbitrary subsequences,
because $\widetilde{\mu}_{xyx} = 0$ while $\widetilde{\mu}_{xx} = 1$.
Hence $\mu \notin \scrp$.
\qed
\end{example}



\begin{proposition}
The class $\scrp$ forms a multiplicative convex cone.
\end{proposition}

\proof
The only nontrivial fact to prove is the stability under convolution.
So let $\mu$ and $\nu$ be clouds in $\scrp$,
let $\eta=(y_0,\ldots,y_k)$ be a marker,
and let $\eta'=(y_{i_0},\ldots,y_{i_r})$
($0\le i_0<\cdots<i_r\le k$) be a subsequence of $\eta$.
For $j=0,\ldots,r$, we set
\be
\sigma_j\ :=\ \eta_0^{y_{i_j}}\ \anc\ \eta  \;.
\ee
Then, for each $\sigma\anc\eta$, we set
\be
\sigma\cap\eta'\ :=\  (y_{i_0},\ldots,y_{i_\ell})  \;,
\ee
where $\ell$ is the largest index such that $\sigma_\ell\anc\sigma$.

Observing that each $\sigma'\anc\eta'$ can be written as 
$\sigma_j\cap\eta'$ for a unique $j$,
we obtain from \reff{eq.cumul.convolution.2}
\begin{equation}
\widetilde{(\mu * \nu)}(\eta') \ = \
  \sum_{j=0}^r \widetilde{\mu}(\sigma_j\cap\eta')
  \left[ \widetilde{\nu}(\eta'\doubleminus(\sigma_j\cap\eta')) - 
\widetilde{\nu}\Bigl((\eta'\doubleminus(\sigma_j\cap\eta'))_-\Bigr)\right].
\end{equation}
Since $\nu \in \scrp \subseteq \scrr$, the square brackets in the 
preceding equation are always nonnegative.
Moreover, since $\sigma_j\cap\eta'$ is a subsequence of $\sigma_j$,
the fact that $\mu$ belongs to $\scrp$ ensures that
\be
\widetilde{\mu}(\sigma_j\cap\eta')\ \le\ \widetilde{\mu}(\sigma_j) \;.
\ee
Hence we get
\begin{eqnarray}
\widetilde{(\mu * \nu)}(\eta')  & \le &
  \sum_{j=0}^r \widetilde{\mu}(\sigma_j)
  \left[ \widetilde{\nu}(\eta'\doubleminus(\sigma_j\cap\eta')) - 
\widetilde{\nu}\Bigl((\eta'\doubleminus(\sigma_j\cap\eta'))_-\Bigr)\right].
\end{eqnarray}
On the other hand, the right-hand side can also be written as
$\widetilde{(\mu' * \nu)}(\eta')$,
where $\mu'$ is any cloud giving mass 
$\widetilde{\mu}(\sigma_j)-\widetilde{\mu}(\sigma_{j-1})$
to the marker $\sigma_j\cap\eta'$ for each $j=0,\ldots,r$.

Using now \reff{eq.cumul.convolution} and again the hypothesis
$\nu\in\scrp$, we get
\begin{eqnarray}
\widetilde{(\mu * \nu)}(\eta') & \le & \widetilde{(\mu' * \nu)}(\eta')
   \nonumber \\[1mm]
 & = & \sum_{j=0}^r \mu'(\sigma_j\cap\eta')  \;
\widetilde{\nu}(\eta'\doubleminus(\sigma_j\cap\eta'))
   \nonumber \\[1mm]
 & \le & \sum_{j=0}^r \mu'(\sigma_j\cap\eta')   \;
\widetilde{\nu}(\eta\doubleminus\sigma_j)
   \nonumber \\
 & = & \sum_{j=0}^r 
\left(\sum_{\sigma_{j-1}\prec\sigma\anc\sigma_j}\mu(\sigma)\right) 
\widetilde{\nu}(\eta\doubleminus\sigma_j)  \;,
  \label{intermediate}
\end{eqnarray}
where $\sigma_{j-1}\prec\sigma$ means $\sigma_{j-1}\anc\sigma$ and 
$\sigma_{j-1}\neq\sigma$, and where we set $\sigma_{-1}:=\emptyset$.
But in each term of~(\ref{intermediate}), we have
$\widetilde{\nu}(\eta\doubleminus\sigma_j) \le 
\widetilde{\nu}(\eta\doubleminus\sigma)$
since $\nu \in \scrp \subseteq \scrr$, hence
\begin{eqnarray}
\widetilde{(\mu * \nu)}(\eta') & \le & 
\sum_{\sigma\anc\sigma_r}\mu(\sigma)  \;
  \widetilde{\nu}(\eta\doubleminus\sigma)
       \nonumber \\[1mm]
& \le & \sum_{\sigma\anc\eta} \mu(\sigma)   \;
\widetilde{\nu}(\eta\doubleminus\sigma)
       \nonumber \\[1mm]
& = & \widetilde{(\mu * \nu)}(\eta)  \;.
\end{eqnarray}
\nopagebreak
\qed

For any cloud $\mu \ge 0$ and any marker $\eta$, let us define
\be
   M_\mu(\eta)  \;\bydef\;
   \sup_{\eta' \colon\, \eta \anc \eta'} \widetilde{\mu}(\eta')
\ee
(this may possibly be $+\infty$);
it is the supremum of the sums of $\mu$
over infinite ascending branches passing through $\eta$.
Obviously
$\widetilde{\mu}(\eta) \le M_\mu(\eta) \le \triplenorm \mu \triplenorm$,
and $M_\mu(\eta)$ is a decreasing function of $\eta$
with respect to the partial order $\anc$.

The clouds belonging to the class $\scrr$ have a remarkable property:

\begin{proposition}
  \label{prop.constant.sums}
If $\mu \ge 0$ belongs to the class $\scrr$
characterized in Proposition~\ref{prop.convolution.right},
then $M_\mu(\eta)$ takes the same value for all markers $\eta$.
\end{proposition}

\begin{corollary}
   \label{cor.constant.sums}
If $\mu \ge 0$ belongs to the class $\scrr$
and is supported on levels $\le N$,
then $\widetilde{\mu}(\eta)$ takes the same value for all markers $\eta$
of level $\ge N$.
\end{corollary}

\proofof{Proposition~\ref{prop.constant.sums}}
Let $\eta_1$ and $\eta_2$ be any two markers.
We can always find a marker $\sigma$ satisfying
$\eta_1\leads\sigma\leads\eta_2$.
Then, for each marker $\eta'$ such that $\eta_2 \anc \eta'$,
we have [using property (c) of Proposition~\ref{prop.convolution.right}]
\be
    \widetilde{\mu}(\eta')  \;\le\;
    \widetilde{\mu}(\eta_1 * \sigma * \eta')   \;\le\;
    M_\mu(\eta_1)  \;,
\ee
hence 
\be
    M_\mu(\eta_2)  \;\le\;  M_\mu(\eta_1)  \;.
\ee
Reversing the roles of $\eta_1$ and $\eta_2$,
we conclude that $M_\mu(\eta_1) = M_\mu(\eta_2)$.
\qed

\begin{definition}
   \label{def_after_prop.constant.sums}
For each (finite) real number $a \ge 0$,
we denote by $\scrs_a$ the class of clouds $\mu \ge 0$
satisfying $M_\mu(\eta) = a$ for all markers $\eta$.
We write $\scrs = \bigcup\limits_{a \ge 0} \scrs_a$.
\end{definition}

\noindent
The class $\scrs_1$ will play a major role in the sequel
(Sections~\ref{sec4.5}--\ref{sec4.8}).

We have just shown that $\scrr \subseteq \scrs$;
let us now show that this inclusion is strict:

\begin{figure}
  \centering
\begin{picture}(0,0)%
\includegraphics{contre_exemple2.pstex}%
\end{picture}%
\setlength{\unitlength}{4144sp}%
\begingroup\makeatletter\ifx\SetFigFont\undefined%
\gdef\SetFigFont#1#2#3#4#5{%
  \reset@font\fontsize{#1}{#2pt}%
  \fontfamily{#3}\fontseries{#4}\fontshape{#5}%
  \selectfont}%
\fi\endgroup%
\begin{picture}(4407,5332)(586,-6042)
\put(2026,-1771){\makebox(0,0)[lb]{\smash{\SetFigFont{10}{12.0}{\rmdefault}{\mddefault}{\updefault}$xy$}}}
\put(1351,-1096){\makebox(0,0)[lb]{\smash{\SetFigFont{10}{12.0}{\rmdefault}{\mddefault}{\updefault}$xxy$}}}
\put(586,-1096){\makebox(0,0)[lb]{\smash{\SetFigFont{10}{12.0}{\rmdefault}{\mddefault}{\updefault}$xxx$}}}
\put(2296,-1096){\makebox(0,0)[lb]{\smash{\SetFigFont{10}{12.0}{\rmdefault}{\mddefault}{\updefault}$xyy$}}}
\put(1853,-1096){\makebox(0,0)[lb]{\smash{\SetFigFont{10}{12.0}{\rmdefault}{\mddefault}{\updefault}$xyx$}}}
\put(4553,-1771){\makebox(0,0)[lb]{\smash{\SetFigFont{10}{12.0}{\rmdefault}{\mddefault}{\updefault}$yy$}}}
\put(1441,-2491){\makebox(0,0)[lb]{\smash{\SetFigFont{10}{12.0}{\rmdefault}{\mddefault}{\updefault}$x$}}}
\put(3968,-2476){\makebox(0,0)[lb]{\smash{\SetFigFont{10}{12.0}{\rmdefault}{\mddefault}{\updefault}$y$}}}
\put(871,-1786){\makebox(0,0)[lb]{\smash{\SetFigFont{10}{12.0}{\rmdefault}{\mddefault}{\updefault}$xx$}}}
\put(3398,-1786){\makebox(0,0)[lb]{\smash{\SetFigFont{10}{12.0}{\rmdefault}{\mddefault}{\updefault}$yx$}}}
\put(2026,-5282){\makebox(0,0)[lb]{\smash{\SetFigFont{10}{12.0}{\rmdefault}{\mddefault}{\updefault}$xy$}}}
\put(1351,-4607){\makebox(0,0)[lb]{\smash{\SetFigFont{10}{12.0}{\rmdefault}{\mddefault}{\updefault}$xxy$}}}
\put(586,-4607){\makebox(0,0)[lb]{\smash{\SetFigFont{10}{12.0}{\rmdefault}{\mddefault}{\updefault}$xxx$}}}
\put(2296,-4607){\makebox(0,0)[lb]{\smash{\SetFigFont{10}{12.0}{\rmdefault}{\mddefault}{\updefault}$xyy$}}}
\put(1853,-4607){\makebox(0,0)[lb]{\smash{\SetFigFont{10}{12.0}{\rmdefault}{\mddefault}{\updefault}$xyx$}}}
\put(4553,-5282){\makebox(0,0)[lb]{\smash{\SetFigFont{10}{12.0}{\rmdefault}{\mddefault}{\updefault}$yy$}}}
\put(3878,-4607){\makebox(0,0)[lb]{\smash{\SetFigFont{10}{12.0}{\rmdefault}{\mddefault}{\updefault}$yxy$}}}
\put(3113,-4607){\makebox(0,0)[lb]{\smash{\SetFigFont{10}{12.0}{\rmdefault}{\mddefault}{\updefault}$yxx$}}}
\put(4823,-4607){\makebox(0,0)[lb]{\smash{\SetFigFont{10}{12.0}{\rmdefault}{\mddefault}{\updefault}$yyy$}}}
\put(4380,-4607){\makebox(0,0)[lb]{\smash{\SetFigFont{10}{12.0}{\rmdefault}{\mddefault}{\updefault}$yyx$}}}
\put(1441,-6002){\makebox(0,0)[lb]{\smash{\SetFigFont{10}{12.0}{\rmdefault}{\mddefault}{\updefault}$x$}}}
\put(3968,-5987){\makebox(0,0)[lb]{\smash{\SetFigFont{10}{12.0}{\rmdefault}{\mddefault}{\updefault}$y$}}}
\put(871,-5297){\makebox(0,0)[lb]{\smash{\SetFigFont{10}{12.0}{\rmdefault}{\mddefault}{\updefault}$xx$}}}
\put(3398,-5297){\makebox(0,0)[lb]{\smash{\SetFigFont{10}{12.0}{\rmdefault}{\mddefault}{\updefault}$yx$}}}
\put(2469,-2521){\makebox(0,0)[lb]{\smash{\SetFigFont{12}{14.4}{\rmdefault}{\mddefault}{\updefault}$\mu$}}}
\put(2469,-5986){\makebox(0,0)[lb]{\smash{\SetFigFont{12}{14.4}{\rmdefault}{\mddefault}{\updefault}$\rho^1*\mu$}}}
\end{picture}
\vspace{5mm}
\caption{
   The clouds $\mu$ and $\rho^1*\mu$ in Example~\ref{example_cloud2}
   take the value 1 (resp.\ 0) on the markers
   indicated by a full (resp.\ empty) circle.
   The cloud $\mu$ belongs to $\scrs_1$ but not to $\scrr$.
}
  \label{fig:contre-exemple 2}
\end{figure}

\begin{example}
  \label{example_cloud2}
\rm
Let $X=\{x,y\}$, and let $\mu$ be the cloud that takes the value 1
on the markers $xx$, $xyx$, $xyy$  and $y$,
and takes the value 0 elsewhere (Figure~\ref{fig:contre-exemple 2}).
Then we have $M_\mu(\eta) = 1$ for every marker $\eta$,
but $\rho^1*\mu\not\audessus\mu$
(consider the cumulative distributions on the marker $xy$). 
So $\mu \in \scrs_1$ but $\mu \notin \scrr$.
\qed
\end{example}


\begin{proposition}
  \label{prop.scrs.cone}
The class $\scrs$ forms a multiplicative convex cone.
More specifically:
\begin{itemize}
   \item[(a)]  If $\mu \in \scrs_a$ and $\nu \in \scrs_b$
      and $s,t \ge 0$, then $s\mu + t\nu \in \scrs_{sa+tb}$.
   \item[(b)]  If $\mu \in \scrs_a$ and $\nu \in \scrs_b$,
      then $\mu * \nu \in \scrs_{ab}$.
\end{itemize}
\end{proposition}

\proof
(a) Given $\mu \in \scrs_a$, $\nu \in \scrs_b$ and $s,t \ge 0$,
we clearly have
$M_{s\mu + t\nu}(\eta) \le sa+tb$
for any marker $\eta$.
The reverse inequality is easily obtained by choosing
$\eta' \ancreverse \eta$ such that 
$\sum_{\sigma\anc\eta'} s\mu(\sigma) \ge sa-\epsilon$
and then choosing $\eta'' \ancreverse \eta'$ such that
$\sum_{\sigma\anc\eta''} t\nu(\sigma) \ge tb-\epsilon$.

(b) Equation \reff{eq.cumul.convolution} gives, for any marker $\eta$,
\be
   \widetilde{(\mu * \nu)}(\eta)
   \ \le \ b\,\sum_{\sigma \anc \eta}  \mu(\sigma) \ \le \ ab  \;.
\ee
For the reverse inequality, we first choose
$\overline{\eta} \ancreverse \eta$ such that
$\sum_{\sigma\anc\overline \eta} \mu(\sigma) \ge a-\epsilon$.
Observe next using \reff{eq.cumul.convolution}
that for any $\eta' \ancreverse \overline{\eta}$, we have
\be
    \widetilde{(\mu * \nu)}(\eta')\ \ge\ 
   \sum_{\sigma \anc \overline\eta} \mu(\sigma) \,
   \widetilde{\nu}(\eta'\doubleminus\sigma)   \;.
 \label{eq.proof.scrs}
\ee
Now, if $\hbox{level}(\overline{\eta}) = k$,
we construct inductively a sequence
$\overline{\eta} \anc \eta'_0 \anc \eta'_1 \anc \ldots \anc \eta'_k$
such that
$\widetilde{\nu}(\eta'_r\doubleminus\sigma) \ge b-\epsilon$
for any $\sigma \anc \overline\eta$ with $\hbox{level}(\sigma) \le r$.
Taking $\eta' = \eta'_k$ in \reff{eq.proof.scrs} gives the desired result.
\qed

In summary, we have introduced four natural classes of clouds,
which are multiplicative convex cones and
satisfy $\scrb \subsetneq \scrp \subsetneq \scrr \subsetneq \scrs$.

%
%


\subsection{A fundamental inequality}   \label{sec4.4}

Let us now substantiate our assertion that the
partial ordering $\audessus$ is related to efficiency of cleaning.

\begin{proposition}[Fundamental comparison inequality]
 \label{prop.fundamental}
\FH
Let $\mu$ and $\nu$ be clouds of finite norm.
If $\mu \audessus \nu$, then 
\be
   T_\mu w \;\le\;  T_\nu w   \;.
 \label{fundamental.inequality.1}
\ee
If in addition $\mu,\nu \ge 0$,
then one has
\be
   \| T_\mu \|_{w \to w}  \;\le\;  \| T_\nu \|_{w \to w}
 \label{fundamental.inequality.operator}
\ee
and, for every vector $c \ge 0$,
\be
   \| c T_\mu \|_w  \;\le\;  \| c T_\nu \|_w   \;.
 \label{fundamental.inequality.vector}
\ee
\end{proposition}

\proof
Because $\mu$ and $\nu$ have finite norm,
$T_\mu w$ (resp.\ $T_\nu w$) is the pointwise limit of
$T_{\mu_N}w$ (resp.\ $T_{\nu_N}w$) as $N \to \infty$,
where for each $N\ge 0$ the cloud $\mu_N$ (resp.\ $\nu_N$)
is supported on levels $\le N$ and
coincides with $\mu$ (resp. $\nu$) on these levels.
Thus it is enough to prove the result
when both $\mu$ and $\nu$ are supported on levels $\le N$.
We shall do this by induction on $N$. 

We clearly have (\ref{fundamental.inequality.1}) if $N=0$.
So let  $N\ge1$, and assume that (\ref{fundamental.inequality.1})
holds whenever $\mu$ and $\nu$ are supported on levels $\le N-1$.
Now let $\mu$ and $\nu$ be clouds supported on levels $\le N$
with $\mu\audessus\nu$. We can suppose that 
\begin{equation}
  \label{eq:hyp1}
  \nu(\eta)=0 \,\mbox{ whenever }\, \mbox{level}(\eta)=N  \;;
\end{equation} 
for if this is not the case, we can simply replace $\mu$ and $\nu$
with the clouds obtained by subtracting $\nu(\eta)$ from
both $\mu(\eta)$ and $\nu(\eta)$ for each marker $\eta$ of level $N$.
Furthermore, we can also suppose that
\begin{equation}
  \label{eq:hyp2}
  \mu(\eta) \ge 0 \,\mbox{ whenever }\, \mbox{level}(\eta)=N  \;;
\end{equation} 
otherwise, for each marker $\eta$ of level $N$ such that $\mu_\eta<0$,
we can replace $\mu_\eta$ with 0:
the new cloud obtained in this way is still $\audessus\nu$.

Next,  we consider the cloud $\mu'$ defined by
\be
 \mu'(\eta)   \;=\;
    \cases{ 0  & if $\hbox{level}(\eta) = N$,       \cr
            \mu(\eta) + \sup_{x_N\in X}\mu({\eta \circ x_N})
               & if $\hbox{level}(\eta)=N-1$,       \cr
            \mu(\eta)    & otherwise.               \cr
          }
\ee
Note that the definition of $\mu'$ is the same as \reff{def_of_nuprime}
in the proof of Proposition~\ref{prop.Tnu.norm}, 
and a similar calculation using the Fundamental Hypothesis gives
\be
    T_{\mu} w \;\le\; T_{\mu'}w  \;.
\ee
Now, $\mu'$ and $\nu$ are supported on levels $\le N-1$,
and the induction hypothesis gives
\be
    T_{\mu'} w \;\le\;  T_{\nu}w  \;.
\ee
Together these prove (\ref{fundamental.inequality.1}).

Finally, \reff{fundamental.inequality.operator}
and \reff{fundamental.inequality.vector}
are easy consequences of (\ref{fundamental.inequality.1})
in the case of nonnegative clouds $\mu,\nu$
and, for the latter, nonnegative dirt vectors $c$.
\qed


\subsection{Clouds carried by $\Lambda$}   \label{sec4.5}

The following definition will play a fundamental role in our analysis:

\begin{definition}
Let $\Lambda \subseteq X$.
We say that a cloud $\mu$ is {\em carried by $\Lambda$}\/
in case $\mu_\eta = 0$ for every marker $\eta = (x_0,x_1,\ldots,x_k)$
having at least one index $j < k$ with $x_j \in \Lambda^c$.
\end{definition}

Physically, this means that $\mu$, interpreted as a cleaning operator
(i.e.\ as acting by convolution on the right), never sends dirt
back into $\Lambda$ from outside $\Lambda$.

Please note that for clouds $\mu$ carried by $\Lambda$,
we have $M_\mu(\eta) = \widetilde{\mu}(\eta)$
whenever $\eta$ has at least one entry outside $\Lambda$.

\begin{lemma}
  \label{lemma.clouds.Lambda}
Let $\Lambda \subseteq X$.  Then:
\begin{itemize}
   \item[(a)]  All clouds $I_f$ are carried by $\Lambda$.
   \item[(b)]  The cloud $\beta_f$ is carried by $\Lambda$
       if and only if $\supp(f) \subseteq \Lambda$.
   \item[(c)]  The cloud $\pi_{\Lambda'}$ is carried by $\Lambda$
       if and only if either $\Lambda' \subseteq \Lambda$ or
       $\Lambda' = X$.
   \item[(d)]  If $\mu$ and $\nu$ are carried by $\Lambda$ and $a,b \in \R$,
       then $a\mu + b\nu$ is carried by $\Lambda$.
   \item[(e)]  If $\mu$ and $\nu$ are carried by $\Lambda$,
       then $\mu * \nu$ is carried by $\Lambda$.
\end{itemize}
\end{lemma}

\proof
(a) The cloud $I_f$ is nonvanishing only on markers of level 0,
and so is manifestly carried by $\Lambda$ for any $\Lambda$.

(b) The cloud $\beta_f$ is nonvanishing only on markers of levels 0 and 1.
For the latter we have $(\beta_f)_{(x_0,x_1)} = f(x_0)$.
It follows that $\beta_f$ is carried by $\Lambda$
if and only if $f(x) = 0$ for all $x \in \Lambda^c$.

(c) follows easily from the definitions.

(d) is trivial.

(e) Suppose that $(\mu * \nu)_{(x_0,\ldots,x_k)} \neq 0$.
Then there must exist an index $j$ ($0 \le j \le k$) such that
$\mu_{(x_0,\ldots,x_j)} \neq 0$ and $\nu_{(x_j,\ldots,x_k)} \neq 0$.
But since $\mu$ and $\nu$ are carried by $\Lambda$,
we must have $x_0,\ldots,x_{j-1} \in \Lambda$ and
$x_j,\ldots,x_{k-1} \in \Lambda$.
\qed

Clouds carried by $\Lambda$ satisfy an identity analogous to
Lemma~\ref{lemma_betagi}:

\begin{lemma}
   \label{lemma_betagi.clouds}
Let $\Lambda \subseteq X$,
and let $\mu_1,\ldots,\mu_n$ be clouds carried by $\Lambda$.
Then the cloud $\mu_1 * \ldots * \mu_n * I_\Lambda$
is supported on markers having all their entries in $\Lambda$.
Moreover, if $h_i \colon\, X \to \R$
are functions satisfying $h_i \restrict \Lambda \equiv 1$,
then
\be
   \mu_1 * \ldots * \mu_n * I_\Lambda  \;\,=\;\,
   I_{h_1} * \mu_1 * \ldots * I_{h_n} * \mu_n * I_\Lambda  \;.
\ee
\end{lemma}

\proof
By Lemma~\ref{lemma.clouds.Lambda},
the cloud $\mu_1 * \ldots * \mu_n$ is carried by $\Lambda$.
It follows that $\mu_1 * \ldots * \mu_n * I_\Lambda$
can be nonzero only on markers $\eta = (x_0,\ldots,x_k)$
with all $x_i \in \Lambda$ ($0 \le i \le k$).
For such $\eta$, the presence of the factors $I_{h_i}$
changes nothing.
\qed

We also have the following property concerning convolutions with
$\pi_\Lambda$:

\begin{lemma}
   \label{lemma.mu_star_pilambda}
Let $\Lambda \subseteq X$, and let $\mu \ge 0$ be a cloud.
Then the following are equivalent:
\begin{itemize}
   \item[(a)]  $\mu$ is carried by $\Lambda$.
   \item[(b)]  $\mu * \pi_\Lambda$ is carried by $\Lambda$.
\end{itemize}
\end{lemma}

\proof
(a) $\implies$ (b) follows from Lemma~\ref{lemma.clouds.Lambda}(c,e).

(b) $\implies$ (a):  Suppose that $\mu$ is not carried by $\Lambda$,
i.e.\ that there exists a marker $\eta = (x_0,x_1,\ldots,x_k)$
and an index $j < k$ such that $x_j \in \Lambda^c \neq \emptyset$
and $\mu(\eta) > 0$.
If $x_k \in \Lambda^c$, we have
$(\mu * \pi_\Lambda)(\eta) \ge \mu(\eta) > 0$;
if $x_k \in \Lambda$, we have
$(\mu * \pi_\Lambda)(\eta') \ge \mu(\eta) > 0$
for $\eta' = (x_0,x_1,\ldots,x_k,y)$ with any $y \in \Lambda^c$;
either way we conclude that $\mu * \pi_\Lambda$ is not carried by $\Lambda$.
\qed

\medskip

{\bf Remark.}  The implication (b) $\implies$ (a) is false
if $\mu$ is not assumed nonnegative.
To see this, consider $X = \{x,y\}$ and $\Lambda = \{x\}$,
and set $\mu(yx) = 1$, $\mu(y x^k y) = -1$ for all $k \ge 1$,
and $\mu = 0$ on all other markers.
Then it is not hard to verify that $\mu * \pi_\Lambda = 0$.
But $\mu$ is not carried by $\Lambda$.
\qed

\bigskip

For certain pairs of clouds carried by $\Lambda$,
we can prove an inequality
going in the {\em opposite}\/ direction
to Proposition~\ref{prop.fundamental},
provided that we look only at markers ending {\em outside}\/ $\Lambda$;
moreover, this inequality holds pointwise.
Let us recall that $\partial\Lambda$ denotes
the set of markers of the form $\eta = (x_0,\ldots,x_k)$ [$k \ge 0$]
for which $x_0,\ldots,x_{k-1} \in \Lambda$ and $x_k \in \Lambda^c$.

\begin{proposition}
   \label{prop.reverse.clouds}
Let $\Lambda \subseteq X$,
and let $\mu, \nu \ge 0$ be clouds carried by $\Lambda$.
Suppose further that $\widetilde{\mu}(\eta) \le \widetilde{\nu}(\eta)$
for all markers $\eta \in \partial\Lambda$.
Then the following are equivalent:
\begin{itemize}
   \item[(a)] $\mu \audessus \nu$
   \item[(b)] $\mu * I_f  \,\audessus\,  \nu * I_f$
      for every $f$ satisfying $\chi_\Lambda \le f \le 1$.
   \item[(c)] $\mu * I_\Lambda  \,\audessus\,  \nu * I_\Lambda$
   \item[(d)] $\widetilde{\mu}(\eta) \,\le\, \widetilde{\nu}(\eta)$
      for every marker $\eta$ having all its elements inside $\Lambda$.
\end{itemize}
Moreover, if $\Lambda \neq X$ and
$\widetilde{\mu}(\eta) = \widetilde{\nu}(\eta)$
for all $\eta \in \partial\Lambda$,
then (a)--(d) are equivalent to
\begin{itemize}
   \item[(e)] $\mu * I_{\Lambda^c}  \,\ge\,  \nu * I_{\Lambda^c}$
\end{itemize}
\end{proposition}

\proof
Since $\mu$ and $\nu$ are carried by $\Lambda$,
the hypothesis implies that
$\widetilde{\mu}(\eta) \le \widetilde{\nu}(\eta)$
for every marker $\eta$ having at least one element outside $\Lambda$.
It follows that (a) is equivalent to (d).

On the other hand, the inequalities expressing (a)--(c) are all identical
when evaluated on markers $\eta$ having all their elements in $\Lambda$:
they simply assert (d).  So we have (b) $\implies$ (c) $\implies$ (d),
and it suffices now to show that (a) implies the inequality (b)
when evaluated on markers $\eta$ having at least one element outside $\Lambda$.
So let $\eta$ be such a marker,
and let $\sigma$ be its unique ancestor in $\partial\Lambda$.
Then
\begin{subeqnarray}
   \widetilde{(\mu * I_f)}(\eta)
   \;=\;
   \widetilde{(\mu * I_f)}(\sigma)
   & = &
   \widetilde{\mu}(\sigma^-) \,+\, f(\hbox{last}(\sigma)) \, \mu(\sigma)
        \\[1mm]
   & = &
   [1 - f(\hbox{last}(\sigma))] \, \widetilde{\mu}(\sigma^-)
   \,+\,
   f(\hbox{last}(\sigma)) \, \widetilde{\mu}(\sigma)
   \qquad\quad
\end{subeqnarray}
and likewise for $\nu$.
[Recall that if $\eta = (x_0,\ldots,x_k)$,
 then $\eta^- \bydef (x_0,\ldots,x_{k-1})$.]
It follows,
by taking a convex combination of inequalities,
that (a) implies (b).

Since $\mu$ is carried by $\Lambda$,
the cloud $\mu * I_{\Lambda^c}$ is supported on $\partial\Lambda$,
and for $\eta \in \partial\Lambda$ we have
\be
   (\mu * I_{\Lambda^c})(\eta)
   \;=\;
   \widetilde{\mu}(\eta) \,-\, \widetilde{\mu}(\eta^-)
\ee
and likewise for $\nu$.
If $\widetilde{\mu}(\eta) = \widetilde{\nu}(\eta)$
for all $\eta \in \partial\Lambda$,
then (e) is equivalent to
\be
   \widetilde{\mu}(\eta^-) \;\le\; \widetilde{\nu}(\eta^-)
\ee
for all $\eta \in \partial\Lambda$.
But since $\Lambda \neq X$, every marker $\sigma$ with all its elements
in $\Lambda$ is of the form $\eta^-$ for some $\eta \in \partial\Lambda$,
so this is equivalent to (d).
\qed

\begin{corollary}
   \label{cor.reverse.clouds}
Let $\Lambda \subseteq X$,
and let $\mu, \nu \ge 0$ be clouds carried by $\Lambda$
and belonging to the class $\scrs_a$
(cf.\ Definition~\ref{def_after_prop.constant.sums})
for the same constant $a$.
If $\mu \audessus \nu$, then
\begin{itemize}
   \item[(i)]  $\mu * I_{\Lambda^c}  \,\ge\,  \nu * I_{\Lambda^c}$
   \item[(ii)] $\mu * I_f  \,\audessus\,  \nu * I_f$
      for every $f$ satisfying $\chi_\Lambda \le f \le 1$.
\end{itemize}
\end{corollary}

\noindent
This is an immediate consequence of Proposition~\ref{prop.reverse.clouds}
together with the definition of $\scrs_a$,
which entails that $M_\mu(\eta) = M_\nu(\eta) = a$ for all markers $\eta$.
Corollary~\ref{cor.reverse.clouds} applies in particular (with $a=1$)
if $\mu$ and $\nu$ are of the form
$\beta_{h_1} * \ldots * \beta_{h_n}$
with all $\supp(h_i) \subseteq \Lambda$.

%
%

\subsection{$\Lambda$-regular clouds}   \label{sec4.6}

  
\begin{definition}
 \label{def.lambdaregular}
We say that cloud $\mu \ge 0$ is {\em $\Lambda$-regular}\/
in case it satisfies the following two conditions:
\begin{subeqnarray}
   \widetilde{\mu}(\eta)
   & \le &  1
     \quad\hbox{for all markers}\; \eta
     \quad\hbox{[i.e.\ $\triplenorm \mu \triplenorm \le 1$]}
           \slabel{def.Lambdaregular.1}  \\[2mm]
   \widetilde{\mu}(\eta)
   & = &   1
     \quad\hbox{for all}\; \eta \in \partial\Lambda
           \slabel{def.Lambdaregular.2}
           \label{def.Lambdaregular.all}
\end{subeqnarray}
\end{definition}


\begin{lemma}
   \label{lemma.equivalent_conditions_Lambda-regularity}
Let $\Lambda \subseteq X$, and let $\mu \ge 0$ be a cloud.
\begin{itemize}
   \item[(a)]  If $\Lambda \neq X$, then $\mu$ is $\Lambda$-regular
      if and only if it is carried by $\Lambda$ and
      belongs to the class $\scrs_1$.
   \item[(b)]  If $\Lambda = X$, then $\mu$ is $\Lambda$-regular
      if and only if $\triplenorm \mu \triplenorm \le 1$.
\end{itemize}
\end{lemma}

\proof
(a)  If $\mu \in \scrs_1$, then \reff{def.Lambdaregular.1} clearly holds.
If, in addition, $\mu$ is carried by $\Lambda$,
we have $\widetilde{\mu}(\eta) = M_\mu(\eta) = 1$
whenever $\eta$ has at least one element outside $\Lambda$,
so that \reff{def.Lambdaregular.2} holds.

Conversely, suppose that $\mu$ is $\Lambda$-regular and that $\Lambda \neq X$.
If $\eta$ has all its entries in $\Lambda$,
then it is a proper ancestor of some $\sigma \in \partial\Lambda$
(since $\Lambda \neq X$),
and we have $M_\mu(\eta) \ge \widetilde{\mu}(\sigma) = 1$.
If $\eta$ has at least one entry outside $\Lambda$,
then it has a (unique) marker $\sigma \in \partial\Lambda$ as an ancestor,
in which case $M_\mu(\eta) \ge \widetilde{\mu}(\eta) \ge
               \widetilde{\mu}(\sigma) = 1$.
On the other hand, by \reff{def.Lambdaregular.1}
we always have $M_\mu(\eta) \le 1$,
hence $M_\mu(\eta) = 1$ for all $\eta$ and thus $\mu \in \scrs_1$.
Finally, if $\eta$ has some proper ancestor $\sigma \in \partial\Lambda$,
then (\ref{def.Lambdaregular.all}a,b) imply
$1 \ge \widetilde{\mu}(\eta) \ge \mu(\eta) + \widetilde{\mu}(\sigma)
   = \mu(\eta) + 1$, hence $\mu(\eta) = 0$,
so that $\mu$ is carried by $\Lambda$.
%
%

(b) If $\Lambda = X$, the condition \reff{def.Lambdaregular.2} is empty.
\qed

\begin{corollary}
   \label{cor.Lambda-regularity}
Let $\Lambda \subseteq X$.  Then:
\begin{itemize}
   \item[(a)]  The cloud $\rho^0$ (indicator of level 0) is $\Lambda$-regular.
   \item[(b)]  The cloud $\beta_f$ is $\Lambda$-regular
       if and only if $\supp(f) \subseteq \Lambda$.
   \item[(c)]  The cloud $\pi_{\Lambda'}$ is $\Lambda$-regular
       if and only if $\Lambda' \subseteq \Lambda$.
   \item[(d)]  If $\mu$ and $\nu$ are $\Lambda$-regular,
       then $\mu * \nu$ is $\Lambda$-regular.
\end{itemize}
\end{corollary}

\proof
This is an immediate consequence of
Lemma~\ref{lemma.equivalent_conditions_Lambda-regularity}
along with Proposition~\ref{prop.scrs.cone}
and Lemma~\ref{lemma.clouds.Lambda};
for part (c) one has to think separately about the cases in which
$\Lambda$ and/or $\Lambda'$ is equal or not equal to $X$.
\qed

\begin{corollary}
   \label{cor.Lambda-regularity.2}
Let $\Lambda \subseteq X$, and let $\mu$ be a $\Lambda$-regular cloud.
Then
\begin{eqnarray}
   I_{\Lambda^c} * \mu           & = &  I_{\Lambda^c}      \label{eq.lr2.1} \\
   I_{\Lambda} * \mu * I_\Lambda & = &  \mu *  I_{\Lambda} \label{eq.lr2.2}
\end{eqnarray}
\end{corollary}

\proof
We begin by proving \reff{eq.lr2.1}.
The claim is trivial if $\Lambda = X$, so assume $\Lambda \neq X$.
By Lemma~\ref{lemma.equivalent_conditions_Lambda-regularity},
$\mu$ is carried by $\Lambda$,
so it must vanish on all markers of the form
$(x_0,\ldots,x_k)$ with $x_0 \in \Lambda^c$ and $k \ge 1$.
But then, since $\mu$ belongs to $\scrs_1$,
it must give mass 1 to each marker of the form $(x_0)$
with $x_0 \in \Lambda^c$.

\reff{eq.lr2.2} follows from \reff{eq.lr2.1}
by convoluting on the right with $I_\Lambda$ and rearranging.
\qed

\begin{corollary}
   \label{cor.reverse.clouds_bis}
Let $\Lambda \subseteq X$,
and let $\mu, \nu$ be $\Lambda$-regular clouds.
If $\mu \audessus \nu$, then
\begin{itemize}
   \item[(i)]  $\mu * I_{\Lambda^c}  \,\ge\,  \nu * I_{\Lambda^c}$
   \item[(ii)] $\mu * I_f  \,\audessus\,  \nu * I_f$
      for every $f$ satisfying $\chi_\Lambda \le f \le 1$.
\end{itemize}
\end{corollary}

\proof
If $\Lambda \neq X$, this is an immediate consequence of
Lemma~\ref{lemma.equivalent_conditions_Lambda-regularity}
and Corollary~\ref{cor.reverse.clouds}.
If $\Lambda = X$, the claims are trivially true.
\qed

\subsection{Comparison of $T_\mu$ with $\Pi_\Lambda$}  \label{sec4.7}



We recall that the cloud $\pi_\Lambda$,
which is the indicator of $\partial\Lambda$,
is associated to the operator $\Pi_\Lambda$.
This cloud plays a special role among $\Lambda$-regular clouds,
by virtue of its minimality with respect to $\audessus$:

\begin{lemma}
   \label{lemma.pilambda}
Let $\Lambda \subseteq X$, and let $\mu$ be a $\Lambda$-regular cloud.
Then:
\begin{itemize}
   \item[(a)]
       $\mu * I_{\Lambda^c} \le \pi_\Lambda \audessus \mu \audessus \rho^0$.
   \item[(b)] If $\Lambda \neq X$ and $\mu$ is supported on levels $\le N$,
      then $\pi_\Lambda \audessus \beta_\Lambda^N \audessus \mu$.
\end{itemize}
\end{lemma}

\proof

(a) Let $\mu$ be any $\Lambda$-regular cloud and
let $\eta$ be any marker.
Then either there exists $\sigma\anc\eta$ with $\sigma\in\partial\Lambda$,
in which case $\widetilde{\mu}(\eta)$ and
$\widetilde{\pi_\Lambda}(\eta)$ are both equal to 1;
or else there does not exist such a $\sigma$,
in which case $\widetilde{\pi_\Lambda}(\eta) = 0 \le \widetilde{\mu}(\eta)$.
This proves that $\pi_\Lambda\audessus\mu$.
Then,
Corollary~\ref{cor.reverse.clouds_bis}
gives
$\mu * I_{\Lambda^c} \le \pi_\Lambda * I_{\Lambda^c} = \pi_\Lambda$.
The inequality $\mu \audessus \rho^0$ is trivial.

(b)  Corollary~\ref{cor.Lambda-regularity}(b,d) shows that
$\beta_\Lambda^N$ is $\Lambda$-regular,
so the relation $\pi_\Lambda \audessus \beta_\Lambda^N$ follows from (a).
For the second inequality, observe first that
for each marker $\eta$ of level $\le N$,
we have $\beta_\Lambda^N(\eta)=\pi_\Lambda(\eta)$,
hence
\be
\widetilde{\beta_\Lambda^N}(\eta)  \;=\;
 \widetilde{\pi_\Lambda}(\eta) \;\le \;
 \widetilde{\mu}(\eta)  \;,
\ee
where the last inequality again uses (a).
On the other hand, if level$(\eta)>N$, since $\beta_\Lambda^N$ and $\mu$
are $\Lambda$-regular clouds supported on levels $\le N$, we always have
\be
\widetilde{\beta_\Lambda^N}(\eta)\; =\; 1\; =\; \widetilde{\mu}(\eta)  \;,
\ee
where the last inequality uses $\mu \in \scrs_1$
from Lemma~\ref{lemma.equivalent_conditions_Lambda-regularity}(a).
This proves that $\beta_\Lambda^N \audessus \mu$.
\qed

\begin{lemma}
   \label{lemma.pilambda.2}
Let $\Lambda \subseteq X$, and let $\mu \ge 0$ be a cloud.
Consider the following statements:
\begin{itemize}
   \item[(a)] $\mu$ is $\Lambda$-regular.
   \item[(b)] $\mu * \pi_\Lambda = \pi_\Lambda$.
   \item[(c)] $\pi_\Lambda * \mu = \pi_\Lambda$.
\end{itemize}
Then (a) implies both (b) and (c);
and if $\Lambda \neq X$, then (b) implies (a) as well.
\end{lemma}

\proof
Suppose first that $\mu$ is $\Lambda$-regular.
Then by Lemma~\ref{lemma.pilambda},
we have $\pi_\Lambda \audessus \mu \audessus \rho^0$.
Since $\pi_\Lambda \in \scrr$, we can right-convolute this inequality
with $\pi_\Lambda$
to obtain $\pi_\Lambda \audessus \mu * \pi_\Lambda \audessus \pi_\Lambda$
(since $\pi_\Lambda * \pi_\Lambda = \pi_\Lambda$), thus proving (b).
Likewise, since $\pi_\Lambda \ge 0$, we can left-convolute with $\pi_\Lambda$
to prove (c).

Now suppose that $\mu * \pi_\Lambda = \pi_\Lambda$.
By Lemma~\ref{lemma.mu_star_pilambda} we deduce that
$\mu$ is carried by $\Lambda$.
Moreover, if $\eta \in \partial\Lambda$, we have
\be
   1 \;=\;  \pi_\Lambda(\eta)  \;=\;  (\mu * \pi_\Lambda)(\eta)
     \;=\;  \sum_{\sigma \ancestor \eta}
                \mu(\sigma) \, \pi_\Lambda(\eta \doubleminus \sigma)
     \;=\;  \sum_{\sigma \ancestor \eta}
                \mu(\sigma)
     \;=\;  \widetilde{\mu}(\eta)
     \;.
\ee
When $\Lambda \neq X$, it easily follows that $\mu$ is $\Lambda$-regular.
%
%
\qed

\medskip

{\bf Remark.}  The implication (c) $\implies$ (a) is false
when $\Lambda \neq \emptyset$:
consider, for instance, $\mu = I_{\Lambda^c}$.
Indeed, (c) does not even imply that $\mu$ is carried by $\Lambda$:
consider, for instance,
$\mu = I_{\Lambda^c} + I_\Lambda * \rho^1 * I_{\Lambda^c} * \rho^1$.
\qed

\bigskip

We now examine the deviation $\mu - \pi_\Lambda$,
and prove analogues of Proposition~\ref{prop.convergence-to-balayage}
and Lemma~\ref{lemma.A_minus_Pilambda}:

\begin{lemma}[Comparison with \mbox{\protect\boldmath $\pi_\Lambda$}]
   \label{lemma.mu_minus_pilambda}
\FH
Let $\mu$ be a $\Lambda$-regular cloud.  Then:
\begin{itemize}
   \item[(a)]  The cloud $\mu - \pi_\Lambda$ can be decomposed in the form
\begin{subeqnarray}
   \mu - \pi_\Lambda
   & = &
   (\mu - \pi_\Lambda) * I_\Lambda  \;+\;
      (\mu - \pi_\Lambda) * I_{\Lambda^c}      \\[2mm]
   & = &
   \mu * I_\Lambda  \;-\;
      (\pi_\Lambda \,-\, \mu * I_{\Lambda^c})  \\[2mm]
   & = &
   \mu * I_\Lambda  \;-\;  \mu * I_\Lambda * \pi_\Lambda
 \label{eq1.lemma.mu_minus_pilambda}
\end{subeqnarray}
where
\begin{itemize}
   \item[(i)]  $\mu * I_\Lambda$ is nonnegative and is supported
     on markers all of whose entries are in $\Lambda$; and
   \item[(ii)] $\pi_\Lambda \,-\, \mu * I_{\Lambda^c}$ is nonnegative,
     is $\audessus \mu * I_\Lambda$, and is supported
     on markers in $\partial\Lambda$.
\end{itemize}
   \item[(b)]  For any vector $c \ge 0$, we have
\be
   \| c (T_\mu - \Pi_\Lambda) \|_w
   \;=\;
   \| c T_\mu I_\Lambda \|_w  \,+\,
   \| c (\Pi_\Lambda - T_\mu I_{\Lambda^c}) \|_w
   \;\le\;
   2 \| c T_\mu I_\Lambda \|_w
   \;.
     \label{eq.mu-pilambda.eqineq}
\ee
   \item[(c)]  We have
\be
  \|T_\mu - \Pi_\Lambda\|_{w\to w}\; \le\; 2\|T_\mu I_\Lambda \|_{w\to w}
   \;.
\ee
\end{itemize}
\end{lemma}

\proof
(a) The equalities (\ref{eq1.lemma.mu_minus_pilambda}a,b) are trivial,
and (\ref{eq1.lemma.mu_minus_pilambda}c) follows by using
Lemma~\ref{lemma.pilambda.2}:
$\pi_\Lambda - \mu * I_{\Lambda^c} =
 \mu * \pi_\Lambda - \mu * I_{\Lambda^c} * \pi_\Lambda =
 \mu * I_\Lambda * \pi_\Lambda$.
The inequalities
$0 \le \pi_\Lambda - \mu * I_{\Lambda^c} \audessus \mu * I_\Lambda$
are an immediate consequence of Lemma~\ref{lemma.pilambda}(a).

(b)  The equality is an immediate consequence of the
identity, sign and support properties from part (a).
The inequality comes from
$0 \le \pi_\Lambda - \mu * I_{\Lambda^c} \audessus \mu * I_\Lambda$
together with Proposition~\ref{prop.fundamental}.

(c) For a general vector $c$, we can always write
\be
   \| c (T_\mu - \Pi_\Lambda) \|_w
   \;\le\;
   \| \, |c| \, |T_\mu - \Pi_\Lambda| \, \|_w
   \;=\;
   \| \, |c| \, T_\mu I_\Lambda \|_w  \,+\,
   \| \, |c| \, (\Pi_\Lambda - T_\mu I_{\Lambda^c}) \|_w   \;,
\ee
and the result then follows from part (b).
\qed

\begin{proposition}[Comparison of cleaners]
  \label{comparison-with-Pi-Lambda}
\FH
Let $\mu$ and $\nu$ be $\Lambda$-regular clouds,
with $\mu\audessus\nu$. Then:
\begin{itemize}
   \item[(a)]  We have
\begin{eqnarray}
   \mu * I_\Lambda      & \audessus & \nu * I_\Lambda   \\[2mm]
   \mu * I_{\Lambda^c}  & \ge  &      \nu * I_{\Lambda^c}   \\[2mm]
   \pi_\Lambda - \mu * I_{\Lambda^c}  & \le  &
        \pi_\Lambda - \nu * I_{\Lambda^c}
\end{eqnarray}
   \item[(b)]  For any vector $c \ge 0$, we have
\begin{eqnarray}
   \| c (T_\mu - \Pi_\Lambda) \|_w
   & \le &
   \| c (T_\nu - \Pi_\Lambda) \|_w
          \label{closer-to-Pi-Lambda} \\[2mm]
   c T_\mu I_{\Lambda^c}
   & \ge &
   c T_\nu I_{\Lambda^c}
     \label{dust-outside-Lambda}
\end{eqnarray}
\end{itemize}
\end{proposition}

\proof
(a) is a restatement of
Corollary~\ref{cor.reverse.clouds_bis}.
(b) then follows by using the equality in \reff{eq.mu-pilambda.eqineq}
together with Proposition~\ref{prop.fundamental},
exploiting the nonnegativity of all the operators in question.
\qed

\begin{proposition}[Comparison of cleaning sequences]
  \label{prop.comparison.sequences}
Let $(\mu_n)_{n \ge 1}$ and $(\nu_n)_{n \ge 1}$ 
be sequences of nonnegative clouds,
with $\mu_n$ $\Lambda$-regular and $\mu_n \audessus \nu_n$ for all $n$.
If $\nu_n$ converges pointwise to $\pi_\Lambda$
[i.e.\ $\lim\limits_{n\to\infty} \nu_n(\eta) = \pi_\Lambda(\eta)$
 for all markers $\eta$],
then $\mu_n$ also converges pointwise to $\pi_\Lambda$.
\end{proposition}

\proof
Since $\mu_n \audessus \nu_n$,
we have $\widetilde{\mu_n} \le \widetilde{\nu_n}$.
The pointwise convergence of $\nu_n$ to $\pi_\Lambda$
is equivalent to the pointwise convergence of
$\widetilde{\nu_n}$ to $\widetilde{\pi_\Lambda}$, where
\be
   \widetilde{\pi_\Lambda}(\eta)
   \;=\;
   \cases{  0   & if $\eta$ has all its entries in $\Lambda$  \cr
            1   & otherwise \cr
         }
\ee
(that is, $\widetilde{\pi_\Lambda} = {\bf 1} - {\bf 1}_\Lambda$).
Hence $\widetilde{\mu_n}(\eta) \to 0$
if $\eta$ has all its entries in $\Lambda$;
and for all other markers $\eta$ we have $\widetilde{\mu_n}(\eta) = 1$
for all $n$ by $\Lambda$-regularity.
Hence $\widetilde{\mu_n} \to \widetilde{\pi_\Lambda}$,
i.e.\ $\mu_n \to \pi_\Lambda$.
\qed

\subsection{Convergence of cleaning operators}
\label{sec4.8}

We are now ready to study the convergence of cleaning operators
$\beta_{h_1} \cdots \beta_{h_n}$ to $\Pi_\Lambda$,
analogously to what was done in Section~\ref{sec.operator.convergence}.
But in the cloud context we can shed additional light on this convergence
by dividing our analysis into two parts:
\begin{itemize}
   \item[1)]  Let $(\mu_n)_{n \ge 1}$ be a sequence of clouds that converges,
in some suitable topology, to a limiting cloud $\mu_\infty$.
Under what conditions can we show that the operators $T_{\mu_n}$
converge, in some correspondingly suitable topology, to $T_{\mu_\infty}$?
   \item[2)]  Consider the special case
$\mu_n = \beta_{h_1} * \ldots * \beta_{h_n}$.
Under what conditions on the sequence $(h_n)$
does $\mu_n$ converge to $\pi_\Lambda$ in the topology needed in part (1)?
\end{itemize}
We shall carry out this two-part analysis in two versions:
\begin{itemize}
   \item[(a)]  Pointwise convergence $\mu_n \to \mu_\infty$
      entails vector-norm convergence $T_{\mu_n} \to T_{\mu_\infty}$.
   \item[(b)]  Uniform-on-levels convergence $\mu_n \to \mu_\infty$
      entails operator-norm convergence $T_{\mu_n} \to T_{\mu_\infty}$.
\end{itemize}
In each case we shall require,
as was done in Section~\ref{sec.operator.convergence},
that $(I_\Lambda \alpha I_\Lambda)^k \to 0$
in a suitable topology.

\subsubsection{Pointwise (vector-norm) version}

We begin with the pointwise (vector-norm) version of the convergence theorems.
So let $(\mu_n)_{n \ge 1}$ be a sequence of clouds
such that, for each marker $\eta$,
the sequence $\mu_n(\eta)$ converges to a limit $\mu_\infty(\eta)$.
In order to control the convergence of the $T_{\mu_n}$,
we shall assume that all the clouds $\mu_n$ are $\Lambda$-regular
(for some fixed set $\Lambda \subseteq X$).
It is immediate from Definition~\ref{def.lambdaregular}
that the limiting cloud $\mu_\infty$ is likewise $\Lambda$-regular.

\begin{theorem}[Convergence theorem, pointwise version]
  \label{thm.main.clouds.AA}
\FH
Let $(\mu_n)_{n\ge1}$ be a sequence of $\Lambda$-regular clouds 
satisfying, for
each marker $\eta$,
\be
\label{convergence-to-mu-infty}
\mu_n(\eta)\; \mathop{\longrightarrow}_{n\rightarrow\infty}\; 
\mu_\infty(\eta)  \;.
\ee
Then, for any vector $c \ge 0$
with $cI_\Lambda\in l^1(w)$ and satisfying
\begin{equation}
\label{to-zero}
\| c(I_\Lambda \alpha I_\Lambda)^k \|_w
  \; \mathop{\longrightarrow}_{k\rightarrow\infty}  \;   0  \;,
\end{equation}
we have
\be
\| c (T_{\mu_n} - T_{\mu_\infty}) \|_w
  \; \mathop{\longrightarrow}_{n\rightarrow\infty} \; 0 \;.
 \label{eq.thm.main.clouds.AA.conclusion}
\ee

\end{theorem}

The proof of Theorem~\ref{thm.main.clouds.AA} will be based
on the following lemma:


\begin{lemma}
\label{comparison-mu-and-nu}
Let $\mu$ and $\nu$ be $\Lambda$-regular clouds. Then
\be
|\mu-\nu|\; \audessus\; 2 |\mu-\nu| * I_\Lambda.
 \label{eq.comparison-mu-and-nu}
\ee
\end{lemma}

\proof
We shall prove that
\be
|\mu-\nu|*I_{\Lambda^c} \;\audessus\; |\mu-\nu|*I_\Lambda \;,
\ee
which is obviously equivalent to \reff{eq.comparison-mu-and-nu}.
Since $\mu$ is $\Lambda$-regular, $\mu*I_{\Lambda^c}$ charges only
markers in $\partial\Lambda$, and for $\sigma\in\partial\Lambda$ we have
\be
(\mu*I_{\Lambda^c})(\sigma)\ =\ \mu(\sigma)\ =\ 
1-\sum_{\sigma' \ancneq\sigma}\mu(\sigma')\ =\ 
1-\widetilde{(\mu*I_\Lambda)}(\sigma)  \;.
\ee
The same result also holds for $\nu$, so that any marker
$\sigma\in\partial\Lambda$ we have
\be
(|\mu-\nu|*I_{\Lambda^c})(\sigma)\ \le\ 
\widetilde{\Bigl(|\mu-\nu|*I_\Lambda\Bigr)}(\sigma)  \;.
\ee
Now, for any marker $\eta$, either $\eta$ has all its entries lying in 
$\Lambda$,
in which case
\be
\widetilde{\Bigl(|\mu-\nu|*I_{\Lambda^c}\Bigr)}(\eta)\ =\ 0\ \le\ 
\widetilde{\Bigl(|\mu-\nu|*I_{\Lambda}\Bigr)}(\eta)  \;,
\ee
or else there exists a unique $\sigma\anc\eta$ with $\sigma\in\partial\Lambda$,
and then we have
\be
\widetilde{\Bigl(|\mu-\nu|*I_{\Lambda^c}\Bigr)}(\eta) \ =\
(|\mu-\nu|*I_{\Lambda^c})(\sigma)\ \le\
\widetilde{\Bigl(|\mu-\nu|*I_{\Lambda}\Bigr)}(\sigma) \ \le\
 \widetilde{\Bigl(|\mu-\nu|*I_{\Lambda}\Bigr)}(\eta)  \;.
\ee
\qed

\proofof{Theorem~\ref{thm.main.clouds.AA}}
Fix any $\epsilon>0$, and choose $k$ so that
\begin{equation}
\label{choice-of-k}
   \| c(I_\Lambda \alpha I_\Lambda)^k \|_w   \;\le\;  \epsilon  \;.
\end{equation}
It is easy to see that for all $j\ge 1$ we have
\be
    (I_\Lambda  * \rho^1 *  I_\Lambda)^{j+1}
    \;\audessus\;
    (I_\Lambda  * \rho^1 *  I_\Lambda)^{j}
    \;\audessus\; \ldots \;\audessus\;  I_\Lambda   \;;
\ee
therefore, by Proposition~\ref{prop.fundamental}, the sequence
$\left(\| c(I_\Lambda \alpha I_\Lambda)^j I_\Lambda \|_w\right)_{j\geq0}$
is decreasing.
In particular, the hypothesis $c I_\Lambda \in l^1(w)$
ensures that all these quantities are finite.
Therefore we can find a finite subset
$\Lambda' \subseteq \Lambda$ such that for every $j$ in the interval
$0\le j\le k$ we have
\begin{equation}
  \label{choice-of-Lambda'}
  \| c(I_\Lambda \alpha I_\Lambda)^j I_{\Lambda\setminus\Lambda'} \|_w
  \;\le\;  \epsilon/(k+1)   \;.
\end{equation}

Using Lemma~\ref{comparison-mu-and-nu} and 
Proposition~\ref{prop.fundamental},
we get
\be
\| c (T_{\mu_n} - T_{\mu_\infty}) \|_w
  \;\le\;  \| c T_{|\mu_n-\mu_\infty|} \|_w
  \;\le\; 2 \| c T_{|\mu_n-\mu_\infty|} I_\Lambda \|_w  \;.
 \label{norm-inequality-AA}
\ee
Furthermore, by the $\Lambda$-regularity of $\mu_n$ and $\mu_\infty$, 
the cloud $|\mu_n-\mu_\infty|*I_\Lambda$ charges only markers
with all their entries lying in $\Lambda$.
Let us now divide the set of such markers into three classes as follows:
\begin{itemize}
  \item[\ ] \emph{Class 1:}  Markers of level $\le k$ with all their
  entries lying in $\Lambda'$.
  \item[\ ] \emph{Class 2:}  Markers of level $\le k$ with at least one
  entry lying in $\Lambda\setminus\Lambda'$. 
  \item[\ ] \emph{Class 3:}  Markers of level $>k$.
\end{itemize}
We can then decompose any cloud $\mu_**I_\Lambda$ ($\mu_*=\mu_n$ or
$\mu_\infty$) in the form
$\mu_*^1+\mu_*^2+\mu_*^3$  so that $\mu_*^i$ charges only markers of 
class~$i$.
The triangle inequality yields
\be
\|  c T_{|\mu_n - \mu_\infty|} I_\Lambda\|_w
\ \le\ \|  c T_{|\mu_n^1 - \mu_\infty^1|} \|_w
    +  \|cT_{\mu_n^2}\|_w + \|cT_{\mu_\infty^2}\|_w
    +  \|cT_{\mu_n^3}\|_w + \|cT_{\mu_\infty^3}\|_w   \;,
\ee
and we will now bound these contributions separately as follows:

{\em Class 1.}\/
Since class 1 is a finite set of markers, (\ref{convergence-to-mu-infty})
ensures that $\|  c T_{|\mu_n^1 - \mu_\infty^1|} \|_w$
can be bounded by $\epsilon$ if $n$ has been chosen large enough.

{\em Class 2.}\/
We compare the cloud $\mu_*^2$ (which stands for either $\mu_n^2$ or 
$\mu_\infty^2$) with the cloud
\be
\nu^2\ :=\ \sum_{j=0}^k (I_\Lambda  * \rho^1 *  I_\Lambda)^{j} *
I_{\Lambda\setminus\Lambda'}  \;.
\ee
For each marker $\eta$, either there exists some $\sigma\anc\eta$ with 
$\sigma\in\partial\Lambda'$ and level$(\sigma)\le k$, in which case
$\widetilde{\nu^2} (\eta) = 1 \ge \widetilde{\mu_*^2} (\eta)$,
or else $\widetilde{\mu_*^2}(\eta)=0$
(remember that $\mu_*^2$ charges only markers of class~2).
Hence, we have
\be
\mu_*^2\ \audessus\ \nu^2  \;.
\ee
Then, Proposition~\ref{prop.fundamental} and~\reff{choice-of-Lambda'} 
together give
\be
\| cT_{\mu_*^2} \|_w\ \le\ \sum_{j=0}^k \| c(I_\Lambda \alpha I_\Lambda)^j
I_{\Lambda\setminus\Lambda'} \|_w  \;\le\;  \epsilon   \;.
\ee

{\em Class 3.}\/
The cloud $\mu_*^3$ (standing for $\mu_n^3$ or $\mu_\infty^3$) satisfies
\be
    \mu_*^3  \;\audessus\; (I_\Lambda  * \rho^1 *  I_\Lambda)^{k}  \;.
\ee
Applying again Proposition~\ref{prop.fundamental} gives,
by~\reff{choice-of-k},
\be
\| cT_{\mu_*^3} \|_w\ \le\  \| c(I_\Lambda \alpha I_\Lambda)^k \|_w
 \;\le\; \epsilon  \;.
\ee
\qed

\medskip

\begin{example}
\rm
It is natural to ask whether the hypothesis of
$\Lambda$-regularity in Theorem~\ref{thm.main.clouds.AA}
can be weakened to assuming that the $\mu_n$ are carried by $\Lambda$,
with $\triplenorm \mu_n \triplenorm = 1$ for all $n$
(or even $\triplenorm \mu_n \triplenorm$ uniformly bounded).
The answer is no, at least when $\Lambda^c$ is infinite
and $\Lambda \neq 0$:  it suffices to choose $x_0 \in \Lambda$
and a sequence of distinct elements $y_1,y_2,\ldots \in \Lambda^c$,
and to let $\mu_n$ be the cloud taking the value 1 on the
marker $(x_0 y_n)$ and 0 on all other markers.
Then $\mu_n \to 0$ pointwise;
but if we take $c = \delta_{x_0}$, $w = {\bf 1}$ and $\alpha = \smhalf I$,
we have $\| c T_{\mu_n} \| _w = \smhalf$ for all $n$.
So the role played by Lemma~\ref{comparison-mu-and-nu}
in handling markers ending in $\Lambda^c$ is apparently crucial.
\qed
\end{example}

\medskip

Here is one natural context in which a sequence of clouds
$(\mu_n)$ has a pointwise limit $\mu_\infty$:
Suppose that $(\mu_n)$ is a sequence of nonnegative clouds
that is decreasing in the sense of $\audessus$,
i.e.
\be
 \label{eq.decreasing_clouds}
\ldots \,\audessus\, \mu_{n+1} \,\audessus\, \mu_n \,\audessus\, \ldots
   \,\audessus\, \mu_1  \;.
\ee
Then for each marker $\eta$ the sequence $\widetilde{\mu_n}(\eta)$
is decreasing and bounded below by 0, hence has a limit.
Since $\mu_n(\eta) = \widetilde{\mu_n}(\eta) - \widetilde{\mu_n}(\eta^-)$,
we deduce that for each $\eta$ the sequence $\mu_n(\eta)$
converges as well, to a limit which we call $\mu_\infty(\eta)$.

In particular, we get a sequence of $\Lambda$-regular clouds
satisfying~\reff{eq.decreasing_clouds} by considering
\be
\label{product-of-beta_h}
\mu_n\ :=\ \beta_{h_1} * \ldots * \beta_{h_n} \;,
\ee
where $(h_n)_{n\ge1}$ is a sequence of functions taking values in
$[0,1]$ and supported on $\Lambda$. We are mainly interested in the
case where the limit $\mu_\infty$ of the sequence~(\ref{product-of-beta_h})
is equal to $\pi_\Lambda$.
The following proposition gives a necessary and sufficient condition
on the functions $(h_n)$ for this to happen:

\begin{proposition}[Conditions on cleaning sequences, pointwise version]
  \label{prop.infinite-sum}
 Let $(h_n)_{n\ge1}$ be a sequence of functions satisfying
$0\le h_n\le \ind{\Lambda}$, and define
\be
\mu_n\ :=\ \beta_{h_1} * \ldots * \beta_{h_n}  \;.
\ee
Then $\mu_n$ converges pointwise to $\pi_\Lambda$
if and only if
\begin{equation}
\label{infinite-sum}
\sum_{n=1}^\infty h_n  \;=\;  \infty \quad \mbox{on }\Lambda \;.
\end{equation}
\end{proposition}

\proof
Let us prove first that the condition~(\ref{infinite-sum}) is necessary.
Suppose that there exists $x\in\Lambda$ such that
\be
\sum_{n=1}^\infty h_n(x)  \;<\;  \infty \;,
\ee
and let $N$ be an integer such that
\be
\sum_{n>N} h_n(x)  \;<\; 1 \;.
\ee
We first claim that we can find a marker $\eta$ with last$(\eta)=x$
such that $\mu_N(\eta)>0$.
Indeed, since $\mu_N$ is supported on finitely many levels,
let us choose a marker $\sigma$ of maximum level with $\mu_N(\sigma)>0$.
Then every marker $\eta$ obtained from $\sigma$ by changing its last entry
--- in particular the one with last$(\eta)=x$ ---
satisfies $\mu_N(\eta)=\mu_N(\sigma)>0$
(this is a consequence of the fact that $\mu_N\in\scrs_1$).

Now, for each $n>N$, we have
\begin{subeqnarray}
\mu_n(\eta) & = & (\mu_N * \beta_{h_{N+1}} * \ldots * \beta_{h_n})(\eta) \\
   & \ge & \mu_N(\eta)\Bigl(1-h_{N+1}(x)\Bigr)\cdots \Bigl(1-h_n(x)\Bigr) \\
   & \ge & \mu_N(\eta)\left(  1-\sum_{k=N+1}^\infty h_k(x) \right) .
\end{subeqnarray}
It follows that
\be
\mu_n(\eta)\ \mathop{\not\!\!\!\longrightarrow}_{n\rightarrow\infty}\ 0\ 
=\ \pi_\Lambda(\eta)  \;.
\ee
This proves the necessity of (\ref{infinite-sum}).

Conversely, assume that~(\ref{infinite-sum}) is satisfied, and let us 
prove that
$\mu_\infty=\pi_\Lambda$. Since each $\mu_n$ is $\Lambda$-regular, it is
enough to verify that for each marker $\eta=(y_0,\ldots,y_k)$ with all
$y_j\in\Lambda$,
we have
\be
\mu_n(\eta)\ \mathop{\longrightarrow}_{n\rightarrow\infty}\ 0.
 \label{eq.mun_convergence}
\ee
We fix such a marker and set
\be
\Lambda' \;\bydef\; \{y_0,\ldots,y_k\}  \;.
 \label{eq.def_Lambdaprime}
\ee
Now, let $N$ be a large enough integer (to be made precise later).
Hypothesis~\reff{infinite-sum} ensures that if $n$ is large enough, we 
can find indices $0 = n_0 < n_1 < \ldots < n_N = n$
such that for each $j=1,\ldots,N$ and each $y\in\Lambda'$, we have
\be
1-\prod_{i=n_{j-1}+1}^{n_j} h_i(y)\ >\ \frac{1}{2}  \;.
\ee
Then Corollaries~\ref{cor.multi-monotonicity.clouds}
and \ref{cor.collapse.clouds} imply that
\be
\mu_n\ \audessus\ \beta_{1-\prod\limits_{i=n_0+1}^{n_1}(1-h_i)} * 
\cdots *  \beta_{1-\prod\limits_{i=n_{N-1}+1}^{n_N}(1-h_i)}\ 
\audessus\ \left(\beta_h\right)^{*N}  \;,
\ee
where
\be
h \;\bydef\; \frac{1}{2}\ind{\Lambda'}  \;.
\ee
Hence
\be
\widetilde\mu_n(\eta)\ \le\ \widetilde{\left(\beta_h\right)^{*N}} (\eta)  \;.
\ee
But an easy induction shows that for each $j=0,\ldots,k$, we have
\be
\left(\beta_h\right)^{*N} (y_0,\ldots,y_j)\ =\ 2^{-N}\ 
\binom{N}{j}\mathop{\longrightarrow}_{N\rightarrow\infty}\ 0  \;.
\ee
Thus, $\mu_n(\eta)$ can be made arbitrarily small if we take $N$ large 
enough.
\qed

Combining Theorem~\ref{thm.main.clouds.AA} and
Proposition~\ref{prop.infinite-sum},
we obtain Theorem~\ref{thm.main.operators} as an immediate corollary.


We can also obtain, as an easy corollary of
Proposition~\ref{prop.infinite-sum},
the following generalization of it:

\begin{corollary}[Cleaning by clouds]
   \label{cor.infinite-sum}
Let $(\nu_n)_{n \ge 1}$ be a sequence of $\Lambda$-regular clouds
satisfying
\be
   \sum_{n=1}^\infty \left[ 1 - \nu_n((x)) \right]
   \;=\;
   \infty
   \quad \mbox{for all } x \in\Lambda \;.
 \label{eq.cor.infinite-sum}
\ee
Then $\nu_1 * \ldots * \nu_n$ converges pointwise to $\pi_\Lambda$.

Conversely, if $\nu_1 * \ldots * \nu_n$ converges pointwise to $\pi_\Lambda$
and the $(\nu_n)_{n \ge 1}$ are all supported on finitely many levels
and belong to $\scrs_1$ (this last is automatic if $\Lambda \neq X$),
then \reff{eq.cor.infinite-sum} holds.
\end{corollary}

\proof
Define $h_n(x) = 1 - \nu_n((x))$.
By Lemma~\ref{lemma.mubetah} we have $\nu_i \audessus \beta_{h_i}$.
By Corollary~\ref{cor.powers} we have
$\nu_1 * \ldots * \nu_n \audessus \beta_{h_1} * \ldots * \beta_{h_n}$.
The result then follows from Propositions~\ref{prop.infinite-sum}
and \ref{prop.comparison.sequences}.

The converse holds under the specified conditions,
by exactly the same proof as in Proposition~\ref{prop.infinite-sum}.
\qed

{\bf Remark.}  For the converse, the hypothesis that the
$(\nu_n)_{n \ge 1}$ be supported on finitely many levels
(and belong to $\scrs_1$ if $\Lambda = X$) is essential:
a counterexample is $\nu_1 = \pi_\Lambda$, $\nu_2 = \nu_3 = \cdots = I$.
\qed

We can also weaken the condition on the sequence $(\nu_n)$
if we look only at the dust that stays inside $\Lambda$:

\begin{corollary}[Generalized cleaning by clouds, behavior inside
                         \mbox{\protect\boldmath $\Lambda$}]
   \label{cor.infinite-sum.generalized}
\quad\par\noindent
Let $(\nu_n)_{n \ge 1}$ be a sequence of nonnegative clouds
that are carried by $\Lambda$ and satisfy
$\triplenorm \nu_n \triplenorm \le 1$ and
\be
   \sum_{n=1}^\infty \left[ 1 - \nu_n((x)) \right]
   \;=\;
   \infty
   \quad \mbox{for all } x \in\Lambda \;.
\ee
Then $\nu_1 * \ldots * \nu_n * I_\Lambda$ converges pointwise to zero.
\end{corollary}
                                                                                
\proof
Note first that, by Lemma~\ref{lemma_betagi.clouds},
the cloud $\nu_1 * \ldots * \nu_n * I_\Lambda$ is supported on markers
having all their entries in $\Lambda$.
Now define $h_n(x) = 1 - \nu_n((x))$;
again by Lemma~\ref{lemma.mubetah} we have $\nu_i \audessus \beta_{h_i}$,
and by Corollary~\ref{cor.powers} we have
$\nu_1 * \ldots * \nu_n \audessus \beta_{h_1} * \ldots * \beta_{h_n}$.
By Proposition~\ref{prop.infinite-sum},
we have $\beta_{h_1} * \ldots * \beta_{h_n} \to \pi_\Lambda$ pointwise,
so that
\be
 \widetilde{(\beta_{h_1} * \ldots * \beta_{h_n})}(\eta) \;\to\; 0
\ee
for all markers $\eta$ having all their entries in $\Lambda$.
Since $\nu_1 * \ldots * \nu_n \audessus \beta_{h_1} * \ldots * \beta_{h_n}$,
it follows that, for such markers,
\be
   \widetilde{(\nu_1 * \ldots * \nu_n)}(\eta)  \;\to\;  0
\ee
and hence also $(\nu_1 * \ldots * \nu_n)(\eta) \to 0$.
\qed

{\bf Remark.}  Here there is no converse:
one could take, for instance,
$\nu_1 = I_\Lambda$ and $\nu_2 = \nu_3 = \cdots = I$.
\qed

\subsubsection{Uniform (operator-norm) version}

Let us now consider the uniform (operator-norm) version of the
convergence theorems.  We shall assume that $\mu_n$ converges to $\mu_\infty$
uniformly on all markers of a given level, i.e.
\be
   \sup\limits_{{\rm level}(\eta)=\ell} |\mu_n(\eta) - \mu_\infty(\eta)|
   \; \mathop{\longrightarrow}_{n\rightarrow\infty} \;
   0
\ee
for each $\ell \ge 0$.
In this case it is not necessary to assume that the $\mu_n$
are $\Lambda$-regular;  it suffices to assume that they are carried
by $\Lambda$ and uniformly bounded in norm.

\begin{theorem}[Convergence theorem, uniform version]
  \label{thm.main.clouds.BB}
\FH
Let $(\mu_n)_{n\ge1}$ be a sequence of clouds carried by $\Lambda$
and satisfying a uniform bound $\triplenorm \mu_n \triplenorm \le M < \infty$.
Suppose further that
\be
\label{uniform-convergence-to-mu-infty}
\mu_n(\eta)\; \mathop{\longrightarrow}_{n\rightarrow\infty}\;
\mu_\infty(\eta)  \quad\mbox{uniformly for all $\eta$ of a given level}
\ee
and that
\be
\label{operator-convergence-to-zero}
\|(I_\Lambda \alpha I_\Lambda)^k \|_{w\to w}\; 
\mathop{\longrightarrow}_{k\rightarrow\infty}\; 0.
\ee
Then
\be
\label{operator-convergence-to-mu-infty}
\| T_{\mu_n} - T_{\mu_\infty} \|_{w\to w}\; 
\mathop{\longrightarrow}_{n\rightarrow\infty}\; 0.
\ee
\end{theorem}

\proof
We proceed as in the proof of {Theorem~\ref{thm.main.clouds.AA}}, 
but here the reasoning is much simpler since
the uniform-convergence hypothesis \reff{uniform-convergence-to-mu-infty}
allows us to avoid the introduction of the subset $\Lambda'$.
We begin by observing that $\mu_\infty$ is carried by $\Lambda$
and satisfies
\be
   \triplenorm \mu_\infty \triplenorm   \;\le\;
   \liminf_{n \to \infty}  \triplenorm \mu_n \triplenorm  \;\le\;  M  \;.
\ee
Now fix any $\epsilon>0$, and choose $k$ so that
\begin{equation}
\label{choice-of-k-BB}
   \| (I_\Lambda \alpha I_\Lambda)^k \|_{w\to w}   \;\le\;  \epsilon/M  \;.
\end{equation}
We divide the set of all markers into two classes:
\begin{itemize}
  \item[\ ] \emph{Class 1:}  Markers of level $\le k$.
  \item[\ ] \emph{Class 2:}  Markers of level $>k$.
\end{itemize}
We decompose any cloud $\mu_*$ ($\mu_*=\mu_n$ or $\mu_\infty$) in the form
$\mu_*^1+\mu_*^2$  so that $\mu_*^i$ charges only markers of class~$i$.
The triangle inequality yields
\be
   \| T_{\mu_n - \mu_\infty} \|_{w \to w}
   \;\le\;
   \| T_{|\mu_n^1 - \mu_\infty^1|} \|_{w \to w}
   \,+\, \| T_{|\mu_n^2|} \|_{w \to w}
   \,+\, \| T_{|\mu_\infty^2|} \|_{w \to w}
   \;,
\ee
and these contributions may be bounded separately as follows:

{\em Class 1.}\/
The cloud $|\mu_n^1-\mu_\infty^1|$ charges only markers of level $\le k$.
By (\ref{uniform-convergence-to-mu-infty}),
if $n$ is large enough we have
\be
|\mu_n^1-\mu_\infty^1|\; \le\; 
\frac{\epsilon}{k+1}(\rho^0+\rho^1+\cdots+\rho^k)  \;,
\ee
so that $\triplenorm \, |\mu_n^1-\mu_\infty^1| \, \triplenorm \le \epsilon$.
Therefore, by Proposition~\ref{prop.Tnu.norm},
$\| T_{|\mu_n^1 - \mu_\infty^1|} \|_{w \to w} \le \epsilon$ as well.

{\em Class 2.}\/
Since $\mu_*$ (= $\mu_n$ or $\mu_\infty$) is carried by $\Lambda$,
any marker $\eta=(x_0,\ldots,x_\ell)$ charged by $\mu_*^2$
must have $\ell > k$ and $x_0,\ldots,x_k \in \Lambda$.
Since $\triplenorm \mu_* \triplenorm \le M$,
it follows that 
\be
   |\mu_*^2| \;\audessus\;  M (I_\Lambda * \rho^1 * I_\Lambda)^k  \;.
\ee
By Proposition~\ref{prop.fundamental} we have
$\| T_{|\mu_*^2|} \|_{w \to w} \le
 M \| (I_\Lambda \alpha I_\Lambda)^k \|_{w \to w} \le \epsilon$
by~\reff{choice-of-k-BB}.
\qed

\medskip

In the uniform case we have the following analogue of
Proposition~\ref{prop.infinite-sum}:

\begin{proposition}[Conditions on cleaning sequences, uniform version]
  \label{prop.infinite-sum.BB}
\hfill\break
Let $(h_n)_{n\ge1}$ be a sequence of functions satisfying
$0\le h_n\le \ind{\Lambda}$, and define
\be
\mu_n\ :=\ \beta_{h_1} * \ldots * \beta_{h_n}  \;.
\ee
Then the following are equivalent:
\begin{itemize}
   \item[(a)] $\mu_n$ converges uniformly-on-levels to $\pi_\Lambda$.
   \item[(b)] $\sum_{i=1}^n h_i$ converges uniformly (as $n\to\infty$)
       to $+\infty$ on $\Lambda$.
   \item[(b\/\textprime)] There exist $\delta>0$ and a sequence
        $0 = n_0 < n_1 < n_2 < \ldots$ satisfying
\be\label{roma.r23.clouds}
\sum_{i=n_{j-1}+1} ^{n_j} h_i \;\ge\; \delta\,\chi_\Lambda
\ee
        for all $j$.
   \item[(b\/\textprime\textprime)] For all $\delta < \infty$,
        there exists a sequence
        $0 = n_0 < n_1 < n_2 < \ldots$ satisfying \reff{roma.r23.clouds}
        for all $j$.
\end{itemize}
\end{proposition}

\proof
Since the $h_i$ are bounded,
it is easy to see that (b), (b\textprime) and (b\textprime\textprime)
are all equivalent.

The proof of (b\/\textprime\textprime) $\implies$ (a)
is a straightforward adaptation of the corresponding proof
in Proposition~\ref{prop.infinite-sum}.
Indeed, we want to show here that the convergence \reff{eq.mun_convergence}
holds uniformly for all markers of a given level.
With the hypothesis (b\/\textprime\textprime),
we can repeat the end of the proof of Proposition~\ref{prop.infinite-sum}
with $\Lambda'=\Lambda$ instead of \reff{eq.def_Lambdaprime},
which gives the result.

Conversely, let us prove that (a) $\implies$ (b).
Suppose that $\sum_{i=1}^n h_i$ does not converge uniformly to $+\infty$
on $\Lambda$.
Then there exists $M < \infty$ such that,
for any $n\ge1$, we can find $x\in\Lambda$ with
\be
\label{non-uniform-convergence}
\sum_{i=1}^n h_i(x)\; <\; M  \;.
\ee
For such $n$ and $x$, we consider the markers $\eta$
whose entries are all equal to $x$,
and we denote by $\eta_k$ the only such marker that lies on level $k$.
For $0\le j\le n$ and $k\ge 0$, let us define
\be
p_j(k)\;\bydef\; \mu_j(\eta_k)  \;,
\ee
where $\mu_0 \bydef \rho^0$.
We have $p_0(0)=1$, $p_0(k)=0$ for each $k>0$, and for $j\ge 1$
\be
p_{j}(k)\;=\; [1-h_j(x)] \, p_{j-1}(k) \,+\, h_j(x) \, p_{j-1}(k-1)  \;.
\ee
(The preceding equation holds also for $k=0$ if for each $j$
 we set $p_j(-1) \bydef 0$.)
Observe that for each $j$ we always have $\sum_k p_j(k)=1$.
Let us introduce the barycenter
\be
G_j  \;\bydef\; \sum_k k \, p_j(k)  \;.
\ee
We obviously have $G_0=0$, and a direct calculation gives, for $j\ge 1$,
\be
G_{j}\;=\;G_{j-1}+h_j(x)  \;.
\ee
Therefore, we have
\be
G_n\;=\;\sum_{j=1}^nh_j(x)\;<\;M  \;.
\ee
By the Markov inequality, this implies that
\be
\sum_{k=0}^{2M}\mu_n(\eta_k)\;=\;
\sum_{k=0}^{2M}p_n(k)\;\ge\; 1/2 \;,
\ee
which clearly prevents $\mu_n$ from converging uniformly to $\pi_\Lambda$.
\qed

Combining Theorem~\ref{thm.main.clouds.BB} and
Proposition~\ref{prop.infinite-sum.BB},
we obtain Theorem~\ref{thm.main.operators.uniform} as an immediate corollary.

\subsubsection{Some final remarks}

There are three natural topologies for convergence of clouds:
$\triplenorm \,\cdot\, \triplenorm$ norm;
uniform-on-levels;  and pointwise.
Norm convergence of the $\mu_n$ implies operator-norm convergence
of the $T_{\mu_n}$ without any hypothesis on $\alpha$
beyond the Fundamental Hypothesis
(cf.\ Proposition~\ref{prop.Tnu.norm}).
Uniform-on-levels convergence is weaker than norm convergence
in that the mass of $\mu_n - \mu_\infty$ can run ``upwards to infinity'',
as in the example $\beta_\Lambda^n \to \pi_\Lambda$.
This is handled in Theorem~\ref{thm.main.clouds.BB} by assuming that
$\| (I_\Lambda \alpha I_\Lambda)^k \| _{w \to w} \to 0$;
then the $T_{\mu_n}$ converge in operator norm.
Finally, pointwise convergence is still weaker than
uniform-on-levels convergence because mass can also run
``outwards to infinity'' (when $X$ is infinite),
as in the example $\beta_{\Lambda_n} \to \beta_\Lambda$
with $\Lambda_n \uparrow \Lambda$ and all $\Lambda_n \subsetneq \Lambda$.
This is handled in Theorem~\ref{thm.main.clouds.AA}
by working on a fixed vector $c \in l^1(w)$
and demanding convergence only in vector norm (not operator norm);
the upwards-running mass is handled by assuming that
$\| c(I_\Lambda \alpha I_\Lambda)^k \|_w \to 0$.

\subsection{Some further identities and inequalities}  \label{sec4.9}

Let us now present some further identities and inequalities
that will play an important role in Section~\ref{sec.converse}.
Our first result is a cloud analogue and extension of
Lemmas~\ref{lemma.12Apr07.identity} and \ref{lemma.12Apr07.inequality}:

\begin{lemma}
   \label{lemma.12Apr07.cloud}

(a) Let $h \colon\, X \to [0,1]$ and set $\Lambda \bydef \supp h$.
Then we have the identity
\be
   \sum_{k=0}^\infty \beta_h^{*k} * I_h  \;=\;  {\bf 1}_\Lambda  \;.
 \label{eq.lemma.12Apr07.cloud.a}
\ee
[Recall that ${\bf 1}_\Lambda$ is the cloud that takes the value 1
 on markers having all their entries in $\Lambda$, and 0 elsewhere.]

(b) Let $\nu \ge 0$ be a cloud satisfying $\triplenorm \nu \triplenorm \le 1$.
Set $h(x) \bydef 1 - \nu((x))$.
Then we have the inequality
\be
   \sum_{k=0}^\infty \nu^{*k} * I_h  \;\le\;  {\bf 1}   \;.
 \label{eq.lemma.12Apr07.cloud.b}
\ee
\end{lemma}

\proof
(a)  The proof of Lemma~\ref{lemma.12Apr07.identity}
is in fact a proof of the cloud identity \reff{eq.lemma.12Apr07.cloud.a}:
no more need be said.
But just for completeness, here is an alternate proof:
Let $\mu$ be any cloud; then by \reff{somme_ancetres3a} we have
\be
   \mu * I_h  \;=\;  \widetilde{\mu} \,-\, \widetilde{\mu * \beta_h}
   \;.
\ee
Now set $\mu = \beta_h^{*k}$ and sum from $k=0$ to $N$:  we get
\be
   \sum_{k=0}^N \beta_h^{*k} * I_h
   \;=\;
   {\bf 1} \,-\, \widetilde{\beta_h^{*(N+1)}}
\ee
(since $\widetilde{\beta_h^{*0}} = \widetilde{\rho^0} = {\bf 1}$).
Now take $N \to\infty$:
by Proposition~\ref{prop.infinite-sum},
$\beta_h^{*(N+1)}$ converges pointwise to $\pi_\Lambda$,
hence $\widetilde{\beta_h^{*(N+1)}}$ converges pointwise to
$\widetilde{\pi_\Lambda} = {\bf 1} - {\bf 1}_\Lambda$.
This proves \reff{eq.lemma.12Apr07.cloud.a}.

(b)
By Lemma~\ref{lemma.10Apr07.a} we have, for any cloud $\mu \ge 0$,
\begin{subeqnarray}
   \widetilde{(\mu * \nu)}(\eta)
    & \le &
   \widetilde{\mu}(\eta^-)
    \:+\:
   \mu(\eta) \, \nu\Bigl((\hbox{last}(\eta) )\Bigr)
        \\[2mm]
    & = &
   \widetilde{\mu}(\eta)
    \:-\:
    \left[ 1 \,-\, \nu\Bigl((\hbox{last}(\eta) )\Bigr) \right]
    \mu(\eta)
\end{subeqnarray}
and hence
\be
   \mu * I_h  \;\le\;  \widetilde{\mu}  \,-\,  \widetilde{\mu * \nu}
   \;.
\ee
Now take $\mu = \nu^{*k}$ and sum from $k=0$ to $N$:
\be
   \sum_{k=0}^N \nu^{*k} * I_h
   \;\le\;
   {\bf 1}  \,-\,  \widetilde{\nu^{*(N+1)}}
   \;\le\;
   {\bf 1}
\ee
(since $\widetilde{\nu^{*0}} = \widetilde{\rho^0} = {\bf 1}$).
Taking $N \to\infty$ gives the result.
\qed

We also have the following curious collection of inequalities and identity
involving the cloud $\mu \bydef {\bf 1} - \widetilde{\nu}$.
These too will play a central role in Section~\ref{sec.converse};
what makes them so powerful is that the inequalities go in
{\em both}\/ directions.

\begin{lemma}
   \label{lemma.16Apr07}
Let $\nu \ge 0$ be a cloud satisfying $\triplenorm \nu \triplenorm \le 1$.
Define $h(x) \bydef 1 - \nu((x))$ and $\mu \bydef {\bf 1} - \widetilde{\nu}$.
\begin{itemize}
   \item[(a)] We have $I_h \le \mu \le {\bf 1}$ and
\be
   \sum_{n=0}^\infty \nu^{*n} * \mu  \;\le\;  {\bf 1}   \;.
 \label{eq.lemma.16Apr07.a}
\ee
   \item[(b)] If $\nu$ is carried by $\Lambda$
       and satisfies $\nu((x)) < 1$ for all $x \in \Lambda$
       [i.e.\ $\supp h = \Lambda$], then
\be
   \sum_{n=0}^\infty \nu^{*n} * \mu  \;\ge\;  {\bf 1}_\Lambda   \;.
 \label{eq.lemma.16Apr07.b}
\ee
   \item[(c)] If $\nu$ is $\Lambda$-regular
       and satisfies $\nu((x)) < 1$ for all $x \in \Lambda$, then
       $\mu$ is supported on markers having all their entries in $\Lambda$
       (so that $I_h \le \mu \le {\bf 1}_\Lambda$) and
\be
    \sum_{n=0}^{\infty} (I_\Lambda * \nu * I_\Lambda)^{*n} * \mu
    \;=\;
    \sum_{n=0}^{\infty} \nu^{*n} * \mu
    \;=\;
    {\bf 1}_{\Lambda}
    \;.
 \label{eq.lemma.16Apr07.c}
\ee
%
   \item[(d)] If $\nu$ belongs to $\scrs_1$ and is supported on levels $\le K$,
      then $\mu$ is supported on levels $\le K-1$.
\end{itemize}
\end{lemma}

Please note that since $\mu \ge I_h$,
Lemma~\ref{lemma.16Apr07}(a) is a strengthening of
Lemma~\ref{lemma.12Apr07.cloud}(b).
Furthermore, if $\nu = \beta_h$, then $\mu = I_h$,
so that Lemma~\ref{lemma.16Apr07}(c) is
a direct generalization of Lemma~\ref{lemma.12Apr07.cloud}(a).

\proofof{Lemma~\ref{lemma.16Apr07}}
(a) Since $\triplenorm \nu \triplenorm \le 1$, we have $\mu \ge 0$.
Since $\nu \ge 0$, we have $\mu \le 1$.
Finally, since $\mu((x)) = 1 - \nu((x)) = h(x)$, we have $\mu \ge I_h$.

Let us now observe that
$\mu \bydef {\bf 1} - \widetilde{\nu} = (\rho^0 - \nu) * {\bf 1}$
and perform the following calculation:
\begin{subeqnarray}
   \sum\limits_{n=0}^{N-1} \nu^{*n} * \mu
   & = &
   \sum\limits_{n=0}^{N-1} \nu^{*n} * (\rho^0 - \nu) * {\bf 1}
        \\[1mm]
   & = &
   (\rho^0 - \nu^{*N}) * {\bf 1}
        \\[1mm]
   & = &
   {\bf 1} \,-\, \widetilde{\nu^{*N}}
       \slabel{eq.proof.lemma.11Apr07}  \\[1mm]
   & \le &
   {\bf 1} \;.
       \slabel{eq.proof.lemma.11Apr07.last}
       \label{eq.proof.lemma.11Apr07.all}
\end{subeqnarray}
Now take $N \to \infty$:
since $\mu,\nu \ge 0$, the left-hand side increases pointwise
to a limiting cloud $\sum\limits_{n=0}^\infty \nu^{*n} * \mu$
whose values lie in $[0,1]$ by \reff{eq.proof.lemma.11Apr07.last}.

(b) Now suppose that $\nu$ is carried by $\Lambda$
and that $\nu((x)) < 1$ for all $x \in \Lambda$.
Then Corollary~\ref{cor.infinite-sum.generalized}
implies that $\nu^{*N}(\sigma) \to 0$ for all markers $\sigma$
ending in $\Lambda$, hence
$\widetilde{\nu^{*N}}(\eta) \to 0$ for all markers $\eta$
having all their entries in $\Lambda$.
By (\ref{eq.proof.lemma.11Apr07.all}a--c)
this implies \reff{eq.lemma.16Apr07.b}.

(c) Since $\nu$ is $\Lambda$-regular,
we have $\widetilde{\nu}(\eta) = 1$
whenever $\eta$ has at least one entry outside $\Lambda$;
hence $\mu$ is supported on markers having all their entries in $\Lambda$.
Since in addition $\nu((x)) < 1$ for all $x \in \Lambda$,
we can apply Corollary~\ref{cor.infinite-sum}
to conclude that $\nu^{*N} \to \pi_\Lambda$ pointwise,
i.e.\ $\widetilde{\nu^{*N}} \to \widetilde{\pi_\Lambda}$ pointwise.
But
\be
   \widetilde{\pi_\Lambda}(\eta)
   \;=\;
   \cases{  0   & if $\eta$ has all its entries in $\Lambda$  \cr
            1   & otherwise \cr
         }
\ee
so that $\widetilde{\pi_\Lambda} = {\bf 1} - {\bf 1}_\Lambda$.
By (\ref{eq.proof.lemma.11Apr07.all}a--c)
this proves the second equality in \reff{eq.lemma.16Apr07.c}.
Since $\mu = I_\Lambda * \mu * I_\Lambda$,
we have $\nu^{*n} * \mu = (I_\Lambda * \nu * I_\Lambda)^{*n} * \mu$
by Lemma~\ref{lemma_betagi.clouds},
which proves the first equality.

(d) If $\nu \in \scrs_1$ is supported on levels $\le K$,
it follows that $\widetilde{\nu}(\eta) = 1$
whenever $\hbox{level}(\eta) \ge K$,
hence $\mu$ is supported on levels $\le K-1$.
\qed

{\bf Remark.}  At a formal level, the computation
\reff{eq.proof.lemma.11Apr07.all} is inspired by the
(admittedly meaningless) ``identity''
$\sum\limits_{n=0}^\infty \nu^{*n} = (\rho^0 - \nu)^{-1}$.
Another, more physical, way of expressing the intuition behind
\reff{eq.proof.lemma.11Apr07.all} is to observe that
\be
   \nu^{*n} * \mu  \;=\;  \widetilde{\nu^{*n}} \,-\, \widetilde{\nu^{*(n+1)}}
   \;,
\ee
so that $(\nu^{*n} * \mu)(\eta)$ measures the mass that lies $\anc \eta$
after $n$ steps but gets pushed above (or out) at the $(n+1)$st step.
Summing over $n$, we should get 1 whenever $\eta$
has all its entries in $\Lambda$,
since by Corollary~\ref{cor.infinite-sum.generalized}
all the mass should eventually be pushed out.
\qed

%
%
%
%
%

\section{Alternate sufficient conditions for cleanability}
     \label{sec.complements}

{\em In this section we do {\bf not} assume the Fundamental Hypothesis.}\/
Rather, our goal is to examine briefly the conditions under which
cleaning can be assured even in the absence of the Fundamental Hypothesis.
Please note that we are entitled to use here those results
of Section~\ref{sec4} that refer only to clouds.
However, we must scrupulously avoid using those results which,
like Propositions~\ref{prop.Tnu.norm} and \ref{prop.fundamental},
refer to operators and therefore depend on the Fundamental Hypothesis.

First, we need a few definitions.
Let $\amat = (a_{xy})_{x,y \in X}$ be any nonnegative matrix
(in practice, we will take $\amat$ to be either
 $\alpha$ or $I_\Lambda \alpha I_\Lambda$).
For any marker $\eta = (x_0,\ldots,x_k)$, we define
\be
   \amat^\eta  \;\bydef\;  a_{x_0 x_1} a_{x_1 x_2} \cdots a_{x_{k-1}x_k}
   \;.
\ee
We then define the spaces of clouds
$l^1(\amat)$ and $l^\infty(\amat)$ by the norms
\begin{eqnarray}
   \| \nu \|_{1,\amat}        & \bydef &
        \sum\limits_\eta  \amat^\eta \, |\nu(\eta)|    \\[2mm]
   \| \nu \|_{\infty,\amat}   & \bydef &
        \sup\limits_{\eta \colon\: \amat^\eta > 0}  |\nu(\eta)|
\end{eqnarray}
If $\mu$ and $\nu$ are clouds, we write
\be
   \< \mu , \nu \>_\amat   \;\bydef\;
   \sum_\eta \amat^\eta \, \mu(\eta) \, \nu(\eta)
\ee
whenever this sum has an unambiguous meaning;
in particular, if $\mu,\nu \ge 0$,
then $\< \mu , \nu \>_\amat$ is well-defined though it may be $+\infty$.

If $\nu$ is a cloud, we would like to define an operator $T_{\nu,\amat}$
by the usual formula
\be
   T_{\nu,\amat}  \;\bydef\;
   \sum_{\eta=(x_0,\ldots,x_k)}  \nu(\eta) \,
             I_{\{x_0\}} \amat I_{\{x_1\}} \amat \cdots
                  \amat I_{\{x_{k-1}\}} \amat I_{\{x_k\}}
   \;.
 \label{def.complements.Tnu}
\ee
[When $\amat=\alpha$, we write simply $T_\nu$.]
The trouble is that, in the absence of the Fundamental Hypothesis
or its equivalent for $\amat$,
it is difficult to guarantee that this sum is well-defined
(compare Proposition~\ref{prop.Tnu.norm}).
We therefore restrict ourselves in this section
to {\em nonnegative}\/ clouds $\nu$
and consider \reff{def.complements.Tnu} as defining
a matrix $T_{\nu,\amat}$ whose elements lie in $[0,+\infty]$.
Of course, when necessary we shall verify {\em a posteriori}\/
that the elements of $T_{\nu,\amat}$ are finite.

Finally, if $c \ge 0$ and $w \ge 0$ are vectors,
we define a nonnegative cloud $\mu_{c,w}$ by
\be
   \mu_{c,w}\bigl( (x_0,\ldots,x_k) \bigr)  \;=\; c_{x_0} w_{x_k}
   \;.
\ee
It is easily checked that for any nonnegative cloud $\nu$ we have
\be
   c T_{\nu,\amat} w  \;=\;  \< \mu_{c,w} , \nu \>_\amat
   \;.
\ee
(We define $0 \cdot \infty = 0$;
 this may be needed to interpret the inner product $c T_{\nu,\amat} w$.)
In particular, taking $\nu = \rho^k$, we have
\be
  c\amat^k w
  \;=\;
  \sum_{\hbox{\scriptsize level}(\eta)=k}  \amat^\eta \, \mu_{c,w}(\eta)
 \label{eq.c_alpham_w}
\ee
and hence that
\be
   \sum_{k=0}^\infty  c \amat^k w   \;=\;
   \| \mu_{c,w} \|_{1,\amat}  \;.
 \label{eq.c_alphasum_w}
\ee

The following theorem and its corollary give a sufficient condition
for a dust vector $c$ to be cleanable in $l^1(w)$ sense
by every cleaning process that visits each site infinitely many times.
The proofs are almost trivial,
now that we have the results of Section~\ref{sec4} in hand.

\begin{theorem}
   \label{thm.complements.1}
Let $\amat$ be a nonnegative matrix,
let $\mu$ be a nonnegative cloud belonging to $l^1(\amat)$,
and let $(\nu_n)_{n \ge 1}$ be a sequence of nonnegative clouds
satisfying a uniform bound $\| \nu_n \|_{\infty,\amat} \le M < \infty$
and tending pointwise $\amat$-a.e.\ to zero
[i.e.\ $\lim_{n \to\infty} \nu_n(\eta) = 0$ for each marker $\eta$
 having $\amat^\eta > 0$].
Then
\be
   \lim_{n \to\infty}  \< \mu , \nu_n \>_\amat   \;=\;   0  \;.
\ee
\end{theorem}

\proof
An immediate consequence of Lebesgue's dominated convergence theorem.
\qed

\begin{corollary}[Sufficient condition for universal cleaning]
   \label{cor.complements.1}
Let $\Lambda \subseteq X$,
and let $c \ge 0$ and $w \ge 0$ be vectors satisfying
\be
   \sum_{k=0}^\infty  c (I_\Lambda \alpha I_\Lambda)^k w   \;<\;  \infty  \;.
 \label{hyp.cor.complements.1}
\ee
If $(f_n)_{n\ge1}$ are functions $X \to [0,1]$, supported in $\Lambda$
and satisfying $\sum_n f_n = +\infty$ everywhere on $\Lambda$, we have
\be
   \lim_{n \to\infty}  c \beta_{f_1} \cdots \beta_{f_n} I_\Lambda w
   \;=\;   0  \;.
 \label{eq.cor.complements.1b}
\ee
More generally, if $(\omega_n)_{n\ge1}$ are nonnegative clouds
carried by $\Lambda$ that satisfy
$\triplenorm \omega_n \triplenorm \le 1$ and
$\sum_n [1 - \omega_n((x))] = +\infty$ for all $x \in \Lambda$,
we have
\be
   \lim_{n \to\infty}  c T_{\omega_1} \cdots T_{\omega_n} I_\Lambda w
   \;=\;   0  \;.
 \label{eq.cor.complements.1c}
\ee
\end{corollary}

\proof
This is simply the special case of Theorem~\ref{thm.complements.1}
in which $\amat = I_\Lambda \alpha I_\Lambda$,
$\mu = \mu_{c,w}$ and
$\nu_n = \omega_1 * \cdots * \omega_n * I_\Lambda$:
by \reff{eq.c_alphasum_w}, the hypothesis \reff{hyp.cor.complements.1}
guarantees that $\mu_{c,w} \in l^1(\amat)$.
[Using Lemma~\ref{lemma_betagi.clouds} to rewrite $\nu_n$ as
 $I_\Lambda * \omega_1 * I_\Lambda * \cdots *
  I_\Lambda * \omega_n * I_\Lambda$,
 we see that $\< \mu_{c,w} , \nu_n \>_\alpha =
              \< \mu_{c,w} , \nu_n \>_{I_\Lambda \alpha I_\Lambda}$.]
We can take $M=1$.
The fact that $\nu_n(\eta) \to 0$ for all $\eta$ is
an immediate consequence of Corollary~\ref{cor.infinite-sum.generalized}.
%
%
\qed

%
%
%
%

In view of Theorem~\ref{thm.main.clouds.AA},
it is natural to wonder whether the hypothesis
$\sum_{k=0}^\infty c \alpha^k w < \infty$ in Corollary~\ref{cor.complements.1}
can be weakened to assuming that
$c \alpha^k w$ is finite for every $k$
and tends to zero as $k \to\infty$.
It turns out that this is not the case;
indeed, the following example shows that
Corollary~\ref{cor.complements.1} is in a certain sense sharp:

\begin{figure}
\begin{picture}(0,0)%
\includegraphics{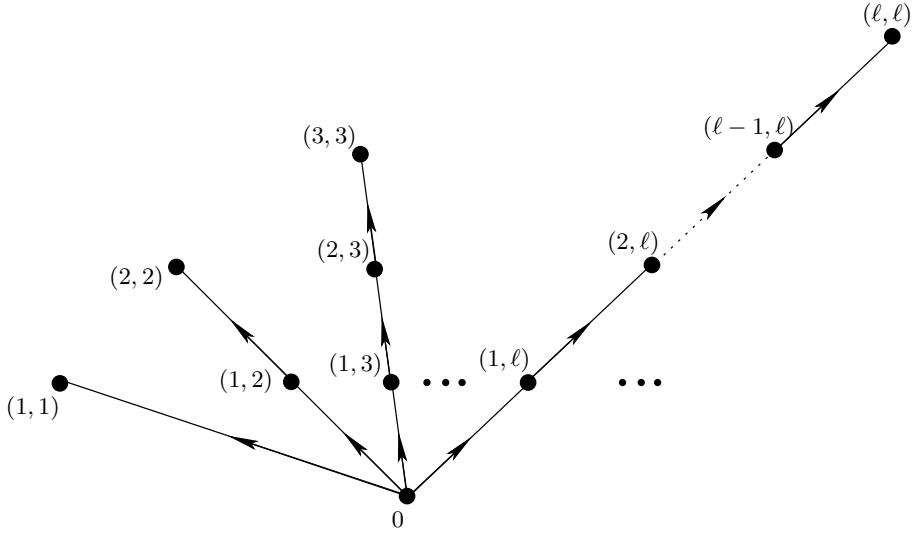}%
\end{picture}%
\setlength{\unitlength}{4144sp}%
\begingroup\makeatletter\ifx\SetFigFont\undefined%
\gdef\SetFigFont#1#2#3#4#5{%
   \reset@font\fontsize{#1}{#2pt}%
   \fontfamily{#3}\fontseries{#4}\fontshape{#5}%
   \selectfont}%
\fi\endgroup%
\begin{picture}(5358,3159)(1269,-8119)
\put(3578,-8088){\makebox(0,0)[lb]{\smash{{\SetFigFont{9}{10.8}{\rmdefault}{\mddefault}{\updefault}$0$}}}}
\put(6391,-5056){\makebox(0,0)[lb]{\smash{{\SetFigFont{9}{10.8}{\rmdefault}{\mddefault}{\updefault}$(\ell,\ell)$}}}}
\put(4094,-7124){\makebox(0,0)[lb]{\smash{{\SetFigFont{9}{10.8}{\rmdefault}{\mddefault}{\updefault}$(1,\ell)$}}}}
\put(4866,-6424){\makebox(0,0)[lb]{\smash{{\SetFigFont{9}{10.8}{\rmdefault}{\mddefault}{\updefault}$(2,\ell)$}}}}
\put(5438,-5736){\makebox(0,0)[lb]{\smash{{\SetFigFont{9}{10.8}{\rmdefault}{\mddefault}{\updefault}$(\ell-1,\ell)$}}}}
\put(3128,-6473){\makebox(0,0)[lb]{\smash{{\SetFigFont{9}{10.8}{\rmdefault}{\mddefault}{\updefault}$(2,3)$}}}}
\put(3043,-5792){\makebox(0,0)[lb]{\smash{{\SetFigFont{9}{10.8}{\rmdefault}{\mddefault}{\updefault}$(3,3)$}}}}
\put(2543,-7248){\makebox(0,0)[lb]{\smash{{\SetFigFont{9}{10.8}{\rmdefault}{\mddefault}{\updefault}$(1,2)$}}}}
\put(1889,-6632){\makebox(0,0)[lb]{\smash{{\SetFigFont{9}{10.8}{\rmdefault}{\mddefault}{\updefault}$(2,2)$}}}}
\put(1269,-7409){\makebox(0,0)[lb]{\smash{{\SetFigFont{9}{10.8}{\rmdefault}{\mddefault}{\updefault}$(1,1)$}}}}
\put(3198,-7153){\makebox(0,0)[lb]{\smash{{\SetFigFont{9}{10.8}{\rmdefault}{\mddefault}{\updefault}$(1,3)$}}}}
\end{picture}%
\caption{
   The directed graph associated to the matrix $\alpha$
   in Example~\ref{example_tree}.
   The vertices of this graph correspond to the points of $X$,
   and the edges correspond to the nonzero matrix elements of $\alpha$.
}
\label{fig_tree}
\end{figure}

\begin{example}
    \label{example_tree}
\rm
Let $(\rho_k)_{k \ge 0}$ be {\em any}\/ sequence of positive numbers
satisfying $\sum_{k=0}^\infty \rho_k = \infty$.
Then, on a countably infinite state space $X$,
we can construct a matrix $\alpha \ge 0$ and vectors $c,w \ge 0$ such that
\begin{itemize}
   \item[(a)]  $c \alpha^k w = \rho_k$ for every $k \ge 0$.
   \item[(b)]  There exists a sequence $x_1,x_2,\ldots \in X$
       in which each element of $X$ occurs infinitely often,
       and for which we have
       $\lim\limits_{n \to \infty} c \beta_{x_1}\cdots\beta_{x_n} w = \infty$.
\end{itemize}
To see this, take $X = \{0\}\cup\{(k,\ell) \colon\, 1\le k\le \ell\}$,
and define the matrix $\alpha$ by
\begin{subeqnarray}
   \alpha_{0,(1,\ell)}           & = &  \rho_1 / (2^\ell \rho_0)   \\
   \alpha_{(k,\ell),(k+1,\ell)}  & = &  2 \rho_{k+1}/\rho_k
                           \qquad\hbox{for } 1 \le k \le \ell-1
\end{subeqnarray}
with all other coefficients of $\alpha$ set to 0.
The state space $X$ can be visualized as a tree with root 0
and branches numbered $\ell=1,2,\ldots$ (see Figure~\ref{fig_tree}).
Now consider the vector $c = \rho_0 \delta_0$,
which puts a mass $\rho_0$ of dirt on site 0 and nothing elsewhere,
and the vector $w = \textbf{1}$.
It is easy to see that, for each $k \ge 1$,
the vector $c \alpha^k$ is supported on sites $(k,\ell)$ [$\ell \ge k$]
and takes there the values $\rho_k/2^{\ell-k+1}$.
Hence $c \alpha^k w = \rho_k$ for all $k$.
In order to construct the sequence $x_1,x_2,\ldots$,
let us first consider the following two sequences of sites:
\begin{eqnarray*}
   A_L   & \bydef  & 
      (1,2) ,\, (1,3) ,\, (2,3) ,\, (1,4) ,\, (2,4) ,\, (3,4) ,\,
      (1,L) ,\,\ldots,\, (L-1,L)
   \\[2mm]
   B_L   & \bydef  & 
      (1,1) ,\, (2,2) ,\, (3,3) ,\,\ldots,\, (L,L)
\end{eqnarray*}
so that $A_L$ (resp.\ $B_L$) sweeps, in order,
the non-summit (resp.\ summit) sites of the branches 1 through $L$.
Now let us define our sweeping process $x_1,x_2,\ldots$ as follows:
$$
   0 ,\, A_{L_1} ,\, B_1 ,\, A_{L_2} ,\, B_2 ,\, A_{L_3} ,\, B_3 ,\,\ldots
$$
where the indices $L_1 < L_2 < \ldots$ will be chosen in a moment.
After the first step ($x_1 = 0$),
we have a mass $\rho_1/2^\ell$ of dust on each site $(1,\ell)$.
After the sequence $A_{L_1}$, we end up with a mass $\rho_\ell/2$
on each site $(\ell,\ell)$ for $1 \le \ell \le L_1$.
The step $B_1$ then destroys the mass $\rho_1/2$ that sat on the site $(1,1)$.
After the sequence $A_{L_2}$, we end up with a mass $\rho_\ell/2$
on each site $(\ell,\ell)$ for $2 \le \ell \le L_2$.
The step $B_2$ then destroys the mass $\rho_2/2$ that sat on the site $(2,2)$.
And so forth.
If we choose each $L_r$ so that
\be
   \sum_{\ell=r}^{L_r} {\rho_\ell \over 2}
   \;\ge\;
   10^r + {\rho_r \over 2}   \;,
\ee
then we are assured of having a total mass at least $10^r$
on the summit sites at each step after the operation $A_{L_r}$
has been completed.
Hence $c \beta_{x_1}\cdots\beta_{x_n} w \to \infty$.
%
\qed
\end{example}

This example shows that $c \alpha^k w$ being finite for every $k$
and tending to zero as $k \to\infty$
is not sufficient to ensure that {\em every}\/ cleaning process
that visits each site infinitely many times will succeed
in removing the dirt in $l^1(w)$ sense.
Nevertheless, this hypothesis turns out to be sufficient to ensure that
there {\em exists}\/ a successful cleaning process:

\begin{theorem}[Sufficient condition for existence of cleaning]
   \label{thm.complements.2}
Let $\Lambda \subseteq X$,
and let $c \ge 0$ and $w \ge 0$ be vectors satisfying
\begin{itemize}
   \item[(a)] $c (I_\Lambda \alpha I_\Lambda)^k w \,<\, \infty$
      for all $k \ge 0$
\end{itemize}
and
\begin{itemize}
   \item[(b)] $\lim\limits_{k\to\infty} c (I_\Lambda \alpha I_\Lambda)^k w
               \,=\, 0$.
\end{itemize}
Then it is possible to find a sequence of sites $x_1,x_2,\ldots \in \Lambda$
and a sequence of numbers $\epsilon_1,\epsilon_2,\ldots \in (0,1]$
such that
\begin{equation}
  \lim\limits_{n \to\infty} c\beta_{\epsilon_1 \delta_{x_1}}
     \cdots \beta_{\epsilon_n \delta_{x_n}} I_\Lambda w
  \;=\;  0
  \;.
\end{equation}
\end{theorem}

\begin{question}
Can we always take $\epsilon_i = 1$?
\end{question}

The proof of this theorem will be based on a simple lemma:

\begin{definition}
We say that a cloud $\nu'\ge 0$ arises from a cloud $\nu\ge 0$ by a
\emph{single-marker update acting at the marker $\eta$}
if there exists a real number $\kappa \in [0, \nu(\eta)]$ such that

\begin{itemize}
\item $\nu'(\eta)=\nu(\eta)-\kappa$,
\item $\nu'(\sigma)=\nu(\sigma)+\kappa$ for each child $\sigma$ of $\eta$,
\item $\nu'(\tau)=\nu(\tau)$ for every other marker $\tau$. 
\end{itemize}
\end{definition}

\begin{lemma}[Imitation Lemma]
Suppose that the cloud $\nu'\ge 0$ arises from the cloud $\nu\ge 0$
by a single-marker update acting at the marker $\eta$.
Let $y = {\rm last}(\eta)$.
Then, for any vector $c \ge 0$,
we can find a number $\epsilon \in [0,1]$
(depending on $c$, $\nu$ and $\nu'$) such that
\begin{equation}
c T_{\nu'} \; = \; c T_\nu \beta_{\epsilon\delta_y}  \;.
  \label{eq.imitation_lemma}
\end{equation}
\end{lemma}

\proof
Let $\kappa = \nu(\eta)-\nu'(\eta)$.
If $\kappa=0$, we can take $\epsilon \bydef 0$;
so let us assume henceforth that $\nu(\eta) \ge \kappa > 0$.
We then have, for all $x\in X$,
\begin{equation}
(c T_\nu \beta_{\epsilon\delta_y})_x\;=\;
   \cases{ (cT_\nu)_x +  \epsilon (cT_\nu)_{y}\alpha_{yx}
                   & if $x \neq y$  \cr
           \noalign{\vskip 2mm}
           (1-\epsilon) (cT_\nu)_{y} +  \epsilon (cT_\nu)_{y}\alpha_{yy}
                   & if $x = y$   \cr
         }
\end{equation}
whereas
\begin{equation}
(c T_{\nu'})_x\;=\;
   \cases{ (cT_\nu)_x + \kappa (cT_\eta)_{y}\alpha_{yx} 
                   & if $x \neq y$  \cr
           \noalign{\vskip 2mm}
           (cT_\nu)_{y} - \kappa (cT_\eta)_{y}
                        + \kappa (cT_\eta)_{y} \alpha_{yy}
                   & if $x = y$   \cr
         }
\end{equation}
Therefore, to satisfy \reff{eq.imitation_lemma}, we proceed as follows:
if $(cT_\eta)_{y} = 0$ we set $\epsilon \bydef 0$;
if $(cT_\eta)_{y} > 0$, we set
\begin{equation}
   \epsilon  \;\bydef\;
   \kappa \, \frac{(cT_\eta)_{y}}{(cT_\nu)_{y}}
   \;.
\end{equation}
In the latter case we have
$(cT_\nu)_{y} \ge \nu(\eta) (cT_\eta)_{y} \ge \kappa (cT_\eta)_{y} > 0$,
so that $\epsilon \le 1$.
\qed

\proofof{Theorem~\ref{thm.complements.2}}
Since $\beta_{\epsilon_1 \delta_{x_1}} \cdots \beta_{\epsilon_n \delta_{x_n}}
       I_\Lambda  =
       I_\Lambda \beta_{\epsilon_1 \delta_{x_1}}  I_\Lambda \cdots I_\Lambda
       \beta_{\epsilon_n \delta_{x_n}} I_\Lambda$
by Lemma~\ref{lemma_betagi},
everything in both hypothesis and conclusion takes place within $\Lambda$,
so we can assume for notational simplicity that $\Lambda=X$.

By virtue of the Imitation Lemma,
it suffices to find a sequence $(\nu_n)_{n \ge 0}$ of nonnegative clouds,
with $\nu_0=\rho^0$, such that each $\nu_{n+1}$ arises from $\nu_n$
by a single-marker update, and which satisfy
\begin{equation}
\label{what_we_want}
cT_{\nu_n}w\;\mathop{\longrightarrow}_{n\to\infty}\;0  \;.
\end{equation}
We shall construct these clouds to be
$\{0,1\}$-valued and lie in $\mathcal{S}_1$.
In particular, updating the marker $\eta$ at step $n$ will mean that
$\nu_n(\eta)=1$ and $\nu_{n+1}(\eta)=0$,
and for each child $\sigma$ of $\eta$,
$\nu_n(\sigma)=0$ and $\nu_{n+1}(\sigma)=1$.

Since for each $k\ge0$ we have
\begin{equation}
   c\alpha^k w \;=\;
   \sum_{\hbox{\scriptsize level}(\eta)=k} \alpha^\eta \, \mu_{c,w}(\eta)
   \;<\; \infty  \;,
\end{equation}
for each $\epsilon>0$ we can find a finite subset $M_{\epsilon,k}$
of the set of all markers of level $k$, such that
\begin{equation}
   \sum_{\begin{scarray}
           \hbox{\scriptsize level}(\eta)=k \\
           \eta\notin M_{\epsilon,k}
         \end{scarray}}
   \alpha^\eta \, \mu_{c,w}(\eta)  \;<\;  \epsilon/2^{k+1}  \;.
\end{equation}
Furthermore, we can arrange for these subsets to be nested,
i.e.\ $M_{\epsilon,k} \subseteq M_{\epsilon',k}$
whenever $\epsilon > \epsilon'$.

We now fix a sequence $(\epsilon_j)$ of positive numbers decreasing to $0$,
and we choose a sequence of integers $N_1 < N_2 < \ldots$
such that
\be
   c\alpha^k w \le \epsilon_j \;\hbox{whenever}\;  k > N_j
   \;.
\ee
Let $S_{j,\ell}$ denote the set of all markers of level $\ell$
which are ancestors of any marker in $\bigcup_{k=0}^{N_j} M_{\epsilon_j,k}$.
By construction, all these sets are finite.
We then define the sequence $(\nu_n)_{n \ge 1}$ of clouds
by a sequence of single-marker updates starting from $\nu_0 = \rho^0$,
visiting the markers in the following order:
$$
   S_{1,0} ,\, S_{1,1} ,\,\ldots,\, S_{1,N_1} ,\,
   S_{2,0} ,\, S_{2,1} ,\,\ldots,\, S_{2,N_2} ,\;\ldots
$$
(the markers within each $S_{j,\ell}$ being updated in an arbitrary order).

Let now $j \ge 1$,
and consider any stage in this process
while we are updating the markers in
$\bigcup_{\ell=0}^{N_{j+1}} S_{j+1,\ell}$.
At such a stage,
the cloud $\nu_n$ is always supported on three kinds of markers:
\begin{itemize}
\item Markers lying on levels $k\le N_j$ but outside $M_{\epsilon_{j},k}$.
   [Since at stage $j$ we updated all the ancestors
    of markers in $\bigcup_{k=0}^{N_j} M_{\epsilon_j,k}$,
    it is impossible for subsequent operations to place mass on
    any of those markers.]
\item Markers lying on two successive levels $N$ and $N+1$
    strictly above $N_j$.
   [This happens while we are updating markers in $S_{j+1,N}$.]
\item Markers lying on levels $N_{j}< k \le N_{j+1}$
   but outside $M_{\epsilon_{j+1},k}$,
\end{itemize}
The first kind of markers contributes to $cT_{\nu_{n}}w$ at most $\epsilon_j$,
the second at most $2\epsilon_j$, and the third at most $\epsilon_{j+1}$.
This proves that (\ref{what_we_want}) holds.
\qed


The following example shows that there cannot be any converse
to Theorem~\ref{thm.complements.2} that refers only to
the behavior of the sequence $(c \alpha^k w)_{k \ge 0}$:

\begin{example}
    \label{example_thierry2}
\rm
Let $(\rho_k)_{k \ge 0}$ be {\em any}\/ sequence of positive numbers
(in particular, it can tend to $+\infty$, or oscillate, or whatever).
Then, on a countably infinite state space $X$,
we can construct a matrix $\alpha \ge 0$ and vectors $c,w \ge 0$ such that
\begin{itemize}
   \item[(a)]  $c \alpha^k w = \rho_k$ for every $k \ge 0$.
   \item[(b)]  There exists a sequence $y_1,y_2,\ldots \in X$
       for which
       $\lim\limits_{n \to \infty} c \beta_{y_1}\cdots\beta_{y_n} w = 0$.
\end{itemize}
We shall use the same state space
$X = \{0\}\cup\{(k,\ell) \colon\, 1\le k\le \ell\}$
as in Example~\ref{example_tree} (cf.\ Figure~\ref{fig_tree})
but shall make a slightly different choice of the matrix $\alpha$.

Start by choosing any sequence of positive numbers
$(\sigma_\ell)_{\ell \ge 1}$ satisfying
$\lim_{\ell \to\infty} \sigma_\ell = 0$ and
$\sum_{\ell=1}^\infty \sigma_\ell = \infty$
(e.g.\ $\sigma_\ell = 1/\ell$ will do).
Then it is not hard to see that we can find positive numbers
$(\gamma_{k\ell})_{1 \le k \le \ell}$ such that
\begin{itemize}
   \item[(i)]   $\gamma_{k\ell} \le \sigma_\ell$ for all $k,\ell$.
   \item[(ii)]  $\sum_{\ell=k}^\infty \gamma_{k\ell} = \rho_k$ for all $k$.
\end{itemize}
(For instance, for each $k$ we could choose inductively
$\gamma_{k\ell} = \min\biggl[ \sigma_\ell, \,
 \smhalf \bigl(\rho_k - \sum_{\ell'=k}^{\ell-1} \gamma_{k\ell'}\bigr)
      \biggr]$
for $\ell \ge k$.)
Now define the matrix $\alpha$ by
\begin{subeqnarray}
   \alpha_{0,(1,\ell)}           & = &  \gamma_{1\ell} / \rho_0   \\
   \alpha_{(k,\ell),(k+1,\ell)}  & = &  \gamma_{k+1,\ell} / \gamma_{k\ell}
                           \qquad\hbox{for } 1 \le k \le \ell-1
\end{subeqnarray}
with all other coefficients of $\alpha$ set to 0.
Choose once again $c = \rho_0 \delta_0$ and $w = \textbf{1}$.
Now choose the sequence $y_1,y_2,\ldots$ to be
$$
      0 ,\, (1,1) ,\, (1,2) ,\, (2,2) \,\ (1,3) ,\, (2,3) ,\, (3,3) ,\,\ldots
$$
It is easy to see that,
at any stage after the site $(\ell,\ell)$ has been swept
and before the site $(\ell+1,\ell+1)$ has been swept,
the total quantity of dirt $c \beta_{y_1}\cdots\beta_{y_n} w$
does not exceed $\sigma_{\ell+1} + \sum_{\ell'=\ell+2}^\infty \gamma_{1,\ell'}$;
and this tends to zero as $\ell\to\infty$.

It is worth remarking that if (and only if)
$\sum_{k=0}^\infty \rho_k = \infty$,
then it is possible to choose first the $(\sigma_\ell)$
and then the $(\gamma_{k\ell})$ so that,
in addition to properties (i) and (ii), we have
\begin{itemize}
   \item[(iii)]  $\sum_{\ell=1}^\infty  \gamma_{\ell\ell} = \infty$.
\end{itemize}
Indeed, let $(\sigma_k)_{k \ge 1}$ be {\em any}\/ sequence of positive numbers
satisfying $\lim_{k \to\infty} \sigma_k = 0$ and
$\sum_{k=1}^\infty \min(\rho_k,\sigma_k) = \infty$.\footnote{
   Given any sequence $(\rho_k)_{k \ge 1}$ of positive numbers
   satisfying $\sum_{k=1}^\infty \rho_k = \infty$,
   we can always find a sequence $(\sigma_k)_{k \ge 1}$ such that
   $0 < \sigma_k \le \rho_k$ for all $k$,
   $\lim_{k \to\infty} \sigma_k = 0$
   and $\sum_{k=1}^\infty \sigma_k = \infty$.
   To see this, note first that $\rho'_k := \min(\rho_k,1)$ also satisfies
   $\sum_{k=1}^\infty \rho'_k = \infty$.
   Now let
   $ \sigma_k  \,:=\,
     \rho'_k  \big/  \big( 1 + \sum_{j=1}^k \rho'_j \big) ^{1/2}$.
   Clearly $\lim_{k \to\infty} \sigma_k = 0$, and
   $ \sum_{k=1}^N \sigma_k
      \,\ge\,
      \big( \sum_{k=1}^N \rho'_k \big)
         \big/
         \big( 1 + \sum_{k=1}^N \rho'_k \big) ^{1/2}
      \,\to\,
      \infty
   $
   as $N \to\infty$.
}
Then the preceding construction yields
$\gamma_{kk} = \min(\sigma_k, \smhalf \rho_k)$,
so that $\sum_{k=1}^\infty  \gamma_{kk} = \infty$.
In this way we can find a matrix $\alpha$ that serves simultaneously
for Example~\ref{example_tree} and the present example:
that is, there exists {\em both}\/ a sequence $x_1,x_2,\ldots \in X$
in which each element of $X$ occurs infinitely often
and such that
$\lim\limits_{n \to \infty} c \beta_{x_1}\cdots\beta_{x_n} w = \infty$,
{\em and}\/ a sequence $y_1,y_2,\ldots \in X$ such that
$\lim\limits_{n \to \infty} c \beta_{y_1}\cdots\beta_{y_n} w = 0$.
\qed
\end{example}

%
%
%

\section{Converse results}   \label{sec.converse}


{\em In this section we do {\bf not} assume the Fundamental Hypothesis.}\/
Rather, our goal is to study what happens in case
the Fundamental Hypothesis fails.
As before, we are entitled to use here those results
of Section~\ref{sec4} that refer only to clouds,
but must avoid using those results that refer to operators.


\subsection{General result}   \label{sec.converse.1}

Our main result is the following:

\begin{theorem}
   \label{thm.converse.2.NEW}
Let $X$ be a finite or countably infinite set,
let $\Lambda \subseteq X$,
and let $c \ge 0$ and $w \ge 0$ be vectors
that are strictly positive on $\Lambda$.
Consider the following conditions on a matrix $\alpha$:
\begin{itemize}
  \item[(a)]  $\sum_{k=0}^\infty c (I_\Lambda \alpha I_\Lambda)^k w < \infty$.
  \item[(b)]  For all $h \colon\, X \to [0,1]$ with $\supp h = \Lambda$
      such that $h \ge \epsilon \chi_\Lambda$ for some $\epsilon > 0$,
      we have
      $\sum_{k=0}^\infty c (I_\Lambda \beta_h I_\Lambda)^k w < \infty$.
  \item[(b${}'$)]  There exists $h \colon\, X \to [0,1]$
      with $\supp h = \Lambda$ such that
      $\sum_{k=0}^\infty c (I_\Lambda \beta_h I_\Lambda)^k w < \infty$.
   \item[(c)]  For every finite sequence $f_1,\ldots,f_m$
      of functions $X \to [0,1]$ with $\supp(f_i) \subseteq \Lambda$
      such that $\sum_i f_i \ge \epsilon \chi_\Lambda$
      for some $\epsilon > 0$, we have
      $\sum_{k=0}^\infty c (I_\Lambda \beta_{f_1} \cdots
                                      \beta_{f_m} I_\Lambda)^k w < \infty$.
   \item[(c${}'$)]  There exists a finite sequence $f_1,\ldots,f_m$
      of functions $X \to [0,1]$ with $\supp(f_i) \subseteq \Lambda$
      such that
      $\sum_{k=0}^\infty c (I_\Lambda \beta_{f_1} \cdots
                                      \beta_{f_m} I_\Lambda)^k w < \infty$.
   \item[(d)]  For every cloud $\nu \ge 0$
      that is carried by $\Lambda$
      and satisfies $\triplenorm \nu \triplenorm \le 1$
      and $\nu((x)) \le 1-\epsilon$ for every $x \in \Lambda$
      (for some $\epsilon > 0$),
      we have $\sum_{k=0}^\infty c (I_\Lambda T_\nu I_\Lambda)^k w < \infty$.
   \item[(d${}'$)]  There exists a cloud $\nu \in \scrs_1$
      that is carried by $\Lambda$
      and supported on finitely many levels,
      such that $\sum_{k=0}^\infty c (I_\Lambda T_\nu I_\Lambda)^k w < \infty$.
\end{itemize}
Then (a) $\iff$ (b) $\iff$ (b${}'$) $\iff$ (c) $\iff$ (d)
  $\implies$ (c${}'$) $\implies$ (d${}'$);
and for matrices $\alpha$ satisfying the additional hypothesis
\be
  \hbox{There exists a constant $C < \infty$ such that
        $I_\Lambda \alpha I_\Lambda w \le C w$}
 \label{hyp.converse.2}
\ee
all seven conditions are equivalent.
\end{theorem}

\proof
(d) $\implies$ (c) $\implies$ (b)
  $\implies$ (b${}'$) $\implies$ (c${}'$) $\implies$ (d${}'$)
is trivial.
So it suffices to prove (a) $\implies$ (d)
and (b${}'$) $\implies$ (a),
as well as to prove (d${}'$) $\implies$ (a) under the hypothesis
\reff{hyp.converse.2}.

(a) $\implies$ (d):
Since $\nu$ is carried by $\Lambda$,
for $k \ge 1$, $(I _\Lambda * \nu * I_\Lambda)^k$
charges only markers with all their entries in $\Lambda$.
Lemma~\ref{lemma.12Apr07.cloud}(b)
then implies that
\begin{equation}
   \sum_{k=0}^\infty (I_\Lambda * \nu * I_\Lambda)^k
   \;\le\;
   \frac{1}{\epsilon} \, {\bf 1}_\Lambda \,+\, I_{\Lambda^c}
   \;\le\;
   \frac{1}{\epsilon} \:
        \sum_{k=0}^\infty (I_\Lambda * \rho^1 * I_\Lambda)^k
 \label{eq.proof.thm.converse.2.NEW.a-e}
\end{equation}
(here $\epsilon \le 1$).
Passing from clouds to operators,
we have
\be
   \sum_{k=0}^\infty c (I_\Lambda T_\nu I_\Lambda)^k w
   \;\le\;
   \frac{1}{\epsilon} \:
   \sum_{k=0}^\infty c (I_\Lambda \alpha I_\Lambda)^k w
   \;<\;
   \infty
   \;.
\ee

(b${}'$) $\implies$ (a) and  (d${}'$) $\implies$ (a):
Let us start with hypothesis (d${}'$),
and observe that hypothesis (b${}'$) is merely the special case
in which $\nu = \beta_h$.

Note first that we must have
$\nu((x)) < 1$ for every marker $(x)$ of level 0 with $x \in \Lambda$,
since otherwise $c (I_\Lambda T_\nu I_\Lambda)^k w$
could not tend to zero as $k\to\infty$.
(Here we have used the strict positivity of $c$ and $w$ on $\Lambda$.)

Now define $\mu \bydef {\bf 1} - \widetilde{\nu}$.
By Lemma~\ref{lemma.16Apr07}(c),
$\mu$ takes values in $[0,1]$ and is
supported on markers having all their entries in $\Lambda$.
By Lemma~\ref{lemma.16Apr07}(d),
if $\nu$ is supported on levels $\le K$,
then $\mu$ is supported on levels $\le K-1$.
Putting these facts together, we obtain
\be
   0  \;\le\;  \mu    \;\le\;
   I_\Lambda * \left( \sum_{k=0}^{K-1} (I_\Lambda * \rho^1 * I_\Lambda)^k
               \right) *  I_\Lambda
   \;.
 \label{eq.lemma.11Apr07.a}
\ee
By Lemma~\ref{lemma.16Apr07}(c) we then have
\begin{subeqnarray}
    {\bf 1}_\Lambda
    & = &
    \left( \sum_{n=0}^{\infty} (I_\Lambda * \nu * I_\Lambda)^n \right)
    * \mu
            \\[1mm]
    & \le &
    \left( \sum_{n=0}^{\infty} (I_\Lambda * \nu * I_\Lambda)^n \right)
    *
    \left( \sum_{k=0}^{K-1} (I_\Lambda * \rho^1 * I_\Lambda)^k \right)
    \;.
\end{subeqnarray}
Passing now to operators, we have
\be
   \sum_{k=0}^\infty (I_\Lambda \alpha I_\Lambda)^k
   \;\le\;
   \left( \sum_{n=0}^{\infty} (I_\Lambda T_\nu I_\Lambda)^n \right)
   *
   \left( \sum_{k=0}^{K-1} (I_\Lambda \alpha I_\Lambda)^k \right)
    \;.
 \label{eq.thm.converse.2.NEW.proof.bimpliesa}
\ee
Now sandwich this between nonnegative vectors $c$ and $w$.
If $\nu = \beta_h$ [case (b${}'$)], then $K=1$,
so that the second large parenthesis on the right-hand side of
\reff{eq.thm.converse.2.NEW.proof.bimpliesa} is the identity operator.
In case (d${}'$), we use hypothesis \reff{hyp.converse.2}.
Either way, we find that
\be
   \sum_{k=0}^\infty c (I_\Lambda \alpha I_\Lambda)^k w
   \;\le\;
   C' \sum_{n=0}^{\infty} c (I_\Lambda T_\nu I_\Lambda)^n w
\ee
for a finite constant $C'$. This completes the proof.
\qed

Please note that hypothesis \reff{hyp.converse.2} is automatic
whenever $\Lambda$ is a {\em finite}\/ set.
On the other hand,
Examples~\ref{example_thm.converse.2.NEW_1} and
\ref{example_thm.converse.2.NEW_1.bis}
show that, when $\Lambda$ is infinite, hypothesis \reff{hyp.converse.2}
cannot be dispensed with in proving that (c${}'$) $\implies$ (a)
[or, {\em a fortiori}\/, (d${}'$) $\implies$ (a)].

The following example shows that the hypothesis
that $\nu$ is supported on finitely many levels
cannot be dispensed with in proving that (d${}'$) $\implies$ (a),
even in the presence of the Fundamental Hypothesis
[which is much stronger than \reff{hyp.converse.2}]:

\begin{example}
    \label{example_thm.converse.2.NEW_2}
\rm
Let $\Lambda = X = \{1,2,3,\ldots\}$.
Let $\alpha$ be right shift (acting on dirt vectors),
i.e.\ $\alpha_{i,i+1} = 1$ for all $i \ge 1$
and all other matrix elements of $\alpha$ are 0.
Set $c = {\bf 1}$ and let $w_j = 1/j^2$.
Note that $\alpha w \le w$.  Then
\be
   (\alpha^k w)_i  \;=\;  w_{i+k}  \;=\; {1 \over (i+k)^2}
   \;,
\ee
so that
\be
   \sum_{k=0}^\infty (\alpha^k w)_i  \;=\;
   \sum_{k=0}^\infty {1 \over (i+k)^2} \;\sim\;
    {1 \over i} \quad \hbox{as } i \to \infty
\ee
and hence $\sum_{k=0}^\infty c \alpha^k w = +\infty$.
On the other hand, let us take $\nu$ to be the cloud
\be
   \nu(\eta)  \;=\;
   \cases{ 1   & if $\hbox{level}(\eta) = \hbox{first}(\eta)$  \cr
           0   & otherwise                                     \cr
         }
\ee
(In other words, for dirt starting at site $i$,
 $\nu$ sends it upwards $i$ levels.)
For this choice of $\alpha$ we have
\be
   (T_\nu^k w)_i  \;=\;  w_{2^k i}  \;=\; {4^{-k} \over i^2}
   \;,
\ee
so that $\sum_{k=0}^\infty c T_\nu^k w < \infty$.
\qed
\end{example}


In view of Theorem~\ref{thm.main.operators},
which guarantees cleaning whenever $h$ is strictly positive on $\Lambda$ ---
without any need for uniformity ---
one might be tempted to remove the uniformity hypothesis
in condition (b), i.e. to replace (b) by
\begin{itemize}
   \item[(b${}^*$)] For all $h \colon\, X \to [0,1]$ with $\supp h = \Lambda$,
      we have
      $\sum_{k=0}^\infty c (I_\Lambda \beta_h I_\Lambda)^k w < \infty$.
\end{itemize}
But neither (a) nor any conceivable stronger hypothesis
can possibly imply this, as Example~\ref{example.thm.main.operators} shows.
%
%
The upshot is that, without a uniformity hypothesis on $h$,
one can conclude in some cases that
$c (I_\Lambda \beta_h I_\Lambda)^k w$ {\em tends to zero as $k\to\infty$}\/
[cf.\ Theorem~\ref{thm.main.operators}
 and Corollary~\ref{cor.complements.1}],
but one cannot conclude anything about the {\em rate}\/ of convergence
--- in particular, one cannot conclude that this sequence is
{\em summable}\/ in $k$.


\subsection{Finite $\Lambda$}

If $\Lambda$ is a {\em finite}\/ subset of $X$
(and in particular if $X$ is finite),
then all choices of strictly positive vectors $c$ and $w$ are equivalent,
and Theorem~\ref{thm.converse.2.NEW} can be rephrased in a simpler form,
in terms of the spectral radii of the various matrices.
The main idea is that if the spectral radius of
$I_\Lambda \alpha I_\Lambda$ is $\ge 1$,
then there is no way to clean the set $\Lambda$ completely.

Let us recall that the spectral radius
(= largest absolute value of an eigenvalue)
of a finite matrix $A$ satisfies
$\spr(A) = \lim\limits_{n\to\infty} \|A^n\|^{1/n}$;
in particular, the latter limit is independent of the choice of norm.
It is well known
(and easy to prove using $\|A^{m+n}\| \le \|A^m\| \, \|A^n\|$) that
$\spr(A) < 1$
$\iff$
$\lim\limits_{n\to\infty} \|A^n\| = 0$
$\iff$
$\sum\limits_{n=0}^\infty \|A^n\| < \infty$.
Finally, on a finite set the convergence of sequences of vectors or matrices
can be understood either pointwise or in norm;
the two notions are equivalent.

Since the uniformity conditions on $h$, etc.\ in
Theorem~\ref{thm.converse.2.NEW} are trivially satisfied
when $\Lambda$ is finite,
we can state the following immediate corollary
of Theorem~\ref{thm.converse.2.NEW}:

%

\begin{corollary}
   \label{cor.converse.finite.lambda}
Let $X$ be a finite or countably infinite set,
and let $\Lambda \subseteq X$ be a {\em finite}\/ subset.
Then the following conditions on a matrix $\alpha$ are equivalent:
\begin{itemize}
  \item[(a)]  $\spr(I_\Lambda \alpha I_\Lambda) < 1$.
  \item[(b)]  For all $h \colon\, X \to [0,1]$ with $\supp h = \Lambda$,
      we have $\spr(I_\Lambda \beta_h I_\Lambda) < 1$.
  \item[(b${}'$)]  There exists $h \colon\, X \to [0,1]$
      with $\supp h \subseteq \Lambda$ such that
      $\spr(I_\Lambda \beta_h I_\Lambda) < 1$.
   \item[(c)]  For every finite sequence $f_1,\ldots,f_m$
      of functions $X \to [0,1]$ with $\supp(\sum_{i=1}^m f_i) = \Lambda$,
      we have
      $\spr(I_\Lambda \beta_{f_1} \cdots \beta_{f_m} I_\Lambda) < 1$.
   \item[(c${}'$)]  There exists a finite sequence $f_1,\ldots,f_m$
      of functions $X \to [0,1]$ with $\supp(f_i) \subseteq \Lambda$
      such that
      $\spr(I_\Lambda \beta_{f_1} \cdots \beta_{f_m} I_\Lambda) < 1$.
   \item[(d)]  For every cloud $\nu \ge 0$
      that is carried by $\Lambda$
      and satisfies $\triplenorm \nu \triplenorm \le 1$
      and $\nu((x)) < 1$ for every $x \in \Lambda$,
      we have $\spr(I_\Lambda T_\nu I_\Lambda) < 1$.
   \item[(d${}'$)]  There exists a cloud $\nu \in \scrs_1$
      that is carried by $\Lambda$
      and supported on finitely many levels,
      such that $\spr(I_\Lambda T_\nu I_\Lambda) < 1$.
   \item[(e)]  For every infinite sequence $f_1,f_2,\ldots$
         of functions $X \to [0,1]$ with $\supp(f_i) \subseteq \Lambda$
         such that $\sum_i f_i = \infty$ everywhere on $\Lambda$, we have
         $\lim\limits_{n\to\infty}
          I_\Lambda \beta_{f_1} \cdots \beta_{f_n} I_\Lambda = 0$.
   \item[(f)]  For every infinite sequence $\nu_1,\nu_2,\ldots$
         of nonnegative clouds carried by $\Lambda$ and satisfying
         $\triplenorm \nu_i \triplenorm \le 1$ 
         and $\sum_i [1-\nu_i((x))] = \infty$ for all $x \in \Lambda$, we have
         $\lim\limits_{n\to\infty}
          I_\Lambda T_{\nu_1} \cdots T_{\nu_n} I_\Lambda = 0$.
\end{itemize}
\end{corollary}

\proof
The equivalence of (a), (b), (b${}'$), (c), (c${}'$), (d) and (d${}'$)
is an immediate consequence of Theorem~\ref{thm.converse.2.NEW}
specialized to $\Lambda$ finite.
The implication (a) $\implies$ (f) follows from
Corollary~\ref{cor.complements.1},
and (f) $\implies$ (e) is trivial.
Finally, (e) $\implies$ (c) [or (f) $\implies$ (d)] is easy:
just consider the sequence $f_1,\ldots,f_m$ repeated infinitely many times
and use Lemma~\ref{lemma_betagi.clouds}.
\qed

The following example shows that, when $X$ is infinite,
the behavior of the matrix $\alpha$
is not completely controlled by that of the matrices
$I_\Lambda\alpha I_\Lambda$ for all finite $\Lambda$:

\begin{example}
   \label{example_sec6_FHfails}
\rm
We shall exhibit a nonnegative matrix $\alpha$ satisfying
$\spr(I_\Lambda\alpha I_\Lambda)<1$
for all finite $\Lambda \subset X$
--- so that in particular the Fundamental Hypothesis holds
for all the matrices $I_\Lambda\alpha I_\Lambda$ ---
but for which the Fundamental Hypothesis does not hold,
i.e.\  we cannot find a vector $w>0$ satisfying $\alpha w\le w$.

Take $X = \{1,2,\ldots\}$ and $0 \le \epsilon<1$,
and set $\alpha_{1,1} \bydef \epsilon$,
$\alpha_{1,j} \bydef 1$ for all $j\ge2$,
$\alpha_{i+1,i} \bydef 2$ for all $i\ge2$,
and $\alpha_{i,j} \bydef 0$ elsewhere:
\begin{equation}
\alpha   \;=\;
\left(\!
\begin{array}{ccccc}
\epsilon & 1 & 1 & 1 & \cdots \\ 
0 & 0 & 0 & 0 & \cdots \\ 
0 & 2 & 0 & 0 & \cdots \\ 
0 & 0 & 2 & 0 & \cdots \\ 
\;\vdots\; & \ddots & \ddots & \ddots & \ddots
\end{array}
\!\right)
\end{equation}
If we could find a $w>0$ such that $\alpha w\le w$,
we would have $w_{i+1}\ge 2w_i$ for all $i\ge 2$.
Since $w_2>0$, this would imply $w_i\to+\infty$ as $i \to\infty$.
But $\alpha w\le w$ implies also that $w_1\ge \sum_{i=2}^\infty w_i$,
which is impossible. 

Now let us show why $\spr(I_\Lambda\alpha I_\Lambda)<1$
for all finite $\Lambda \subset X$.
Set $m \bydef \max \Lambda$.
Then, for any vector $v$ and any $j\ge m-1$,
we have $(I_\Lambda\alpha I_\Lambda)^j v \in \mathop{\mbox{span}}(e_1)$
where $e_1 \bydef (1,0,0,\ldots)$.
Since $\alpha e_1=\epsilon e_1$, it follows that
\be
   \spr(I_\Lambda \alpha I_\Lambda)
   \;=\;
   \cases{  \epsilon   & if $1 \in \Lambda$  \cr
            0          & if $1 \notin \Lambda$ \cr
         }
\ee
In particular, by taking $\epsilon=0$ we can even arrange to have
$\spr(I_\Lambda \alpha I_\Lambda) = 0$
[i.e., $I_\Lambda \alpha I_\Lambda$ is nilpotent]
for every finite $\Lambda$.
\qed
\end{example}

\section*{Acknowledgments}

We wish to thank Claude Dellacherie for valuable conversations
on probabilistic potential theory;
Uriel Rothblum for a valuable correspondence
concerning Perron--Frobenius theorems;
and Victor Blondel, Kousha Etessami and Alex Scott
for valuable correspondence and/or conversations concerning
computational complexity and decidability.
Last but not least, we wish to thank Marina Papa Sokal
for emphasizing to us the importance of the {\em experimental}\/ study
of the floor-cleaning problem.

This research was supported in part by
U.S.\ National Science Foundation grants PHY--0099393 and PHY--0424082.
One of us (A.D.S.)\ also wishes to thank the Universit\'e de Rouen
for kind hospitality and financial support during several visits.

\clearpage


\begin{thebibliography}{99}

\bibitem{Aizenman_unpub}  M. Aizenman and E.H. Lieb, Dobrushin revisited,
   unpublished manuscript (October 1984).

\bibitem{Berman_79}  A. Berman and R.J. Plemmons,
   {\em Nonnegative Matrices in the Mathematical Sciences}\/
   (Academic Press, New York, 1979).

\bibitem{Blondel_03}  V.D. Blondel and V. Canterini,
   Undecidable problems for probabilistic automata of fixed dimension,
   Theory Comput. Systems {\bf 36}, 231--245 (2003).

\bibitem{Blondel_97}  V.D. Blondel and J.N. Tsitsiklis,
   When is a pair of matrices mortal?,
   Inform. Process. Lett. {\bf 63}, 283--286 (1997).

\bibitem{Blondel_00}  V.D. Blondel and J.N. Tsitsiklis,
   A survey of computational complexity results in systems and control,
   Automatica {\bf 36}, 1249--1274 (2000).

\bibitem{Blondel_05}  V.D. Blondel and J.N. Tsitsiklis,
   Stable matrix in a semigroup generated by
   nonnegative matrices is undecidable,
   unpublished manuscript (March 2005).

\bibitem{Dobrushin_68}  R.L. Dobrushin,
The description of a random field by means of conditional probabilities
and conditions of its regularity,
Theory Probab. Appl. {\bf 13}, 197--224 (1968)
[= Teor. Verojatnost. i Primenen. {\bf 13}, 201--229 (1968)].

\bibitem{Dobrushin_70}  R.L. Dobrushin,
Prescribing a system of random variables by conditional distributions,
Theory Probab. Appl. {\bf 15}, 458--486 (1970)
[= Teor. Verojatnost. i Primenen. {\bf 15}, 469--497 (1970)].


\bibitem{Dobrushin-Shlosman}  R.L. Dobrushin and S.B. Shlosman,
Constructive criterion for the uniqueness of Gibbs field,
in {\em Statistical Physics and Dynamical Systems}\/ (K\"oszeg, 1984),
edited by J. Fritz, A. Jaffe and D. Sz\'asz
(Birkh\"auser, Boston MA, 1985), pp.~347--370.

\bibitem{Follmer_82}  H. F\"ollmer, A covariance estimate for Gibbs measures,
   J. Funct. Anal. {\bf 46}, 387--395 (1982).

\bibitem{Georgii_88} H.-O. Georgii,
    {\em Gibbs Measures and Phase Transitions}\/
    (de Gruyter, Berlin--New York, 1988), Section 8.1.

\bibitem{Kemeny_66}  J.G. Kemeny, J.L. Snell and A.W. Knapp,
   {\em Denumerable Markov Chains}\/ (Van Nostrand, Princeton, 1966),
   especially Chapters 7--9.

\bibitem{Lanford_73}
O.E. Lanford III, Entropy and equilibrium states in classical statistical
mechanics, in {\em Statistical Mechanics and Mathematical Problems}\/
(Lecture Notes in Physics \#20), edited by A. Lenard
(Springer-Verlag, Berlin, 1973), pp.~1--113.

\bibitem{Nummelin_84}  E. Nummelin, {\em General Irreducible Markov Chains
   and Non-Negative Operators}\/
   (Cambridge University Press, Cambridge, 1984),
   especially Chapter 3.


\bibitem{Revuz_75}  D. Revuz,  {\em Markov Chains}\/
   (North-Holland, Amsterdam, 1975),
   especially Chapters 1 and 2.

\bibitem{Rothblum_private}  U. Rothblum, private communication (April 2004).

\bibitem{Seneta_81}  E. Seneta, {\em Non-Negative Matrices and
   Markov Chains}\/, 2nd ed.
   (Springer-Verlag, New York--Heidelberg--Berlin, 1981).

\bibitem{Simon_93}
B. Simon, {\em The Statistical Mechanics of Lattice Gases}\/, vol.~I
(Princeton University Press, Princeton, 1993), Sections V.1 and V.3.

\bibitem{Spitzer_76}  F.L. Spitzer, {\em Principles of Random Walk}\/,
   $2^{nd}$ ed. (Springer-Verlag, New York, 1976).

\bibitem{Stanley_99}  R.P. Stanley, {\em Enumerative Combinatorics}\/,
       vols.~1 and 2 (Cambridge University Press, Cambridge--New York, 1999).

\bibitem{Vasershtein_69}  L.N. Vasershtein,
Markov processes over denumerable products of spaces describing
large system of automata,
Problems Inform. Transmission {\bf 5}, no.~3, 47--52 (1969)
[= Problemy Pereda\v{c}i Informacii {\bf 5}, no.~3, 64--72 (1969)].

\bibitem{Weitz_05}  D. Weitz,
Combinatorial criteria for uniqueness of Gibbs measures,
Random Struct. Algorithms {\bf 27}, 445--475 (2005).


\end{thebibliography}
\end{document}